\documentclass{gtpart}

\usepackage[all]{xy}
\usepackage{graphicx}
\usepackage{hyperref}
\usepackage{rotating}
\usepackage{pinlabel}

\author{Tim Perutz}
\givenname{Tim}
\surname{Perutz}
\address{DPMMS, Centre for Mathematical Sciences, University of Cambridge, Wilberforce Road, Cambridge CB3 0WB, UK}
\email{tim.perutz@cantab.net}

\title{Lagrangian matching invariants for fibred four--manifolds: I}



\arxivreference{math.SG/0606061}
\arxivpassword{}

\volumenumber{11}
\issuenumber{}
\publicationyear{2007}
\papernumber{}
\startpage{}
\endpage{}
\doi{}
\MR{}
\Zbl{}
\received{7 June 2006}
\revised{}
\accepted{27 March 2007}
\published{}
\publishedonline{}
\proposed{Peter Ozsv\'ath}
\seconded{Simon Donaldson, Ron Fintushel}
\corresponding{}
\editor{}
\version{}

\theoremstyle{plain}
\newtheorem{Thm}{Theorem}[section]
\newtheorem{Lem}[Thm]{Lemma}
\newtheorem{Prop}[Thm]{Proposition}
\newtheorem{Cor}[Thm]{Corollary}

\newtheorem{Conj}[Thm]{Conjecture}

\newtheorem{Rk}[Thm]{Remark}
\newtheorem{Ex}[Thm]{Example}

\newtheorem{MonThm}[Thm]{Monodromy theorem}
\newtheorem*{CorrThm}{Theorem A}
\newtheorem*{MatchThm}{Theorem B}
\newtheorem*{TheoremC}{Construction Theorem C}
\newtheorem*{TheoremD}{Theorem D}
\newtheorem*{TheoremE}{Theorem E}
\newtheorem*{TheoremF}{Theorem F}

\theoremstyle{definition}
\newtheorem{Defn}[Thm]{Definition}

\newenvironment{pf}{ \begin{proof} }{ \end{proof} }

\DeclareMathAlphabet\EuScript{U}{eus}{m}{n}
\SetMathAlphabet\EuScript{bold}{U}{eus}{b}{n}

\DeclareFontFamily{U}{eus}{\skewchar\font'60}%
\DeclareFontShape{U}{eus}{m}{n}{<-6>eusm5<6-8>eusm7<8->eusm10}{}%
\DeclareFontShape{U}{eus}{b}{n}{<-6>eusb5<6-8>eusb7<8->eusb10}{}%

\newcommand{\PS}{\mathbb{P}}
\newcommand{\T}{\mathbb{T}}

\DeclareMathOperator*{\SO}{\mathrm{SO}}
\DeclareMathOperator*{\SU}{\mathrm{SU}}
\DeclareMathOperator*{\U}{\mathrm{U}}
\DeclareMathOperator*{\GL}{\mathrm{GL}}
\DeclareMathOperator*{\Sp}{\mathrm{Sp}}

\DeclareMathOperator{\im}{im}

\DeclareMathOperator{\Hom}{\mathrm{Hom}}
\DeclareMathOperator{\coker}{\mathrm{coker}}

\DeclareMathOperator{\sym}{Sym}
\DeclareMathOperator{\Pic}{Pic}

\DeclareMathOperator{\Ppic}{\EuScript{P}}
\DeclareMathOperator{\hilb}{Hilb}

\DeclareMathOperator{\imag}{Im}
\DeclareMathOperator{\real}{Re}

\DeclareMathOperator{\aut}{Aut}
\DeclareMathOperator{\diff}{Diff}
\DeclareMathOperator{\ham}{Ham}

\DeclareMathOperator{\supp}{Supp}
\DeclareMathOperator{\spinc}{\mathrm{Spin}^c}

\DeclareMathOperator{\hol}{hol}
\DeclareMathOperator{\torus}{T}

\DeclareMathOperator{\Div}{div}

\newcommand{\ii}{\mathrm{i}}
\newcommand{\id}{\mathrm{id}}
\newcommand{\Tv}{T^{\mathrm{v}}}
\newcommand{\Th}{T^{\mathrm{h}}}

\newcommand{\univ}{\mathrm{univ}}
\newcommand{\crit}{\mathrm{crit}}

\def\theheight{4.2ex}
\def\HMto{\smash{\mskip2mu\raise\theheight\rlap{%
      \begin{turn}{180}
          {$\widehat{\phantom{\mathit{HM}}}$}
       \end{turn}} \mskip-2mu
            \mathit{HM}
}{\vphantom{\widehat{\mathit{HM}}}}{}}

\begin{document}

\begin{abstract}
In a pair of papers, we construct invariants for smooth four--manifolds equipped with `broken fibrations'---the singular Lefschetz fibrations of  Auroux, Donaldson and Katzarkov---generalising the Donaldson--Smith invariants for Lefschetz fibrations.

The `Lagrangian matching invariants' are designed to be comparable with the Seiberg--Witten invariants of the underlying four--manifold;
formal properties and first computations support the conjecture that equality holds. They fit into a field theory which assigns Floer homology groups to three--manifolds fibred over $S^1$.

The invariants are derived from moduli spaces of pseudo-holomorphic sections of relative Hilbert schemes of points on the fibres, subject to Lagrangian boundary conditions. Part I---the present paper---is devoted to the symplectic geometry of these Lagrangians.
\end{abstract}

\maketitle

\setcounter{page}{1}

\section{Introduction}
The Seiberg--Witten invariants of a symplectic four--manifold can be calculated, according to Taubes' famous theorem \cite{Ta3},
as Gromov invariants enumerating embedded pseudo-holomorphic curves and their unramified coverings. In the presence of a \emph{symplectic Lefschetz fibration}, the Donaldson--Smith invariants \cite{DS,Smi} mediate
between the gauge-theoretic and symplectic viewpoints. They are counts of pseudo-holomorphic multisections of the fibration, within a chosen homology class---more properly, of pseudo-holomorphic sections of an associated family of symmetric products of the non-singular fibres, appropriately compactified over the singular fibres. The motivating observation is that an embedded
pseudo-holomorphic curve in the four--manifold, not having a fibre as a component, will have positive intersections with the
fibres and so define a section of a family of symmetric products.

The equality of the Donaldson--Smith and Gromov invariants, for fibrations of high degree, has been proved by Usher \cite{Ush}.
In the other direction, the link between symmetric products and gauge theory arises from the fact that the dimensionally reduced
Seiberg--Witten equations on a surface are the abelian vortex equations. The moduli space of solutions to these equations is a
symmetric product of the surface; in the `adiabatic limit' as the base of the fibration is expanded, solutions to the Seiberg--Witten
equations approximate pseudo-holomorphic curves in the family of vortex moduli spaces \cite{Sa2, Ta3}.

In this paper and its sequel we extend the Donaldson--Smith construction in two directions. First, we generalise it to
\emph{singular Lefschetz fibrations} in the sense of Auroux, Donaldson and Katzarkov \cite{ADK},
objects which we shall refer to as `broken fibrations'. These are available on any four--manifold with $b_2^+>0$---the point being that these are the manifolds which admit \emph{near-symplectic forms}. We construct an invariant of broken fibrations, the {\bf Lagrangian matching invariant}, which can be compared to the Seiberg--Witten invariant of the underlying manifold.
We conjecture that equality holds. It still counts pseudo-holomorphic sections of the associated families of symmetric products, but these are now subject to certain Lagrangian boundary conditions. Much of this paper is concerned with the construction of these Lagrangians, and with teasing out their properties.

Second, we show that the Lagrangian matching invariant arises from a field theory: a $(1+1)$--dimensional TQFT coupled to
singular surface-bundles. To a three--manifold $Y$ with a fibre bundle $\pi\colon Y\to S^1$ and a $\spinc$--structure $\mathfrak{t}$ (subject to certain restrictions) we assign a symplectic Floer homology group $HF_*(Y,\pi, \mathfrak{t})$, and when $Y$ is the boundary of a broken fibration there is a relative Lagrangian matching invariant in $HF_*(Y,\pi,\mathfrak{t})$.

Another use of our Lagrangian boundary conditions is to define a Floer homology group $HF_*(Y,\pi,\mathfrak{t})$ when $\pi$ is an $S^1$-valued Morse function without local extrema. Its Euler characteristic is the Turaev torsion. We shall explain the construction in a separate paper.

\subsection{Relation to near-symplectic geometry}
The construction of Lagrangian matching invariants was guided by Taubes' programme \cite{Ta1} to obtain
the Seiberg--Witten invariants of a near-symplectic four--manifold as generalised Gromov invariants. However, making a rigorous comparison presents considerable challenges (besides the matter of precisely defining the generalised Gromov invariants); there is no `tautological correspondence' in the sense of Usher \cite{Ush}.

It may be worth emphasising that, in contrast to Taubes' framework, the technical difficulties in the pseudo-holomorphic theory underlying our invariants from our moduli spaces are rather mild, at least if one does not aim for the greatest conceivable generality. Rather, the difficulty in formulating these invariants was in finding good moduli spaces to consider.

It seems to be typical of theories based on symmetric products that there are technical gains in the pseudo-holomorphic theory, and in manifest functoriality, but a loss in explicitness which can make computations hard. (A strength of Heegaard Floer homology is that it strikes an effective balance between these aspects.)

\subsection{Relation to Seiberg--Witten invariants}
The conjectural equality of Lagrangian matching and Seiberg--Witten invariants has the flavour of an `Atiyah--Floer' conjecture. For a start, the Lagrangian boundary conditions, which we construct by direct symplectic means, may well have a gauge-theoretic interpretation arising from the Seiberg--Witten equations on a three--manifold with boundary (see Remark \ref{gauge lags}).  It would be more interesting, though, to find a \emph{formal} (TQFT) reason for equality, in the vein of Donaldson's argument in \cite{Don}.

\subsection{Relation to Heegaard Floer homology}
There is a rather direct link with Heegaard Floer homology, and also with Yi-Jen Lee's programme to relate it to monopole Floer homology \cite{Lee}, which will be developed in a future article.

As mentioned above, the Lagrangian boundary conditions studied in the present paper, and used to define the Lagrangian matching invariants in its sequel, can also be employed to define symplectic Floer homology groups for a three--manifold with an $S^1$-valued Morse function all of whose critical points have indefinite indices ($1$ or $2$). Now, given  a \emph{self--indexing} Morse function $f$ on $Y^3$, the connected sum $Y'= Y\# (S^1\times S^2)$ carries such an $S^1$-valued Morse function: one thinks of $ Y'$ as the result of removing two balls in $Y$ containing the maximum and minimum, and gluing back $[-1,1]\times S^2$. One considers only those $\spinc$--structures $\mathfrak{t}$ on $Y'$ such that $ c_1(\mathfrak{t})$ evaluates as $2$ on a 2--sphere in the added handle. The Lagrangian boundary conditions which we use to define our Floer homologies then reduce (up to smooth isotopy, and probably also up to Hamiltonian isotopy) to the Heegaard tori $\mathbb{T}_\alpha$, $\mathbb{T}_\beta$ in $\sym^g(\Sigma)$, where $\Sigma =f^{-1}(3/2) \subset Y$. This viewpoint may give some insight into the cobordism maps in Heegaard Floer homology, which can be computed as Lagrangian matching invariants for broken fibrations of a particular kind. In this framework there is no need to decompose cobordisms into their elementary pieces.

\subsection{Broken pencils and fibrations}
\begin{Defn}\label{BF}
{\bf A broken fibration} $(X,\pi)$ on a smooth, compact, oriented, four--manifold $X$ (possibly with boundary) is a smooth map $\pi\colon X\to S$ to a compact surface $S$ such that the set of critical points $X^\crit$ is the union of a discrete set $D\subset \mathrm{int}(X)$ and a one-submanifold $Z \subset \mathrm{int}(X)$. These are constrained as follows:
\begin{itemize}
\item
For each $x\in D$, there exist positively oriented local coordinate charts 
\[\psi_x\colon (\C^2,0) \to (X,x),\quad \xi_{\pi(x)}\colon (\C,0) \to (S,\pi(x)),\] 
such that $ \xi_{\pi(x)}^{-1} \circ \pi \circ \psi_x $ coincides, near the origin, with the map $(z_1,z_2)\mapsto z_1 z_2$.
\item
For each point $z\in  Z$, there exist positively oriented local coordinate charts 
\[\phi_z\colon (\R\times \R^3,0) \to  (X,z),\quad \xi_{\pi(z)}\colon (\R^2,0) \to (S,\pi(z)),\] such that $ \xi_{\pi(z)}^{-1} \circ \pi \circ \phi_x $ coincides, near the origin, with one of the two maps
\begin{equation}\label{model} (t; x_1,x_2,x_3 )\mapsto (t, x_1^2 + x_2^2 - x_3^2 ), \quad (t; x_1,x_2,x_3 )\mapsto (t, -x_1^2 - x_2^2 + x_3^2 ). \end{equation}
\item
$\pi(Z)\subset S$ is an embedded 1--submanifold disjoint from $\pi(D)$. Furthermore, $\pi$ maps each component of $Z$ diffeomorphically to its image.
\item
There exists $w\in H^2(X;\R)$ such that $\langle w, h \rangle >0$ for every $h\in H_2(X)$ which is represented by a connected component of a fibre of $\pi$.
\end{itemize}
\end{Defn}
\begin{Rk}
One can contemplate many variations on this definition. For instance, there is no good reason to exclude the possibility that $Z$ contains arcs transverse to $\partial X$, though it was expedient to do so here. We believe that the techniques presented here can absorb the greater generality that that would entail, and we plan to explain this in a future paper.
\end{Rk}
These fibrations were introduced (under the name `singular Lefschetz fibrations') by Auroux, Donaldson and Katzarkov \cite{ADK}. As they showed, the cohomological hypothesis implies that there exists a \emph{near-symplectic form} $\omega \in Z^2(X)$, a closed two--form such that $\omega^2_x>0$ for $x\in X\setminus Z$ and $\omega_z=0$ for $z\in Z$, positive on the fibres at regular points. (If $Z= \emptyset$, this reduces to an observation of Gompf's about Lefschetz fibrations.)

\begin{figure}[t!]
\labellist
\small\hair 2pt
\pinlabel $W$ at 150 150
\pinlabel $\Sigma$ [l] at -20 218
\pinlabel $\Sigma_0$ at 170 218
\pinlabel $\bar{\Sigma}$ at 345 218
\pinlabel $L$  at 40 176
\pinlabel {Morse function}  at 300 40
\endlabellist
\centering
\includegraphics[scale=0.6]{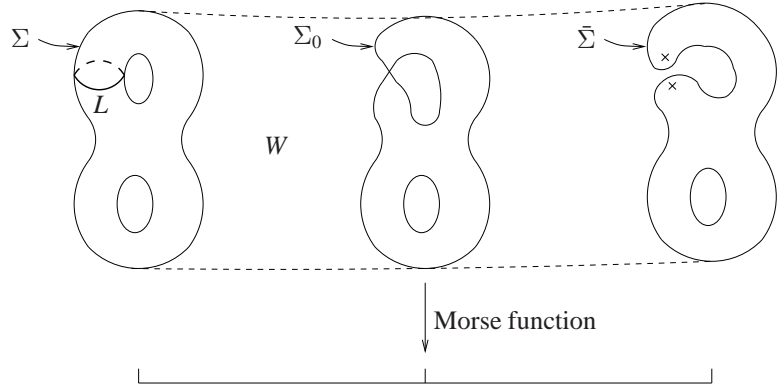}
\caption{\label{surgery fig} As one crosses a circle of critical values, the topology of the fibres changes by surgery. The singular fibres have conical singularities.}
\end{figure}

\begin{Ex}
Suppose that $Y$ is a closed, oriented three--manifold, and $f\colon Y \to S^1$ a Morse function such that all critical points have index 1 or 2. Then $f\times \id_{S^1}\colon Y\times S^1\to S^1\times S^1$ is a broken fibration. The cohomological condition (4) holds because every regular point lies on a loop $\gamma$ such that $df(\dot{\gamma})>0$.
\end{Ex}
The topology of the fibres of a broken fibration changes as one crosses a circle of critical values. The preimage of a transverse arc is a three--manifold with boundary, equipped with a Morse function, as illustrated in Figure \ref{surgery fig}.

Near each circle of critical points, there is a `attaching surface'\footnote{$Q$ is the surface along which one would attach
a round two-handle $S^1 \times D^1 \times D^2$ in giving a round handle decomposition of the broken fibration.} $Q$, cut out in the coordinates
$(t;x_1,x_2,x_3)$ as $\{ x_3=0;\, x_1^2+ x_2^2 =  \epsilon\}$. This can be a torus or a Klein bottle.
In the last example these surfaces are tori, but in the next one we get a Klein bottle.
\begin{figure*}[t!]
\small\hair 2pt
\labellist
\pinlabel $\gamma_2$ at 145 270
\pinlabel $\gamma_1$ at 220 270
\pinlabel $\delta$ at 165 230
\pinlabel $v_1=b$ at 545 170
\pinlabel $v_2=b+2a$ at 695 240
\pinlabel $a$ at 657 200
\pinlabel $a$ at 365 195
\pinlabel $=$ at 465 207
\pinlabel $=$ at 115 52
\pinlabel {monodromy $= \tau_{v_2}\circ \tau_{v_1}$} at 480 300
\pinlabel {monodromy $= \id$} at 50 130
\endlabellist
\centering
\includegraphics[scale=0.55]{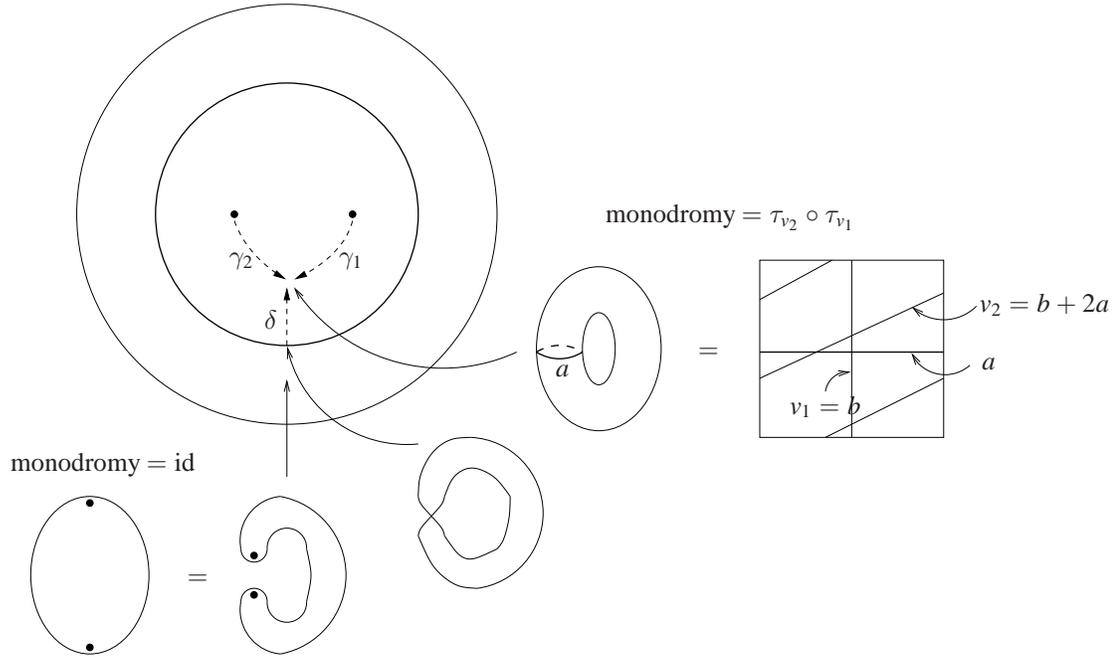}
\caption{\label{Klein fig} A broken fibration over $S^2$, with tori and spheres as regular fibres. Only one hemisphere of the base is shown. The loop $a$ sweeps out a Klein bottle (over the `tropic of Capricorn') which collapses to a circle of critical points over the equator.}
\end{figure*}

\begin{Ex}\label{klein}
(See Figure \ref{Klein fig}.) Let $T$ be a $2$--torus, and $a$, $b$ the standard loops generating $H_1(T;\Z)$. Let $\pi\colon E \to \Delta$ be a Lefschetz fibration over the disc with regular fibre $T=\pi^{-1}(1)$ and two critical points whose  vanishing cycles $v_1$ and $v_2$ (for vanishing paths $\gamma_1$ and $\gamma_2$) represent the homology classes $[b]$ and $[b]+2[a]$. We may then suppose that the  $v_i$ both meet $a$ transversely at a single point, as shown in the figure. 

The monodromy of the fibration is the composite of positive Dehn twists: $m=\tau_{v_2}\circ \tau_{v_1}$. By the Picard--Lefschetz formula, $m_*([a])=-[a]$; hence $m(a)$ is isotopic to the curve $a$, with  reversed orientation. Therefore $a$ sweeps out a Klein bottle $Q\subset \partial E$. One can extend $E$ to a broken fibration over a larger disc $D$ in which $Q$ collapses to a circle of critical points, mapping to a circle parallel to $\partial D$. This is what is depicted in the figure ($\delta$ is a `vanishing path' to the critical circle; associated with it is the attaching circle $a$, shown as a heavy line in the two pictures of the fibre). The fibres over $\partial D$ are two--spheres, and they contain a distinguished `braid' (two points on each fibre, corresponding to the points $\{ (t; 0,0, \pm \epsilon) \}$ in the local model \ref{model}). This braid is isotopic to a trivial braid in $S^1\times S^2$. We can complete the fibration to a broken fibration $\pi'\colon X\to S^2$ by gluing on a trivial $S^2$-bundle over the north-polar disc, using a gluing map which trivialises the braid.

The resulting four--manifold $X$ is simply connected with Euler characteristic $\chi(X)=4$ and signature $\sigma(X)=0$ (because $\pi'$ has a section). Indeed, it has a square-zero section, and since this and the fibre form a $\Z$-basis for homology, the intersection form must be even. Hence $X$ is a homotopy--$S^2\times S^2$. 
\end{Ex}
Since it is of some interest for the uniqueness problems for broken fibrations, we now verify that the manifold $X$ of the last example really \emph{is} $S^2\times S^2$.
\begin{Prop}
$X$ is diffeomorphic to $S^2\times S^2$. 
\end{Prop}
\begin{pf}
Think of the base as $\C\cup \{\infty\}$, with the critical values lying in $\C$ in the manner depicted in Figure \ref{Klein fig}. 
Consider the circles $c_x=\{ \real(z)  = x \}\cup \{\infty\} \subset \C\cup\{\infty\}$; these appear as vertical rulings of the diagram. For  $C\gg 0$, the preimages $\pi'^{-1}(\bigcup_{x\leq -C}{c_x})$ and  $\pi'^{-1}(\bigcup_{x\geq C}{c_x})$ are both diffeomorphic to $D^2\times S^2$. We claim that $\pi'^{-1}(\bigcup_{-C\leq x\leq C}{c_x})$ is a trivial cobordism. To see why the claim holds, let $x_0$ be the least $x$ such that $c_x$ hits the critical circle. For small positive $\epsilon$, $\pi'^{-1}(c_{x_0+\epsilon})$ is diffeomorphic to $(S^1\times S^2 ) \cup e_1$, where $e_1$ is a 1-handle; and 
 $\pi'^{-1}(\bigcup_{x_0-\epsilon \leq x\leq x_0+\epsilon}{c_x})$ is an elementary cobordism. This is left to the reader to see, but we do point out the belt-sphere of the handle attachment (a $2$--sphere in $\pi'^{-1}(c_{x_0+\epsilon})$): it has a circle in a torus fibre, isotopic to $a$, as its equator, which is then pinched off at the two critical points of $\pi'$ lying over $c_{x_0+\epsilon}$. For some $x_1>x_0$, $c_{x_1}$ passes through the first isolated critical value, and $\pi'^{-1}(c_x)$ changes by adding a 2-handle along $v_1$ (with framing given by the fibre-framing minus 1). This cancels with the 1-handle, because its attaching circle $v_1$ intersects the belt sphere transversely at a point. Thus
$\pi'^{-1}(\bigcup_{-C\leq x\leq x_1+\epsilon}{c_x})$ is a trivial cobordism. For the same reason, $\pi'^{-1}(\bigcup_{x_1+\epsilon \leq x\leq C }{c_x})$ is a trivial cobordism. The claim follows.

We deduce that $X'$ is diffeomorphic to  $(D^2 \times S^2) \cup_\phi (D^2 \times S^2)$, for some orientation-reversing self-diffeomorphism $\phi$ of $S^1\times S^2$. We can finish the proof rapidly but heavy-handedly by invoking Hatcher's theorem on $\diff(S^1\times S^2)$ \cite{Hat}, which implies that $\phi$ is isotopic to one of two standard maps. These maps yield $S^2\times S^2$ and $\mathbb{CP}^2\# \overline{\mathbb{CP}}^2$, the second of which does not have even intersection form. 
\end{pf} 

A notion closely related  to broken fibrations is that of `broken pencils':
\begin{Defn}
A {\bf broken pencil} is a triple $(X,B,\pi)$, where $B\setminus \mathrm{int}(X)$ is a  discrete subset,
and $\pi\colon X\setminus B\to S^2$ a map whose critical points conform to the models in points (1) and (2) of Definition \ref{BF},
and satisfy the condition (3). The model near a point of $B$ is the projectivisation map
$\C^2\setminus \{0\} \to \mathbb{CP}^1$, $(z_1,z_2)\mapsto (z_1: z_2)$. Condition (4) is also imposed,
where by a `fibre' of $\pi$ we mean the closure in $X$ of $\pi^{-1}(\mathrm{pt})$.
\end{Defn}
Again, broken pencils are near-symplectic, and the remarkable result of \cite{ADK} is that the converse is true: on a closed near-symplectic four--manifold $(X,\omega)$, there exist broken pencils $\pi$ whose one-dimensional critical set $Z$ coincides with $\omega^{-1}(0)$. Moreover, one can take $\pi(Z)$ to be a single circle in $S^2$,
and $\pi$ to be `directional' in the sense that only one of the two models in (2) is invoked.

After blowing up $X$ along $B$, the composite of the blow-down map $\widehat{X} \to X$ with $\pi$ extends smoothly to a broken fibration $\widehat{X}\to S^2$. The exceptional spheres are sections of it. Since near-symplectic forms exists as soon as $b^+(X)>0$, the conclusion is that, for any non-negative-definite $X$, $X \# N\overline{\mathbb{CP}}^2$ admits a broken fibration for any sufficiently large $N$.
(After the basepoints are exhausted, one can go on blowing up using a simple procedure which produces reducible nodal fibres.)

Roughly speaking, broken pencils are to near-symplectic forms as Lefschetz pencils are to symplectic forms. However, basic questions remain unanswered:
\begin{itemize}
\item
Which connected four--manifolds support broken pencils with connected fibres? (The fibres over one hemisphere are certainly connected.)
\item
Auroux--Donaldson--Katzarkov's sequences of broken pencils seem not to be `asymptotically unique'. How are different sequences of pencils on the same manifold related? A more elementary problem is to write down a set of `moves' which generate many broken fibrations starting from a given one.  Example \ref{klein} can easily be adapted to give a procedure which  
increases by 1 the genus of the fibres of a given fibration over a small disc, and does not change the diffeomorphism-type. This provides one such move.
\item
Which smooth four--manifolds support broken pencils if we drop the cohomological condition (4)? (All of them?)
\end{itemize}
\begin{Rk}
A development relevant to the second and third questions has occured since the first version of this paper appeared: Gay and Kirby \cite{GK} have produced \emph{achiral} broken fibrations (not usually satisfying condition (4)) on arbitrary smooth closed four--manifolds. Any embedded surface of self-intersection zero can be realised as a fibre,
and its framing can apparently be chosen at will; so, for example, both homotopy classes of maps $S^4 \to S^2$ are represented by achiral broken
fibrations.
\end{Rk}
Our invariants will be defined for broken fibrations rather than pencils. One \emph{could} define invariants of broken pencils from ours simply by blowing up the base locus $B$ (and relating $\spinc$ structures on $X$ and $\widehat{X}$ in the usual way). Then the invariants of $X$ and $\widehat{X}$ would be related by the blow-up formula familiar from Seiberg--Witten theory. However, we will have no more to say about this.

\subsection{Outline of the construction}
Our construction has three stages:
\begin{enumerate}
\item[(I)]
Constructing Lagrangian boundary conditions for pseudo-holomorphic curves.
\item[(IIa)]
Properties of the moduli spaces of pseudo-holomorphic curves (transversality, compactness, orientation, etc.).
\item[(IIb)]
Algebraic formulation of the invariants. 
\end{enumerate}
Stage (I), the subject of the present paper, is perhaps the most interesting, for it is not obvious how to generalise the moduli spaces studied by Donaldson and Smith. For (IIa), carried out in the sequel \cite{Pe1}, standard techniques from the theory of pseudo-holomorphic curves suffice, at least if one is prepared to make numerical hypotheses on the $\spinc$--structures considered. Stage (IIb) is mostly `soft' topology.

\subsubsection{The Lagrangian correspondences} 
\begin{figure*}[t!]
\labellist
\small\hair 2pt
\pinlabel $x^+$ at 135 500
\pinlabel $x^-$ at 135 485
\pinlabel $x$ at 135 285
\pinlabel $\Delta$ at 35 110
\pinlabel $0$ at 180 63
\pinlabel $1$ at 425 63
\pinlabel $L$ at 348 275
\pinlabel $\Sigma=E_1$ at 460 240
\pinlabel $\bar{\Sigma}\cong \widetilde{E}_0$ at 290 500
\pinlabel $E_0$ at 270 240
\pinlabel $n$ at 190 345
\pinlabel {normalisation map} at 280 345
\endlabellist
\centering
\includegraphics[scale=0.6]{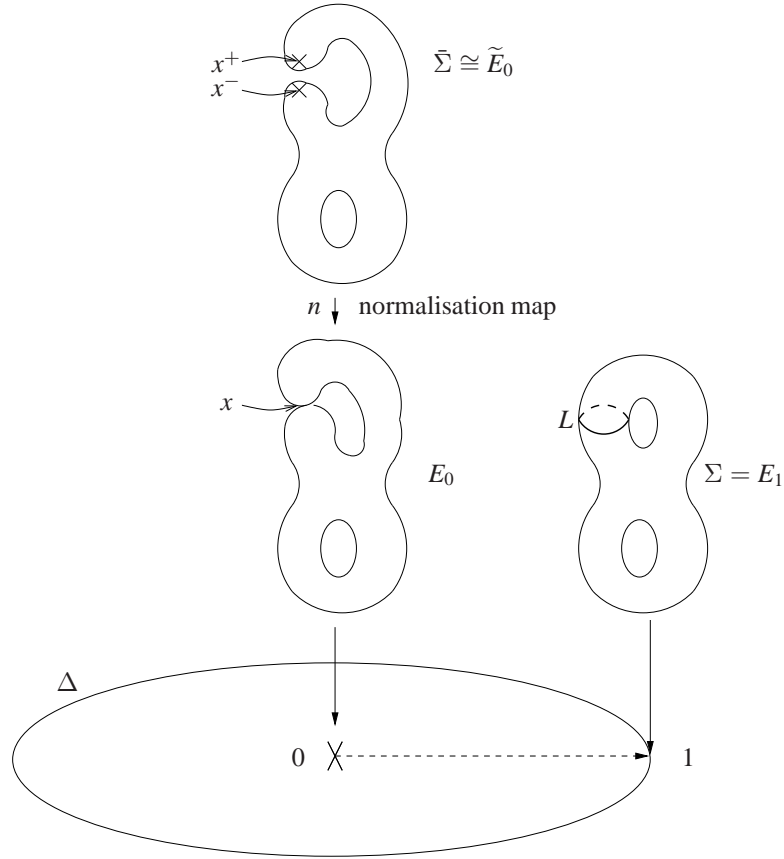}
\caption{\label{lefschetz figure} An elementary Lefschetz fibration $E\to \Delta$ with smooth fibre $\Sigma$ and vanishing cycle $L$. We recover $\bar{\Sigma}$ as the normalisation of its central
fibre.}
\end{figure*}

For (I), the key observation (Theorem A, in Section \ref{Lag corr}) is that,
if $\Sigma$ is a closed Riemann surface, and $\bar{\Sigma}$ a Riemann surface obtained by surgering out a circle
$L\subset \Sigma$ (i.e. by excising an annulus-neighbourhood of $L$ and gluing in two discs) then there is a
distinguished Hamiltonian isotopy-class of Lagrangian correspondences
\[\widehat{V}_L \subset \sym^n(\Sigma)\times \sym^{n-1}(\bar{\Sigma}),\quad n= 1, 2,\dots  : \]
Lagrangian for a symplectic form of shape $(-\omega)\oplus\bar{\omega}$, where $\omega$ and $\bar{\omega}$ are K\"ahler forms lying in
certain cohomology classes. $\widehat{V}_L$ arises as a \emph{vanishing cycle} for a symplectic degeneration of $\sym^n(\Sigma)$, as follows.
Form a holomorphic Lefschetz fibration $(E,\pi)$ over the closed unit disc $\Delta$, with smooth fibre $\Sigma=E_1$ and
vanishing cycle $L\subset E$, as in Figure \ref{lefschetz figure}.
The normalisation of the nodal fibre can then be identified with $\bar{\Sigma}$ by a diffeomorphism which is canonical up to isotopy.

The family $\sym^n_\Delta(E)\to\Delta$ of symmetric products of the fibres is a globally singular space, but,
as observed by Donaldson and Smith, it has a resolution of singularities $\hilb^n_\Delta(E)\to \sym^n_\Delta(E)$,
the \emph{relative Hilbert scheme of $n$ points}, which fits into a commutative diagram
\[\xymatrix{
\hilb^n_\Delta(E) \ar[rr] \ar[dr]  && \sym^n_\Delta(E)\ar[dl]  \\
& \Delta.
}\]
A point on the Hilbert scheme is a pair $(s,\mathcal{I})$, where $s\in \Delta$ and $\mathcal{I}$
is an ideal sheaf in $\mathcal{O}_{E_s}$ such that $\sum_{x\in E_s}{\dim_\C(\mathcal{O}_{E_s,x}/\mathcal{I}_x)}=n$.
The natural `cycle map' $\hilb^n_\Delta(E)\to \sym^n_\Delta(E)$ is bijective except over $0\in \Delta$; it partly resolves
the fibre over $0$. The crucial observation, which is mentioned in passing by Smith \cite[Proposition 3.7]{Smi}, is that the critical manifold
of the natural map $\hilb^n_\Delta(E)\to \Delta$ (that is, the normal crossing divisor in the zero-fibre $\hilb^n(E_0)$) is
naturally biholomorphic to $\sym^{n-1}(\widetilde{E}_0)$.

\emph{This hidden link between $\sym^n(\Sigma)$ and $\sym^{n-1}(\bar{\Sigma})$ is the starting point for our construction.}

Choose a K\"ahler form on $\hilb^n_\Delta(\Sigma)$.  The Lagrangian correspondence $\widehat{V}_L$ between $\sym^n(\Sigma)$ and $\sym^{n-1}(\bar{\Sigma})$ is defined as the graph of symplectic parallel transport, over the ray $[0,1]\subset \Delta$, into the critical set $\sym^{n-1}(\widetilde{E}_0)$.

The relative Hilbert scheme belongs to a class of symplectic degenerations which we call `symplectic Morse--Bott fibrations'. Their geometry is developed in Section 2. The specific geometry of $\widehat{V}_L$ is explored in Section 3. For a bare definition of the invariants, a shorter treatment would suffice, but our aim is to get a grip on these decidedly slippery correspondences, both for the sake of intuition and as a foundation for future work. For example, the link with Heegaard Floer homology requires a good deal of control over the vanishing cycles.

\begin{Rk}\label{gauge lags}
One can also construct (embedded? non-singular?) Lagrangian correspondences between symmetric products---equipped with their canonical K\"ahler forms arising from their interpretation as vortex moduli spaces---via the Seiberg--Witten equations on (a metric completion of) the elementary 3-dimensional cobordism, compare \cite{Sa3}. It would be very interesting to understand the relation of these correspondences with ours; this could also be a first step in relating Lagrangian matching and Seiberg--Witten invariants.
\end{Rk}

\subsubsection{A moduli space} For concreteness, we will focus here on an `elementary broken fibration' $(X_0,\pi)$ over a cylinder $C=S^1\times [-1,1]$. This has just one circle of critical points, mapping diffeomorphically to $C^\crit = S^1\times \{0\}$. Let $Y= \pi^{-1}(S^1\times \{-1\})$ and $\bar{Y}=\pi^{-1}(S^1\times \{1\})$. Let us suppose that the fibres $Y_t$ have genus $g$, and that the fibres $\bar{Y}_t$ genus $g-1$. (One can perfectly well reverse the orientation of $C$, but things run a little more smoothly if we suppose that the fibres are connected.)

Inside $Y$, there is the attaching surface $Q$---a torus or a Klein bottle---formed from circles in the fibres $Y_t$ which shrink to points in the circle $X_0^\crit$.\footnote{I am indebted to Paul Seidel for the idea that this surface could serve as a boundary condition for holomorphic curves.}
Let
\[ Y^{[n]} =\sym_{S^1}^{n}(Y),\quad \bar{Y}^{[n]}=\sym^{n-1}_{S^1}(\bar{Y}) \]
be the `associated bundles' of symmetric products of the fibres. These become smooth when one chooses complex structures on the fibres of $Y\to S^1$ and $\bar{Y}\to S^1$. By applying the construction of Lagrangian correspondences \emph{simultaneously} to all the circles $Q\cap Y_t \subset Y_t$, we obtain a sub-fibre bundle
\[ \EuScript{Q} \subset Y^{[n]} \times_{S^1}\bar{Y}^{[n-1]}. \]
There exist closed, fibrewise K\"ahler forms $\Omega$ on $Y^{[n]}$ and $\bar{\Omega}$ on $\bar{Y}^{[n-1]}$ such that $\EuScript{Q}$ is globally isotropic with respect to $(-\Omega)\oplus \bar{\Omega}$. This construction is the content of Theorem B (in Section \ref{LMC from BF}).

We obtain a moduli spaces of holomorphic curves as follows. Consider the two spaces
\[ X^{[n]}= Y^{[n]}\times (-\infty,0],\quad \bar{X}^{[n]}=\bar{Y}^{[n]}\times [0,\infty),  \]
equipped with the forms obtained from $\Omega$ and $\bar{\Omega}$ by pulling back. These spaces fibre respectively over $C^-:=S^1\times (-\infty,0]$ and over $C^+:= S^1\times (\infty,0]$. Choose asymptotically translation-invariant almost complex structures on the vertical tangent bundles of $X^{[n]}$ and $\bar{X}^{[n]}$, compatible with the fibrewise symplectic forms. Consider pairs $(u, \bar{u})$ where $u$ (resp. $\bar{u}$) is a pseudo-holomorphic section of  $X^{[n]}\to C^-$ (resp. $\bar{X}^{[n]} \to S^1\times C^+$). Both are required to have finite energy. Neither $u$ nor $\bar{u}$ is \emph{individually} constrained by a boundary condition over $S^1\times \{0\}$. However, the \emph{pair} $(\partial u ,\partial\bar{u}): = (u|S^1\times\{0\}, \bar{u}|S^1\times\{0\})$, which is a section of $Y^{[n]} \times_{S^1}\bar{Y}^{[n-1]}$, is required to be a section of the sub-fibre bundle $\EuScript{Q}\to S^1$.

It is possible to interpret this boundary condition as a Lagrangian boundary condition for a pseudo-holomorphic half-cylinder in a symplectic manifold. In particular, its linearisation is Fredholm.

There are three sources of non-compactness for the moduli space of pairs $(u,\bar{u})$:
(i) a sequence may converge to a `broken trajectory' in the sense of Floer theory;
(ii) there may be bubbles in interior fibres;
(iii) there may be boundary bubbles, i.e. holomorphic discs in the fibres of $Y^{[n]} \times_{S^1}\bar{Y}^{[n-1]}$
with boundary on  $\EuScript{Q}$. Of these, (i) is an essential feature, built into the algebra of
Floer homology; we will rule out (ii) and (iii) on `monotonicity' grounds providing that $n\geq g(\Sigma)$,
and on `weak monotonicity' grounds when $n\leq (g-1)/2$, using the methods due to Hofer--Salamon and Lazzarini.
The range $(g-1)/2< n < g-1$ would take us into the realm of virtual moduli chains;
boundary-bubbling would be a serious issue. The case $n=g-1$ is a little better, but there could be complicated wall-crossing phenomena.

\begin{Rk}
How are we to understand the boundary condition $\EuScript{Q}$? A section $\gamma$ of $\EuScript{Q}$ projects to an $n$-fold multi-section $\gamma$ of $Y$, and an $(n-1)$-fold multi-section $\bar{\gamma}$ of $\bar{Y}$; these have homology classes $[\gamma]\in H_1(Y;\Z)$ and $[\bar{\gamma}]\in H_1(\bar{Y};\Z)$. It turns out that there is a standard relative homology class $\beta\in H_2(X_0,\partial X_0 \cup Z;\Z)$, where $Z=X_0^\crit$ is the circle of critical points, such that
\[ \partial \beta = [\bar{\gamma}] - [\gamma] + [Z].\]
Thus there is a surface in $X_0$ which `tunnels' between $\gamma$ and $\bar{\gamma}$ and has one further boundary component, along $Z$. This can be chosen to be a disjoint union of cylinders. Loosely speaking, $\EuScript{Q}$ encodes the condition that $\gamma$ `matches' with $\bar{\gamma}$ in the sense that there is such a surface joining them.
This notion of matching boundaries explains the name of our invariant, which is intended to be suggestive rather than literal. It also makes a weak link with Taubes' programme---at the level of homology, not of moduli spaces.
\end{Rk}

As Floer homology experts will realise, there are symplectic Floer homology modules
\[ HF_*(Y^{[n]},\Omega),\quad HF_*(\bar{Y}^{[n-1]},\bar{\Omega}) \]
associated with the bundles of symmetric products $Y^{[n]}$ and $\bar{Y}^{[n-1]}$.
In general, these are modules over the universal $\Z$-Novikov ring $\Lambda_{\Z}$. Our moduli space gives rise to a homomorphism
\[ \Phi\colon HF_*(Y^{[n]},\Omega )\to HF_*(\bar{Y}^{[n-1]},\bar{\Omega}). \]
This is a raw form of the Lagrangian matching invariant of an elementary broken fibration. (If one decomposes the Floer homology groups by topological sectors ($\spinc$--structures) it is possible to choose the two--forms in such a way that the group associated with a fixed sector is defined and finitely generated over $\Z$ rather than $\Lambda_{\Z}$.)

By considering the same broken fibration over the orientation-reversed cylinder, one gets another map $\Psi$ running in the opposite direction. This map is the adjoint of $\Phi$ with respect to the Poincar\'e duality on Floer homology. 

It seems highly likely that computational consequences in Floer homology can be extracted from Lagrangian matching invariants, possibly in the form of exact triangles (a precise conjecture is stated in Part II).

\subsubsection{Lagrangian matching invariants for compact manifolds}

The Lagrangian boundary condition $\EuScript{Q}$ serves equally well when one has a broken fibration over a compact surface. For example, given a broken fibration $(X,\pi)$ over $S^2$, with just one circle of critical points $Z$, mapping to an `equator' $\pi(Z)$, one parametrises a neighbourhood $C$ of $\pi(Z)$ as a cylinder. Then $S^2\setminus \mathrm{int}(C)= D^+\cup D^-$, where $D^\pm $ are closed discs. Let $X^\pm=\pi^{-1}(D^\pm)$, and suppose that the regular fibres over $D^+$ (resp. $D^-$) have genus $g$ (resp. $g-1$). Then one considers pairs $(u^+,u^-)$, where $u^+$ (resp. $u^-$) is a pseudo-holomorphic section of the relative Hilbert scheme $\hilb^n_{D^+}(X^+)$ (resp. $\hilb^{n-1}_{D^-}(X^-)$. The pair of boundary values $(\partial u^+,\partial u^-)$ is again required to lie on a Lagrangian $\EuScript{Q}$.

We should keep track of the `topological sectors' of the moduli space, namely, the homotopy classes of pairs of smooth sections $(u^+,u^-)$ with boundary on $\EuScript{Q}$. Each topological sector $\beta$ distinguishes a $\spinc$--structure $\mathfrak{s}_\beta\in \spinc(X)$ (that is, the homotopy class of a lift of the classifying map for the tangent bundle from $B\GL^+(4,\R)$ to $B\spinc(4)$).
\begin{TheoremD}\cite{Pe1}
The expected dimension of the moduli space for the sector $\beta$ is the number
\[ d(\mathfrak{s}_\beta)=\frac{1}{4}(c_1(\mathfrak{s}_\beta)^2-2\chi(X)-3\sigma(X)).\]
\end{TheoremD}
This number is familiar to gauge theorists as the dimension of the Seiberg--Witten moduli space.

One can cut down the moduli space to zero dimensions by insisting that $u^+$ passes through certain cycles in marked fibres of the relative Hilbert scheme. The cut down moduli space is (for certain $\mathfrak{s}_\beta$) compact. The Lagrangian matching invariant is a count of its points. It also keeps track of homological information about the cycles used to cut down the moduli space.

\begin{Rk}
We should warn the reader of the potentially confusing point that relative Hilbert schemes play two distinct roles in this story. We study pseudo-holomorphic sections of relative Hilbert schemes of points on Lefschetz fibrations; however, the boundary conditions for these sections are also derived from relative Hilbert schemes of points on elementary Lefschetz fibrations.
\end{Rk}

\subsubsection{Formulation of the invariants}
Let $\spinc(X)$ denote the $H^2(X;\Z)$--torsor of (isomorphism classes of) $\spinc$--structures on the oriented four--manifold $X$. When $X$ is given a structure of broken fibration, with $F\in H_2(X;\Z)$ the class of a regular fibre, we write
\[ \spinc(X)_k = \{ \mathfrak{s}\in \spinc(X) : \langle c_1(\mathfrak{s}), F \rangle = 2k , (*) \}, \]
where $(*)$ is the condition that for any connected component $\Sigma$ of a regular fibre, one has $\langle c_1(\mathfrak{s}),[\Sigma]\rangle \geq \chi(\Sigma)$. 

\begin{Defn} $k\in \Z$ is {\bf admissible} for $(X,\pi)$ if either (i) the fibres are all connected and $k>0$, or  (ii) $\chi(X_s) /2 < k < - \chi(X_s)/2$ for all $s\in S^{\mathrm{reg}}$. A $\spinc$--structure $\mathfrak{s}$ is admissible if $\mathfrak{s}\in \spinc(X)_k$ with $k$ admissible.
\end{Defn}

\begin{TheoremC}\cite{Pe1}
To a broken fibration $(X, \pi)$ over $S^2$, such that $\pi|X^{\crit}$ is injective, one can associate an invariant $\EuScript{L}_{(X,\pi)}$, the {\bf Lagrangian matching invariant}. 
This is a map
\[   \bigcup_{k \text{ admissible}} \spinc(X)_{k} \to  \mathbb{A}(X), \quad \mathfrak{s} \mapsto \EuScript{L}_{(X,\pi)}(\mathfrak{s}). \]
Here $\mathbb{A}(X)$ is the graded abelian group
\[ \Z[U] \otimes_{\Z} \Lambda ^*H^1(X;\Z),\quad \deg(U)=2. \] 
The element $ \EuScript{L}_{(X,\pi)}(\mathfrak{s}) $ is homogeneous of degree $d(\mathfrak{s})$. It is invariant under isotopies of $\pi$ through fibrations of the same type, and equivariant under isomorphisms $(X,\pi)\cong (X',\pi')$.

The same holds when the base is an arbitrary surface $S$, providing one replaces $\mathbb{A}(X)$ by $\mathbb{A}(X,\pi)= \Z[U]\otimes_{\Z}\Lambda^* \Hom(I_\pi,\Z)$, where $I_\pi\subset H_1(X;\Z)$ is the subgroup of classes supported on a fibre of $\pi$.
\end{TheoremC}

\begin{Rk} The condition that $\pi|X^{\crit}$ should be injective can always be achieved by perturbing a given fibration so that a multiply-covered circle of critical values becomes a family of parallel circles. Unfortunately, we do not prove here that different perturbations give the same $\EuScript{L}_{(X,\pi)}$, so we cannot remove this restriction. What is needed is a commutativity property for Lagrangian correspondences (Conjecture \ref{commutativity}), which should be provable by fine-tuning of symplectic forms.
\end{Rk}
The format of this invariant is familiar. The Seiberg--Witten invariant of a four--manifold with $b^+\geq 1$ can be formulated as a map  
\[ \mathrm{SW}_X \colon \spinc(X)\to \mathbb{A}(X),\] 
where $\mathrm{SW}_X(\mathfrak{s})$ is homogeneous of degree $d(\mathfrak{s})$. In this formulation, $\mathrm{SW}_X(\mathfrak{s})$ is the fundamental homology class of the Seiberg--Witten moduli space in the ambient configuration space $\EuScript{B}^*_{X,\mathfrak{s}}$, under isomorphisms
\[ H_*(\EuScript{B}^*_{X,\mathfrak{s}};\Z) = H_*(B\EuScript{G}_X;\Z) \cong \mathbb{A}(X). \]
Here $\EuScript{G}_X=\mathrm{Map}(X,S^1)$ is the gauge group.
We should make two standard caveats: (i) The overall sign of $\mathrm{SW}_X$ depends on a homology orientation for $X$; (ii) if $b^+(X)=1$, one has to choose a `chamber' in the space of auxiliary parameters.

\begin{Conj}
Let $\mathfrak{s}\in\spinc(X)_k$ with $k\geq 1 $ admissible. Then 
\[ \EuScript{L}_{X,\pi}(\mathfrak{s}) = \pm   i  (\mathrm{SW}_X (\mathfrak{s})) \] 
on $\spinc(X)_k$, $k\geq 1$, where $i\colon \mathbb{A}(X)\to \mathbb{A}(X,\pi)$ is the map induced by $I_\pi \hookrightarrow H_1(X;\Z)$, and the sign is independent of $\mathfrak{s}$.
In particular, $\EuScript{L}_{X,\pi}$ depends only on $X$ and not on $\pi$.  
(When $b_2^+(X)=1$, the right-hand side is calculated in the `Taubes chamber' of a compatible near symplectic-form.)
\end{Conj} 
The assumption $k\geq 1$ is somewhat arbitrary; the conjecture might well hold for many (if not all) of the other $\spinc$--structures for which $\EuScript{L}_{X,\pi}$ is defined.

It is already known, by combining theorems of Usher and Taubes, that the Donaldson--Smith invariant is the Seiberg--Witten invariant for high-degree Lefschetz pencils. In the case of broken fibrations which arise from Morse functions by crossing with $S^1$, the Lagrangian matching invariant can be understood via $(2+1)$--dimensional TQFT methods. 

\begin{TheoremE}\cite{Pe1}
Let $M_K$ be a 3-manifold obtained by zero-surgery on a knot, $f\colon M_K\to S^1$ a Morse function with critical critical points of index one and two such that $f^*[S^1]$ is a generator for $H^1(M_K;\Z)$. Consider the broken fibration $\id\times f\colon S^1\times M_K\to S^1\times S^1$, and let $\mathfrak{s}$ be an admissible $\spinc$--structure. If $c_1(\mathfrak{s})$ is the pullback of $k$ times the generator of  $H^2(M_K;\Z)$, then $d(\mathfrak{s})=0$ and 
\[ \pm \EuScript{L}_{S^1\times M_K,\id\times f}(\mathfrak{s}) = \sum_{i\geq 1}{ i a_{k+i}} = \pm \mathrm{SW}_{S^1\times M_K}(\mathfrak{s}), \]
where $a_0 + \sum_{i}{a_i (t^i+ t^{-i})}$ is the normalised Alexander polynomial of $K$.
\end{TheoremE}
There is also a vanishing theorem for the connected sum (not fibre sum!) of broken fibrations:
\begin{TheoremF}\cite{Pe1}
Let $(X,\pi)$ and $(X',\pi')$ be broken fibrations over the same base. Then there is a broken fibration on the connected sum $X\# X'$ whose invariants vanish for any admissible $\spinc$--structure.
\end{TheoremF}
(The admissible $\spinc$--structures are fewer than one would like here, because there are disconnected fibres.)
\subsubsection{The field theory} This is constructed using the theory of Floer homology for symplectic automorphisms, applied to relative symmetric products of fibred three--manifolds. See also Usher's closely related work \cite{Ush2}. The symplectic forms involved are chosen carefully so as to get a theory which is finitely generated over $\Z$ (for a fixed topological sector).
\begin{itemize}
\item
Let $(Y,\pi)$ be a closed, oriented 3-manifold fibred over a closed, oriented, one-manifold $T$. Let $ \mathfrak{t}\in\spinc(Y)_k$ where $k$ is $\pi$-admissible. To $(Y,\pi,\mathfrak{t})$ is assigned a finitely generated abelian group $HF_*(Y,\mathfrak{t})$ (which might depend on $\pi$, despite the notation). Write $HF^{\Q}_*(Y,\mathfrak{t})$ for $HF_*(Y,\mathfrak{t})\otimes \Q$. If $\iota\colon T\to T$ is an orientation-reversing diffeomorphism, and $\iota^*\pi\colon -Y\to T$ the resulting fibration, one has
\begin{equation}\label{duality} HF_*^{\Q} (-Y, -\mathfrak{t}) = HF^{\Q}_*(Y, \mathfrak{t})^*. \end{equation}
If $T= T_1 \amalg T_2$,  and $Y_i = Y|T_i$, then 
\begin{equation}\label{tensor}
HF^{\Q}_*(Y,\mathfrak{t})  = HF^{\Q}_*(Y_1,\mathfrak{t} |T_1) \otimes_{\Q}  HF^{\Q}_*(Y_2,\mathfrak{t}|T_2) .\end{equation} 
\item
$HF_*(Y,\mathfrak{t} )$ is graded by the $\Z$-set $J(Y,\mathfrak{t})$ of homotopy classes of oriented two-plane fields $\xi \subset TY$ underlying $\mathfrak{t}$. That is, $HF_*(Y,\mathfrak{t})$ is a direct sum $\bigoplus_{j\in J(Y,\mathfrak{t})}HF_j(Y,\mathfrak{t})$. (As observed by Kronheimer et al. in \cite[section 2.4] {KMOS}, $J(Y,\mathfrak{t})$ is naturally a transitive $\Z$-set with stabiliser $\Div(c_1(\mathfrak{t}))\Z$, where $\Div(c)$ is the divisibility of $c$ in $H^2(Y;\Z)$.)
\item
$HF_*(Y,\mathfrak{t})$ is a graded module over the graded ring $\Z[U]\otimes_{\Z}\Lambda^*I_\pi$, where $U$ has degree $-2$, and $I_\pi$ has degree $-1$. That is, $U \cdot HF_j \subset HF_{j-2}$ and $a\cdot HF_j \subset HF_{j-1}$ for $a\in I_\pi$.
\item
Let $(X,\pi)$ be a broken fibration over a compact surface-with-boundary $S$. Let $Y= \pi^{-1}(\partial S)$, and suppose $\mathfrak{s}\in \spinc(X)$ is admissible. Then 
there is a relative invariant, namely an element
\[ \EuScript{L}_{(X,\pi)} ( \mathfrak{s}) \in  HF_*(Y,\mathfrak{t}),\quad \mathfrak{t}=\mathfrak{s}|Y .   \]
$\EuScript{L}_{(X,\pi)}( \mathfrak{s})$ has degree zero in the sense that it lies in
$HF_{j} (Y,\mathfrak{t}) $, where $j\in J(Y,\mathfrak{t})$ is characterised by the existence of an almost complex structure on $X$, representing $\mathfrak{s}$ and preserving $j$ on $Y$.
\item
When $(X,\mathfrak{s})$ is a cobordism from $(Y_1,\mathfrak{t}_1)$ to $(Y_2,\mathfrak{t}_2)$, and $Y_i = \pi^{-1}(S_i)$ for a decomposition $\partial S= \partial S_1 \amalg \partial S_2$, we usually rewrite the relative invariant, using formulae (\ref{duality}, \ref{tensor}), as a group homomorphism
\[ \EuScript{L}_{(X,\pi)}(\mathfrak{s}) \in \Hom (HF^{\Q}_*(Y_1,\mathfrak{t}_1), HF^{\Q}_*(Y_2,\mathfrak{t}_2)) . \]
It then intertwines the action of $\Z[U]$. If classes $\alpha_1\in I_{\pi_1}$ and $\alpha_2\in I_{\pi_2}$ become homologous in $X$, then $\alpha_2\cdot \EuScript{L}_{(X,\pi)}(\mathfrak{s})(x) =\EuScript{L}_{(X,\pi)}(\mathfrak{s}) ( \alpha_1\cdot x) $. 
\item 
When $S$ separates along an embedded circle $T\subset S^{\mathrm{reg}}$, $S$ decomposes as $S= S_1\cup S_2$ and $X$ as $X= X_1\cup_{Y} X_2$, where $Y= \pi^{-1}(T)$. Orienting $T$ as the boundary of $S_1$, one has
\begin{equation}\label{composition} 
\EuScript{L}_{(X,\pi)}(\mathfrak{s}) = \EuScript{L}_{(X_2,\pi|X_2)}(\mathfrak{s}|X_2)  \circ \EuScript{L}_{(X_1,\pi|X_1)} (\mathfrak{s}|X_1) .   \end{equation} 
\item
Continue with the last point but now suppose $S$ is closed. One can use the duality property of the groups under orientation-reversal to write $\EuScript{L}_{(X_2,\pi|X_2)}(\mathfrak{s}|X_2)$ as a homomorphism
\[   \EuScript{L}^\vee_{(X_2,\pi|X_2)}(\mathfrak{s}|X_2)  \colon HF^{\Q} (Y,\mathfrak{s}|Y)\to \Q .\]
Then the invariant for the closed manifold $X$ is computed using the module structure as
\[  \EuScript{L}_{(X,\pi)} ( \mathfrak{s})  (U^p \otimes \lambda )  = \EuScript{L}^\vee_{(X_2,\pi|X_2)} \circ U^p \circ \lambda \circ ( \EuScript{L}_{(X_1,\pi|X_1) } \mathfrak{s}|X_1 ),  \]
where $\lambda \in \Lambda^* I_{\pi}$ is in the image of $\Lambda^* I_{\pi|Y}$.
\end{itemize}
\begin{Rk}
We have brought in rational coefficients only so as to avoid distracting Ext and Tor terms.
\end{Rk}
\begin{Rk}
Floer homology aficionados will want to know how $HF_*(Y,\mathfrak{t})$ compares with the groups arising in other Floer theories. D Salamon conjectured that the symplectic Floer homology groups of relative symmetric products of $(Y,\pi)$, defined using a closed two--form $\Omega\in \Omega^2(\sym^n_{S^1}(Y))$, should be isomorphic to the perturbed Seiberg--Witten monopole Floer homology of $Y$ \cite{Sa2}. Different forms $\Omega$ will correspond to different perturbations; one can make this precise by comparing the periods, cf. Y Lee's article \cite{Lee}.  The canonical, finitely-generated group $HF(Y,\mathfrak{t})$ should be isomorphic to  monopole Floer homology with `monotone perturbations':
\[ HF_*(Y,\mathfrak{t}) \stackrel{?}{\cong}  \HMto_\bullet(Y , \mathfrak{t}; [w] ). \]
Here $w$ is a closed two--form which is used to perturb the Chern-Simons-Dirac functional, chosen so that $[w]=\lambda c_1(\mathfrak{s})$, where $\lambda> 2\pi$ (compare \cite{Lee}). For such a $[w]$ there are no reducibles, so $\overline{HM}_\bullet=0$. Hutchings' periodic Floer homology groups (which, like $HF_*(Y,\mathfrak{t})$, involves the map $\pi$) are also thought to be isomorphic; Usher has made some progress in this direction \cite{Ush2}. Ozsv\'ath--Szab\'o's $HF^+$ is closely related, but not always the same (our groups can be non--zero for infinitely many $\mathfrak{t}$). 

It is natural to conjecture that, under the isomorphism between $HF_*$ and monopole Floer theory, the cobordism--maps should also coincide, i.e. that these are equivalent field theories on the cobordism--category of broken fibrations equipped with $\spinc$--structures.
\end{Rk}

\subsubsection{Further directions} There are other gauge theories which one could attempt to mimic using methods similar to those in this paper, notably $\SO(3)$ instanton theory. Indeed, Chris Woodward and Katrin Wehrheim are working on a closely related theory based on Lagrangian correspondences between spaces of flat connections.  The $\SU(2)$ instanton theory is more problematic; one needs a good way of dealing with the singularities in moduli spaces of flat connections over surfaces. It would also be interesting to develop the knot Floer homology of Ozsv\'ath--Szab\'o and Rasmussen from a `Lagrangian matching' point of view. As we have mentioned, Lagrangian matching invariants have applications within symplectic Floer homology; we plan to develop this by studying the exactness of certain triangles of Floer homology groups.

\subsubsection{Navigation} The main results of this paper---the ones which have a major bearing on its sequel---are Theorem A (which occurs near the beginning of Section \ref{Lag corr}) and Theorem B (in Section \ref{LMC from BF}). The results of Section \ref{Lag topology} will also be needed, though their role is of secondary importance. The `structure theorem' \ref{structure} is an important staging post in establishing Theorems A and B. Section 2 sets the stage; it also contains a non-trivial result, the `monodromy theorem' \ref{monod} which, however, is not part of the main logical thread. 

\subsubsection{Acknowledgements} This paper is a refined and reworked version of material from the author's Ph.D. thesis \cite{Per} (Imperial College London, 2005) under Simon Donaldson. It was he who proposed the problem of constructing invariants for broken fibrations. I owe many thanks to him for his generosity with time and ideas and for his patient encouragement. I am indebted to Paul Seidel, who has made invaluable suggestions at several points in the project, and to Ivan Smith for many useful discussions and comments. I acknowledge the strong influence on this work of papers by these same three people (chiefly \cite{DS}, \cite{SS} and \cite{Sei}). My thanks also to Michael Usher for telling me about his related work \cite{Ush2} and pointing out a good way of using a cohomology-class calculation from \cite{Per}; to Inanc Baykur for helpful comments; and to an eagle--eyed referee.

I also acknowledge support from EPSRC grant 41602.

\section{Fibred symplectic Picard--Lefschetz theory}

A leading role in this paper is played by certain \emph{Lagrangian correspondences}
\[\widehat{V}\subset ( \sym^n(\Sigma)\times \sym^{n-1}(\bar{\Sigma}) , \omega\oplus -\bar{\omega})\]
between symmetric products of Riemann surfaces $\Sigma$, $\bar{\Sigma}$, where $\chi(\bar{\Sigma})-\chi(\Sigma)=2$. The natural setting for these submanifolds is Picard--Lefschetz theory: they arise as vanishing cycles for algebro-geometric degenerations of $\sym^n(\Sigma)$. The degenerations are globally smooth, but have critical fibres with normal crossing singularities. Our approach emphasises the symplectic geometry of these degenerations; for us, they are examples of `symplectic Morse--Bott fibrations'. In this section we develop the geometry of such fibrations. Only the basics are required for the definition of Lagrangian matching invariants. The results on monodromy (\ref{twists}, \ref{monodromy}) are not, though they may well prove useful in computations. 

The work of Seidel and Smith \cite{SS} exploited closely related geometries. However, by working with affine algebraic varieties with $\C^*$-actions, they were able to circumvent several difficulties (but had to contend with an additional one, concerning parallel transport).

\begin{Defn}
{\bf(a)} A {\bf symplectic Morse--Bott fibration} $(E,\pi,\Omega,J_0,j_0)$  consists of a manifold $E^{2n+2}$ (possibly with boundary) and a smooth proper map $\pi\colon E\to S$ to an oriented surface $S$, mapping $\partial E$ submersively to $\partial S$; a closed two--form $\Omega$ on $E$; an almost complex structure $J_0$ in a neighbourhood of the set of critical points of $\pi$, $E^{\crit}\subset E$; and a positively oriented complex structure $j_0$ in a neighbourhood of the set of critical values $S^{\crit}\subset S$.

It is required that $\pi$ is $(J_0,j_0)$--holomorphic near $E^{\crit}$; that $E^{\crit}$ is a smooth submanifold of $E$; that the complex Hessian form
\[ H_x := \frac{1}{2}(D^2\pi)_x \colon N_x \otimes_{\C} N_x \to T_{\pi(x)}S\]
is non-degenerate as complex bilinear form, for each fibre $N_x$ of the normal bundle $N \to E^{\crit}$; and that $\Omega$ is non-degenerate on the vertical tangent distribution $\Tv E = \ker (D\pi)$, compatible with $J_0$, and `normally $J_0$--K\"ahler' (defined momentarily) near $E^{\crit}$.

{\bf (b)} An {\bf elementary symplectic Morse--Bott fibration} is one where the base $S$ is a closed disc $\bar{D}(r)$, having $S^{\crit} = \{0\}$, with connected critical set $E^{\crit}$.
It has {\bf rank} $k$ when $E^{\crit}$ has real codimension $2k+2$ in $E$.

{\bf (c)} We explain the term {\bf normally K\"ahler}. This means that a neighbourhood of $E^{\crit}$ is foliated by $J_0$-complex normal slices $\{S_x\}_{x\in E^\crit}$, such that $J_0|S_x$ is integrable and $\Omega|S_x$ K\"ahler.

\end{Defn}

In our applications, $J_0$ will usually be integrable, and $\Omega$ a $J_0$-K\"ahler form near $E^{\crit}$. This implies that it is normally $J_0$-K\"ahler: for example, one can construct the leaves $S_x$ as the fibres of the tubular neighbourhood embedding of a disc bundle in $N$ induced by the K\"ahler metric.

The condition of being `normally K\"ahler' is a technical convenience. It could probably be eliminated by means of an argument to show that after perturbing $J_0$ and $\Omega$ it can always be satisfied, but we do not pursue this point.

\begin{Defn}
A {\bf locally Hamiltonian fibration} (LHF) is a triple $(E,\pi,\Omega)$, where $\pi\colon E\to S$ is a smooth fibre bundle, and $\Omega$ a closed two--form such that $\Omega|\ker(D\pi)$ is non-degenerate.\footnote{The nomenclature is not standard but there does not seem to be a generally accepted term for these objects. Some other names that have been used, such as `symplectic bundle', are not really accurate. Some authors insist on normalisation conditions, but these are irrelevant here.}
\end{Defn}
A point of crucial importance is that an LHF has a canonical symplectic connection, whose horizontal distribution $\Th E$ is the symplectic complement of $\Tv E =\ker(D\pi)$. Locally in $S$, this connection effects a reduction of structure group from the symplectic automorphism group $\aut(E_s,\Omega|E_s)$ to the Hamiltonian group $\ham(E_s,\Omega|E_s)$.

\subsection{Lefschetz fibrations} A symplectic Morse--Bott fibration with discrete critical locus $E^{\crit}\subset E$ is called a  symplectic Lefschetz fibration; these are the subject of symplectic Picard--Lefschetz theory. In \cite{Sei}, Seidel gives a complete account of the part of this this theory which is local in the base. Briefly, this goes as follows.

With an elementary Lefschetz fibration, whose smooth fibre is $(M,\omega)= (E_{r},\Omega|E_{r})$, one associates its {\bf vanishing cycle} $(L,[f])$, which is a Lagrangian submanifold $L\subset M$ together with a framing (a diffeomorphism $f\colon S^n\to L$, up to reparametrisation by orthogonal transformations). Conversely, given a framed Lagrangian sphere $(L,[f])$ in $(M,\omega)$, there is a standard elementary Lefschetz fibration which has $(L,[f])$ as vanishing cycle.

Also associated with $(L,[f])$ is the Dehn twist $\tau_{(L,[f])}\in \aut(M,\omega) $, which is a symplectomorphism determined up to Hamiltonian isotopies supported near $L$. Given an elementary Lefschetz fibration, the monodromy of the symplectic fibration $\partial E \to \partial \bar{D}(r)$ is in the same Hamiltonian-isotopy class as the Dehn twist about its vanishing cycle.

Our aim here is to explain how these ideas generalise to symplectic Morse--Bott fibrations. For the most part, the generalisation is straightforward (and was known to Seidel circa 1998---unpublished notes). There are two new points. One concerns Hamiltonian deformation invariance of the vanishing cycles, which can be efficiently handled using a lemma of Wei-Dong Ruan concerning deformations of fibred coisotropic submanifolds. The second point, which requires some work, is how to deform $\Omega$ to a two--form which is of a standard kind in some tubular neighbourhood of $E^{\crit}$. Such an isotopy is an essential step in computing the symplectic monodromy of $(E,\pi,\Omega)$. The strategy---making a preliminary deformation so as to make $\Omega$ manageable along $E^{\crit}$ itself---is my implementation of a suggestion of Paul Seidel, and I am grateful for his advice. 

{\bf Notation.} We collect here our customary notation concerning (elementary) symplectic Morse--Bott fibrations $(E^{2n+2},\pi,\Omega,J_0,j_0)$. We write
\begin{itemize}
\item $M$ for the smooth fibre $\pi^{-1}(r)$, $\bar{M}$ for the critical manifold $E^\crit \subset \pi^{-1}(0)$;
\item $\omega$ (resp. $\bar{\omega}$) the restriction of $\Omega$ to $M$ (resp. $\bar{M}$);
\item $N$ for the normal bundle $N_{\bar{M}/E}\to \bar{M}$;
\item $V\subset M$ for the coisotropic vanishing cycle; this comes with a submersion $\rho\colon V\to \bar{M}$;
\item $\widehat{V}\subset (M,-\omega)\times (\bar{M},\bar{\omega})$ for the graph of $\rho$ (the Lagrangian vanishing cycle).
\item $\Tv E\subset TE$ is the vertical tangent distribution $\ker(D\pi)$; $\Th E$ is the $\Omega$-horizontal tangent distribution (defined on $E\setminus E^\crit$).
\end{itemize}

\subsection{Preliminaries}\label{prelim}
\subsubsection{Tubular neighbourhoods}
Let $(E^{2n+2},\pi,\Omega,J_0,j_0)$ be a symplectic Morse--Bott fibration with critical locus $\bar{M}: = E^{\crit}$ and rank $k$.
To simplify the notation, we will assume in the present discussion that there is a single critical value $s\in S$.
Choose a holomorphic chart $\xi\colon (D(r),0) \to (S,s)$.

The normal bundle $N\to \bar{M}$ is a complex vector bundle of rank $k+1$ (it carries a complex structure $J_0^N$ obtained by
linearising $J_0$). It has a natural non-degenerate complex quadratic form, the Hessian form of $\pi$, so its structure group
is reduced to $\mathrm{O}(k+1,\C)$.

\begin{Defn}
A {\bf Morse--Bott tubular neighbourhood} for $\bar{M}$ is a smooth tubular neighbourhood embedding
$\iota\colon D_\epsilon N_{\bar{M}/E} \to E$ such that, for each $x\in \bar{M}$,
(i) the map $\iota_x\colon D_\epsilon N_x \to E$ is $(J_0^N, J_0)$--holomorphic;
(ii) $\iota_x^*\Omega$ is a K\"ahler form; and
(iii) $D\xi^{-1}\circ \pi\circ \iota \colon D_\epsilon N \to \C$ is equal to the Hessian form $H_x$.
\end{Defn}

\begin{Lem}
A Morse--Bott tubular neighbourhood always exists.
\end{Lem}

\begin{pf}
By assumption, there is a family of holomorphic normal slices $S_x\subset E$ through $x\in \bar{M}$,
foliating a neighbourhood of $\bar{M}$, such that $\Omega|S_x$ is K\"ahler. Shrinking the $S_x$ if necessary,
we may suppose that they define a locally trivial fibration in which the fibres are complex manifolds
biholomorphic to a ball $B^{2n}(0;\epsilon)\subset \C^n$.

The holomorphic Morse lemma implies that, for fixed $x\in \bar{M}$, there is a holomorphic embedding
$\iota_x' \colon D_\epsilon N_x \to S_x $, mapping $0$ to $0$, such that $D\xi^{-1}\circ \pi\circ \iota' = H_x$.
Moreover, the parametric (or Morse--Bott) version of the lemma says that we can find a smooth family of such
embeddings as $x$ ranges over a ball $B\subset \bar{M}$.

Fixing $x$ again, any other such embedding differs from $\iota'_x$ by an element of $\mathrm{O}(H_x)$,
since the only holomorphic automorphisms of the ball which preserve a non-degenerate quadratic form are its orthogonal transformations. We can pin down $\iota_x$ \emph{uniquely} by saying that the derivative $D_0 \iota_x\colon N_x \to T_x S_x = N_x $ should be the identity map on $N_x$
(this is in any case part of the definition of a tubular neighbourhood embedding).
The maps $\iota_x$ then depend smoothly on $x$, for over the ball $B\subset \bar{M}$ they evidently differ from the initial choice $x\mapsto \iota_x'$ by a smooth gauge transformation. \end{pf}

\begin{Rk}
The normal bundle $p \colon N\to\bar{M}$ has a totally real subbundle
$N_{\R} = \{ v \in N: (D^2\pi)_{p(v)} (v,v) \in \R\}$ which, like $N$ itself, carries a non-degenerate quadratic form.
Thus the structure group of $N$ is reduced to $\mathrm{O}(k+1) \subset \mathrm{O}(k+1,\C)$: there is a
principal $\mathrm{O}(k+1)$-bundle $P\to \bar{M}$ (the orthonormal frames of $N_{\R}$) and an isomorphism
$P \times_{\mathrm{O}(k+1)}\C^{k+1} \cong N$.
\end{Rk}

\subsubsection{Symplectic associated bundles}
One way to construct associated bundles in the symplectic category is the following (compare e.g. Guillemin--Sternberg \cite[example 2.3]{GS} or the useful discussion in Seidel--Smith \cite[section 4.3]{SS}). Take a Hamiltonian action of a compact Lie group $G$ (with Lie algebra $ \mathfrak{g}$) on the
symplectic manifold $(F,\zeta)$, generated by the moment map $\mu\colon F \to \mathfrak{g}^*$,
so $d\mu(\xi) = -\iota(X_\xi)\zeta$ for $\xi \in \mathfrak{g}$. Suppose $p\colon P\to B$ is a principal $G$-bundle over a
smooth manifold $B$. Choose a connection one--form $\alpha \in \Omega^1(P;\mathfrak{g})$.
Pulling back $\alpha \in \Omega^1(P;\mathfrak{g})$ and $\mu \in \Omega^0(F;\mathfrak{g}^*)$ to the product $P \times F$ and contracting
via the pairing $\mathfrak{g}^* \otimes \mathfrak{g} \to \R$ results in an ordinary one--form
$ \langle \mu ,\alpha \rangle\in\Omega^1(P)$. The two--form
\[    d  \langle \mu ,\alpha \rangle + \zeta \in \Omega^2(P\times F) \]
is invariant under the action of $G$ on $P\times F$ given by $g\cdot (x ,z) = (x \cdot g^{-1}, g\cdot z)$,
so it descends to the associated bundle $P\times_G F$. It is closed, and non-degenerate on the fibres of $P\times_G F \to B$.
If $\bar{\omega}$ is a symplectic form on $B$ then
\begin{equation}\label{assoc symp form}
\Omega:= p ^* \bar{\omega} + d  \langle \mu ,\alpha \rangle + \zeta \in \Omega^2(P\times_G F)
 \end{equation}
is symplectic in an open set $U_\epsilon = P \times_G \mu^{-1}(B_\epsilon) \subset P\times_G F$ when
$B_\epsilon \subset \mathfrak{g}^*$ is a ball of sufficiently small radius $\epsilon$.
This becomes apparent as soon as one writes down the value of $\Omega$ on a pair of tangent vectors:
take $(u_1,u_2)\in T_{p(x)} B $ and lift them to $\alpha$-horizontal vectors $(u_1^{\natural}, u_2^{\natural})\in T_p P$.
Choose vertical vectors $(v_1,v_2)\in T_z F$. Then the vectors $u_i^\natural + v_i \in T_{[x,z]} (P\times_G F)$ satisfy
\[ \Omega( u_1^\natural+  v_1, u_2^\natural + v_2  )  =  \bar{\omega}(u_1,u_2)+    \langle  \mu(z) ,
	d\alpha(u_1^\natural,u_2^\natural)  \rangle +  \zeta (v_1,v_2)  . \]
We call $\Omega\in \Omega^2(U_\epsilon)$ an {\bf associated symplectic form}.

A case to keep in mind is $F=\C^{k+1}$, with $G=\mathrm{O}(k+1)$ acting linearly. The moment map is $\mu\colon \C^{k+1}\to \mathfrak{o}(k+1)^*,\; \mu(x)=  (\xi \mapsto \frac{1}{2}( x, \xi x ))$, so $\mu^{-1}(0)=\{0\}$. Its associated symplectic forms
\begin{equation}\label{assoc forms}
\Omega =  p ^* \bar{\omega} + d  \langle \mu ,\alpha \rangle + \omega_{\C^{n+1}}.
 \end{equation}
will appear periodically in this paper.

\subsection{Vanishing cycles}\label{vc}
In a symplectic Lefschetz fibration $(E,\pi,\Omega)$, one associates with a path $\gamma\colon [a,b]\to S$ leading
to a critical value  its vanishing cycle. This is a Lagrangian sphere in $E_{\gamma(a)}$. In a symplectic Morse--Bott fibration,
the vanishing cycle is rather a \emph{Lagrangian correspondence} between $E_{\gamma(a)}$ and the singular locus in $E_{\gamma(b)}$.

\subsubsection{Symplectic parallel transport} As observed above, an LHF $(E,\pi,\Omega)$ has a natural connection:
the horizontal subspace $\Th _x E \subset T_x E$ is defined to be the set of vectors $u$ such that $\iota(u)\Omega$ is zero on
$\Tv_x E $. The connection is symplectic is the sense that any horizontal vector field $h$ satisfies $\mathcal{L}_{h}\Omega  =0 $.
When $\pi$ is proper, the connection can be integrated over any smooth path $\gamma\colon [a,b]\to S$ in the base, so $\gamma$ has a
parallel transport map
\[ \rho_\gamma \colon E_a \to E_b, \]
a symplectomorphism between the fibres.

Now suppose that $(E,\pi,\Omega)$ is symplectic Morse--Bott. Parallel transport obviously still make sense for paths in the open
set of regular values $S^{\mathrm{reg}}\subset S$. If $\gamma\colon [a,b]\to S$ is a path satisfying
$\gamma^{-1}(S^{\crit})=\{b\}$, define
\begin{equation}V_\gamma = \{ x \in E_a : \lim_{ t\to b^-} {\rho_{\gamma|[a,t]}(x)} \in E^{\crit}  \},  \end{equation}
the set of points for which the limiting parallel transport exists and lands in the singular locus of $E_b$.
Put $M= E_a$ and $\bar{M} = E^{\crit}\cap E_b$, and denote by $\rho$ the limiting parallel transport map $V_\gamma \to \bar{M}$.
\begin{Lem}
$V_\gamma\subset M$ is a submanifold and the map $\rho\colon V_\gamma\to \bar{M}$ a smooth fibre bundle.
The fibres are spheres $S^k$, where $k$ is the rank of $(E,\pi)$.
The structure group of $\rho$ is reduced, in a canonical way, to $\mathrm{O}(k+1)$,
via an isomorphism of $V_\gamma$ with the unit sphere bundle in $N_{\R}$.
\end{Lem}
\begin{pf}
See \cite[Lemma 1.13]{Sei}. The only difference is that we must use Morse--Bott tubular neighbourhoods instead of holomorphic Morse charts.
\end{pf}
Notice that there is even a well-defined parallel transport map $\rho_\gamma\colon E_a \to E_b$, the pointwise limit of
$\rho_\gamma|{[a,b']}$ as $b'\to b $ from below. Consequently, an elementary symplectic Morse--Bott fibration deformation-retracts to its critical fibre.

The fibre bundle $\rho\colon V_\gamma \to \bar{M}$, together with its reduction to $\mathrm{O}(k+1)$ and embedding $V_\gamma \to M$ is called the {\bf vanishing cycle} associated with $\gamma$.

The restriction of the $\Omega$ to $V_\gamma$ is also the pullback of $\Omega$ from the critical set:
\begin{equation} \rho^* (\Omega|\bar{M}) = \Omega|V_\gamma .\end{equation}
This follows readily from the fact that the parallel transport maps $\rho_{\gamma|[a,t]}$ are symplectic.

The symplectic complement of $T V_\gamma$ is $\ker(D \rho_\gamma)\subset TV_\gamma$. The fibres of $\rho$ are isotropic spheres.
Thus $V_\gamma$ is a {\bf fibred coisotropic submanifold} of $M$: $TV_\gamma$ contains its own symplectic complement, and the isotropic
foliation of $V_\gamma$ is a fibration: each leaf $F$ has a neighbourhood diffeomorphic to $F\times B^{2n-k}$ by a
diffeomorphism which takes the isotropic foliation to the product foliation with leaves $F \times \{z\}$.

\subsubsection{Good two--forms}
It is sometimes convenient to have at one's disposal a space of two--forms, not necessarily closed, which have well-defined vanishing cycles in the smooth (not symplectic) category. We therefore make the following \emph{ad hoc} definition.

Consider a symplectic Morse--Bott fibration $(E,\pi,\Omega_0,J_0,j_0)$.
\begin{Defn} \label{good forms}
Fix a Morse--Bott tubular neighbourhood $N$ of $E^\crit$, and consider it as an associated $\C^{k+1}$-bundle of a principal $\mathrm{O}(k+1)$-bundle $p\colon P\to E^\crit$. Fix also an almost complex structure $J$ on $E$, extending $J_0$, such that $ D\pi \circ J= \ii \circ  D\pi $, tamed by $\Omega_0$. A two--form $\Omega$ is {\bf good} if (a) it is non-degenerate on $\Tv E$ (and so defines a connection away from $E^\crit$); (b)
it is an `associated form' $p ^* \bar{\omega} + d  \langle \mu,\alpha \rangle + \omega_{\C^{k+1}}$ on $N$, as in formula (\ref{assoc forms}) above; and (c) it `tames' $J$, i.e. $\Omega(u,Ju)>0$ when $u\neq 0$.
\end{Defn}
So $\alpha$ is a connection one--form on $P$; $\mu$ the moment map for $\mathrm{O}(k+1)$ acting on $\C^{k+1}$;
and $\bar{\omega}$ a non-degenerate (but not necessarily closed) two--form on $E^\crit$.

Parallel transport is well-behaved for a good two--form $\Omega$ (the only potential problem is near the critical set, and there it is the closed forms $ d  \langle \mu,\alpha \rangle+ \omega_{\C^{k+1}}$ which control the transport, so everything works as usual). It therefore defines a vanishing cycle which is an $S^k$-bundle over $E^\crit$.

\begin{Lem}
Once the tubular neighbourhood and $J$ are fixed, (i) the space of good forms is contractible, hence any two vanishing cycles are smoothly isotopic; (ii) if one has a good form $\Omega$ defined on an open subset $U$ of $E$ containing $N$, and if $U' \subset\subset  U$ is an open set whose closure lies within $U$, then there is a  globally-defined good form which agrees with $\Omega$ on $U'$.
\end{Lem}
\begin{pf}
The space of good forms is convex, which gives (i). For (ii), we can patch locally-defined forms; the patching works because of the taming condition.
\end{pf}

\subsubsection{Hamiltonian deformations}
A solution to the infinitesimal deformation problem for fibred coisotropic manifolds has been given by W. Ruan \cite{Rua}.

Let $(M,\omega)$ be a symplectic manifold, and $V_0\subset M$ a fibred coisotropic submanifold.
Let $(\bar{M}_0,\bar{\omega}_0)$ be the reduced space of isotropic leaves, and $\rho_0\colon V_0\to \bar{M}_0$ the quotient map.

Given an isotopy $\{ V_t\}_{t\in [0,1]} $, one chooses diffeomorphisms $\phi_t \colon V_0\to V_t$ covering
diffeomorphisms $\bar{\phi}_t\colon \bar{M}_0\to \bar{M}_t$ of the reduced spaces. Let $X_t = \dot{\phi}_t$ be the
generating vector fields on $V_t$, and put $\beta_t = \iota(X_t)(\omega|V_t)$. Then for any vector field $Y_t$
on $V_t$ tangent to the isotropic distribution, one has $\iota(Y_t)\mathcal{L}_{X_t}(\omega|V_t)=0$. But
$\iota(Y_t)\mathcal{L}_{X_t}(\omega|V_t)=\iota(Y_t)d \iota(X_t) (\omega|V_t) = \iota(Y_t) d\beta_t$.
Thus both $\iota_{Y_t}d\beta_t$ and $\mathcal{L}_{Y_t} d\beta_t = d \iota_{Y_t}d\beta_t=0$ vanish.
Hence $d\beta_t = \rho_t^* \gamma_t $ for closed two--forms $\gamma_t \in \Omega^2(\bar{M}_t)$.
It follows that $\phi_t^*\beta_t$ represents a class
\begin{equation}\label{flux class}
[\phi_t^*\beta_t]\in H^0(\bar{M}_0;\mathcal{H}^1),\end{equation}
where $\mathcal{H}^1$ is the natural local system on $\bar{M}_0$ with fibres $\mathcal{H}^1_x = H^1( \rho_0^{-1}(x) ;\R)$. Let us call the class (\ref{flux class}) the {\bf flux} of $\{(V_t, \phi_t)\}$.

Ruan shows that the isotopy $\{V_t\}$ is symplectic---that is, it is induced by a global symplectic flow on $M$---if and only
if $\delta [\phi_t^*\beta_t]  = 0$ for each $t$, where $\delta\colon H^0(\bar{M}_0;\mathcal{H}^1)\to H^2(\bar{M}_0;\R)$ is the
Leray--Serre differential. He also proves the following lemma; we take the liberty of reproducing the proof.
\begin{Lem}[Ruan \cite{Rua}]
$\{V_t\}$  is a Hamiltonian isotopy if and only if the flux $[\phi_t^*\beta_t]$ is zero for all $t$.
\end{Lem}
\begin{pf}
Necessity is clear. For sufficiency, observe that the vanishing of $[\phi_t^*\beta_t]$ means that
there are functions $K_t\in C^\infty(V_t)$ such that $dK_t - \beta_t$ vanishes on the isotropic fibres of $V_t$.
For any vector field $Y_t$ on $V_t$ which is tangent to the isotropic distribution, we have
$\iota(Y_t)(dK_t-\beta_t)=0$, and
$\mathcal{L}_{Y_t}(dK_t-\beta_t)= -\iota(Y_t)d\beta_t=0$.
So $dK_t-\beta_t = \rho_t^* \beta_t'$ for some $\beta_t'\in\Omega^1(\bar{M})$. Now, $\beta_t' = \iota(Z_t)\omega$ for a vector field $Z_t$ along $V_t$ which is tangent to $V_t$.
Let $X_t' = X_t - Z_t$. Integrating $X_t'$, one gets a new flow
$\phi_t ' \colon V_0\to V_t$. Since $\iota(X_t')\omega = dK_t$, $\phi_t'$ globalises to a Hamiltonian isotopy.
\end{pf}

We apply the lemma to the vanishing cycles of $(E,\pi,\Omega)$ for varying $\Omega$.
Let $\gamma\colon [0,1]\to S$ be a path in the base, with $\gamma^{-1}(S^\crit)=\{1\}$.
Let $M=\pi^{-1}(\gamma(0))$ be the smooth fibre, $\bar{M}\subset \pi^{-1}(\gamma(1))$ the critical manifold,
and $ V\subset M$ the coisotropic vanishing cycle associated with $\gamma$. Parallel transport defines an
$S^k$-bundle $\rho\colon V\to \bar{M}$ with isotropic fibres. One can also identify this map with the quotient
map to the reduced space of isotropic leaves in $V$.

Consider a path $\Omega_t$ of two--forms on $E$, each making $\pi$ a symplectic Morse--Bott fibration,
such that $d\Omega/dt$ is exact for all $t$.
Suppose additionally that $\Omega_s$ is constant on $M$ and on $\bar{M}$.
Let $(V_t, \rho_t)$ be the vanishing cycle associated with $\gamma$, defined via $\Omega_t$.

\begin{Lem}\label{ham isotopy}
There is a Hamiltonian isotopy $t\mapsto \Phi_t\in \ham(M,\omega)$ with $\Phi_t(V_0)=V_t$.
\end{Lem}
\begin{Rk}
The folowing proof applies for any $k$; however, when $k>1$, the result follows immediately from Ruan's lemma,
since $H^1(S^k;\R)=0$.
\end{Rk}
\begin{pf}
Consider the `thimble' $W_t \subset E$ associated with the path $\gamma$ and the form $\Omega_t$, i.e. the closure of the manifold swept out by $V_t$ via $\Omega_t$-parallel transport along $\gamma$. We may assume (by adding a form pulled back from the base) that $\Omega_t$ is symplectic on $E$. Then $W_t\subset (E_t,\Omega_t)$ is fibred coisotropic, and its isotropic fibres are $(k+1)$-disks. Choose diffeomorphisms $\Phi_t\colon W_0 \to W_t$ covering the identity map on $\bar{M}$, with generating vector field $X_t = \dot{\Phi}_t$ along $W_t$. Let $\alpha_t = \iota(X_t)\Omega_t \in \Omega^1(W_t)$. Then $\Phi_t^*\alpha_t$ has a flux $[\Phi_t^*\alpha_t]\in H^0(\bar{M}; \mathcal{H}^1)$. But the fibre $\mathcal{H}^1$ is $H^1(\bar{D}^{k+1};\R)=0$, so $\mathcal{H}^1$ is the zero local system and $[\Phi_t^*\alpha_t]=0$.

Now, $\Phi_t$ restricts to $V_0=\partial W_0$ to give a diffeomorphism $\phi_t \colon V_0 \to V_t$, generated by the vector field $Y_t=X_t|V_t$. Put $\beta_t =\iota(Y_t)\omega$. The flux $[\phi_t^*\beta_t]$ lies in $H^0(\bar{M}; \mathcal{K}^1)$, where the fibre of the local system $\mathcal{K}^{1}$ is $H^1(S^k;\R)$. It is the image of of $[\Phi_t^*\alpha_t]$ under the natural restriction map, and is therefore zero. The result now follows from Ruan's lemma.
\end{pf}

There is a useful repackaging of this result in the language of fibre bundles, which allows us to dispense with the assumption that the path $\Omega_t$ is constant on $M$ and $\bar{M}$.

Suppose that $\EuScript{M}\to [0,1]$ is a fibre bundle, and $\zeta \in \Omega^2(\EuScript{M})$
a closed, fibrewise-symplectic form. Let $\EuScript{V}\subset \EuScript{M}$ be a sub-bundle such that the fibres $\EuScript{V}_t\subset\EuScript{M}_t$ are fibred coisotropic. There is then a bundle $\bar{\EuScript{M}}\to[0,1]$ of reduced spaces, and a quotient map $\rho \colon \EuScript{V}\to \bar{\EuScript{M}}$. One has the easy
\begin{Lem}
The following two conditions are equivalent:
\begin{enumerate}
\item  $\EuScript{V}$ is coisotropic; moreover, $\zeta|\EuScript{V}= \rho^*\bar{\zeta}$ for a closed two--form $\bar{\zeta}$ on $\bar{\EuScript{M}}$;
\item  the sub-bundle $\EuScript{V}$ is preserved by the $\zeta$-parallel transport maps over intervals $[a,b]\subset [0,1]$.
\end{enumerate}
When this happens, we shall simply say that $\EuScript{V}$ is `globally coisotropic'.
\end{Lem}
The $\zeta$-parallel transport $\phi_t\colon \EuScript{M}_t\to\EuScript{M}_0$ over $[0,t]$ trivialises the symplectic fibration $\EuScript{M}\to [0,1]$. The fibre $(M,\omega)=(\EuScript{M}_0,\omega_0)$ contains coisotropic submanifolds $\phi_t(\EuScript{V}_t)$ which are easily seen to be Hamiltonian isotopic. Conversely, if $V_t\subset (M,\omega)$ are coisotropic submanifolds which are Hamiltonian isotopic, then $\bigcup_{t\in [0,1] }{V_t\times \{t\}}\subset M\times[0,1]$ is globally coisotropic with respect to a closed form of shape $\omega+d(H_t dt)$ on $M\times [0,1]$.

Choose a connection $\nabla$ on $\EuScript{V}\to [0,1]$, and define $\beta\in \Omega^1(\EuScript{V})$ by
\[ \beta = \iota(\partial_t^\natural - \widetilde{\partial}_t)\Omega.  \]
Here $\partial_t^\natural$ is the vector field of $\Omega$-horizontal lifts of $\partial_t$ along $\EuScript{V}$, and $\widetilde{\partial}_t$ the field of $\nabla$-horizontal lifts of $\partial_t$, tangent to $\EuScript{V}$. (If $\EuScript{V}$ were globally coisotropic, we could define $\nabla$ to be the connection obtained from $\Omega$ by restriction, and $\beta$ would be zero.) Then $\beta$ defines a flux
\begin{equation}\label{bundle flux} [\beta] \in H^0(\bar{\EuScript{M}},\mathcal{H}^1),\end{equation}
where $\mathcal{H}^1$ is now the natural local system on $\bar{\EuScript{M}}$ with fibres $\mathcal{H}^1_x = H^1(\rho^{-1}(x);\R)$. To see this, simply observe that $\Omega$-parallel transport trivialises $\EuScript{M}$ and reduces the situation to the one already considered.

\begin{Prop}
If $[\beta]$ vanishes, we can find a different closed two--form $\zeta'$
such that $\zeta-\zeta'$ vanishes on the fibres, for which $\EuScript{V}$ is globally coisotropic.
\end{Prop}
\begin{pf}
This is Ruan's lemma, translated into fibre bundle language.
\end{pf}
Let us apply this to families of vanishing cycles:
\begin{Lem}\label{globally coisotropic}
Suppose that $E$ is equipped with two--forms $\Omega_s$, $s\in[0,1]$, each making $(E,\pi,\Omega_s)$ a symplectic Morse--Bott fibration, and that the $\Omega_s$ are the restrictions of a closed two--form $\zeta$ on $E\times[0,1]$. Then the forms $\omega_s=\Omega|E_1$ are the restrictions of a closed two--form $\eta\in\Omega^2(M\times [0,1])$ such that the union of the coisotropic vanishing cycles,
\[ \EuScript{V} = \bigcup_{s\in[0,1]}{ V_s \times\{s\} } \subset M\times [0,1] \]
is globally coisotropic. Moreover, we can take $\eta$ to equal $\zeta| M\times [0,1]$ outside a small neighbourhood of $\EuScript{V}$.
\end{Lem}
Indeed, our previous argument shows that the flux vanishes.
\begin{Rk}
We can also allow the base to be $S^1$ rather than $[0,1]$.
Any LHF $(\EuScript{M} \to S^1, \Omega)$ is isomorphic to the \emph{mapping torus} of its monodromy.
A mapping torus is a bundle $\torus(\phi)$ obtained from a symplectic automorphism $\phi\in \aut(M,\omega)$:
\[\torus(\phi) = (M\times[0,1]) /(\phi (x) ,0)\sim (x,1).\]
This space maps naturally to $S^1$.  The two--form is $\omega_\phi$, the unique form whose pullback to $M\times[0,1]$ is the pullback of $\omega$ from $M$. If $\EuScript{V}$ is a fibrewise-integral-coisotropic subbundle whose flux vanishes, we can replace $\omega_\phi$ by a form of shape $\omega_\phi+ d(t\wedge H_t)$ which makes $\EuScript{V}$ globally coisotropic. Notice that the replacement form is cohomologous to $\omega_\phi$.
\end{Rk}

\subsection{Fibred Dehn twists}\label{twists}
\emph{Neither this subsection, nor the next, is required for the definition of Lagrangian matching invariants.}

Seidel explains in \cite[section 1.2]{Sei} that the symplectic manifold $(T^* S^n, d\lambda_{\mathrm{can}})$ carries a distinguished class of compactly supported symplectomorphisms, the \emph{model Dehn twists}. Define
\begin{equation}  \mathcal{H} = \{ h \in C^\infty(\R,\R): h\text{ is even and}\;  h(t) = - |t|/2\text{ for } t\gg 0  \}.  \end{equation}
The model Dehn twist $\delta_h$ is an automorphism associated with an element $h\in \mathcal{H}$. The map $\mathcal{H}\to \aut_c(T^*S^n,d\lambda_{\mathrm{can}}) $, $h\mapsto \delta_h$, is constructed as follows.

There is a Hamiltonian $S^1$--action on $T^*S^n\setminus S^n$, generated by the moment map $\mu ( v,x ) = |v|$, where $|\cdot | $ is the norm inherited from $\R^{n+1}$ on $T^*S^n = \{  (v,x) \in \R^{n+1}\times\R^{n+1}: |x| =1, \,\langle v,  x\rangle  = 0 \}$. The Hamiltonian functions $H_1 = h\circ \mu$, $H_2 = \frac{1}{2}\mu$ Poisson--commute; hence $\phi_t^{H_1+ H_2} = \phi_t^{H_1}\circ \phi_t^{H_2}$. The model Dehn twist $\delta_h$ is defined (on $T^*S^n\setminus S^n$) to be the time-$2\pi$ Hamiltonian flow $\phi_{2\pi}^{H_1+ H_2} = \phi_{2\pi}^{H_1}\circ \phi_{2\pi}^{H_2}$. Because $h$ is a smooth function of $t^2$, $H_1$ extends smoothly over the zero-section $S^n$, and its flow is trivial there; $\phi^{H_2}_{2\pi}$ also extends smoothly over the zero-section, where it acts as the antipodal map. Thus $\delta_h \in \aut_c(T^*S^n,d\lambda_{\mathrm{can}}) $, and $\delta_h (0,x) = (0,-x)$. Since $\mathcal{H}$ is is convex, $\delta_{h'}$ differs from $\delta_h$ by a canonical, compactly supported Hamiltonian isotopy.

Now consider a principal $\mathrm{O}(n+1)$-bundle $p\colon P\to \bar{M}$. The linear action of $\mathrm{O}(n+1)$ on $S^n$ induces an action on $T^*S^n$, which is Hamiltonian with moment map $\nu$. We can form the associated bundle $T = P\times_{\mathrm{O}(n+1)}T^*S^n$, with associated closed forms,
\[   d \big( \lambda_{\mathrm{can}}  + \langle \nu ,\alpha \rangle\big) \in \Omega^2(T), \]
defined via connection one--forms $\alpha \in \Omega^1(P;\mathfrak{o}_{n+1})$. If $\bar{\omega}\in \Omega^2(\bar{M})$ is symplectic then $\omega:= p^*\bar{\omega} +  t d ( \lambda_{\mathrm{can}}  + \langle \nu ,\alpha \rangle) \in \Omega^2(T)$ is symplectic on a neighbourhood of the zero-section $U= P\times_{\mathrm{O}(n+1)}T^*S^n_{\leq \lambda}$.

The moment map $\mu\in C^\infty( T^*S^n \setminus S^n )$ globalises to a map $\mu_P\in C^\infty (T \setminus \bar{M})$. Given $h\in\mathcal{H}$, we can form $H_1 = h\circ \mu_P$, as well as $H_2 = \mu_P/2$. Let $\mathcal{H}_U= \{h\in \mathcal{H}: \mathrm{supp}(H_1)\subset U\}$. Given $h\in \mathcal{H}_U$, we can consider the $\omega$-Hamiltonian flow $\delta_h^U:=\phi_{2\pi}^{H_1+H_2}$ on $U$. The vertical vector field $X_{H_i}$, given on each fibre as the $d\lambda_{\mathrm{can}}$-Hamiltonian vector field of $H_i$, is easily checked to be the \emph{global} $\omega$-Hamiltonian vector field of $H_i$. Hence the flow $\phi_{2\pi}^{H_1+H_2}$ preserves the fibres of $p$. Thus $\delta_h^U$ is just the map obtained by applying $\delta_h$ to each fibre, and as such it is globally smooth, with proper support inside $U$.

To sum up, we have:
\begin{Lem}
There is a map $\mathcal{H}_U \to \aut(U, \omega)$, $h\mapsto \delta_h^U$, such that (i) $\delta_h^U$ covers the identity map on $\bar{M}$; (ii) $\delta_h$ has proper support inside $U$; (iii) $\delta_h^U$ acts on the zero section $P\times_{\mathrm{O}(n+1)}S^n$ as the antipodal map $[p,x]\mapsto [p,-x]$; and (iv) $\delta^U_{h'}$ is canonically Hamiltonian isotopic to $\delta_{h}^U$.
\end{Lem}
We call $\delta_h^U$ the {\bf model fibred Dehn twist} associated with $h$ on the associated bundle $P\times_{\mathrm{O}(n+1)}T^*S^n$.

{\bf Model fibred Dehn twists as monodromy maps.} Model fibred Dehn twists arise as monodromy maps for certain symplectic Morse--Bott fibrations over the disc. Actually this cannot \emph{literally} be true, because according to our definition, the fibres of a symplectic Morse--Bott fibration are closed manifolds. However, there is an obvious extension of the definition to allow fibres which have boundary (similar to that of Lefschetz fibrations in \cite{Sei}) where the total space has a codimension-two corner where the `horizontal' and `vertical' parts of the boundary meet.

A basic example of a Lefschetz fibration $(E_r,\pi,\omega_{\C^{n+1}},J_0,j_0)$ in this broader sense is the following from \cite[Lemma 1.10]{Sei}: let
\begin{equation}\label{model fibration}
 E_r = \{ z\in  \C^{n+1}: |q(z)| \leq r,\, \| z\|^4-|q(z)|^2\leq 4\lambda^2\},\quad q(z)=\sum_i{z_i^2}.  \end{equation}
This maps to the disc $\bar{D}(r)$ by $\pi:=q|E_r$. The complex structures are the standard ones.

\begin{Lem}[Seidel]
There is an $\mathrm{O}(n+1)$-invariant one--form $\alpha$, supported in $E_r\cap \{ z: \| z\|^4-|q(z)|^2\geq 3 \lambda^2\}$ such that the following holds. There is a canonical isomorphism $\phi\colon (E_r,\omega_{\C^{n+1}}+ d\alpha ) \to ( (T^*S^n)_{\leq \lambda}, d\lambda_{\mathrm{can}})$, and the (positive) monodromy $\rho$ of $(E_r, \omega_{\C^{n+1}}+ d\alpha)$ around $\partial \bar{D}(r) $ has the property that $ \phi\circ \rho \circ \phi^{-1 }$ is a model Dehn twist.
\end{Lem}

The orthogonal group $\mathrm{O}(n+1)$ preserves $E\subset \C^{n+1}$ and preserves the fibres of $q$,
so we can form an associated bundle $E_P := P\times_{\mathrm{O}(n+1)} E\to \bar{M}$ and a map $q_P\colon E_P\to \bar{D}(r)$,
$[p,z]\mapsto q(z)$. A choice of connection $\alpha$ on $P$ gives an associated two--form
\[ \Omega: =  p^*\bar{\omega} + \zeta + d\langle \alpha,\nu \rangle,  \]
which, if we choose $\lambda$ and $r$ small, is symplectic on $E_P$. The fibre $(E_P)_r=q_P^{-1}(r)$ is also a bundle over $\bar{M}$. Moreover, the map $\phi\colon E_r\to (T^*S^n)_{\leq \lambda} $ is $\mathrm{O}(n+1)$-equivariant,\footnote{This follows from the interpretation of $\phi$ as a symplectic parallel transport map, \cite[Lemma 1.10]{Sei}.} hence $\phi_P\colon (E_P)_r \to P\times_{\mathrm{O}(n+1)}T^*S^n_{\leq \lambda}$, $[p,z]\mapsto [p,\phi(z)]$ is well-defined and pulls back
$\omega = p^*\bar{\omega}+ d\langle \nu,\alpha \rangle + d\lambda_{\mathrm{can}}$ to  the restriction of $\Omega$.

Thus we obtain the
\begin{Lem}
There is a canonical isomorphism
\[\phi _P \colon ((E_P)_r,\Omega|(E_P)_r) \to ( P\times_{\mathrm{O}(n+1)} T^*S^n_{\leq \lambda}, d\lambda_{\mathrm{can}}).\]
Moreover, the monodromy $\rho_P$ of $E_P$ around $\partial \bar{D}(r) $, taken in the positive direction,
has the property that $ \phi_P\circ \rho_P \circ \phi^{-1 }_P$ is a model fibred Dehn twist.
\end{Lem}

\subsubsection{Fibred Dehn twist along a coisotropic sphere--bundle}
The model fibred Dehn twists, and the standard fibrations $E_P$, can be transplanted into other manifolds.

Suppose given
\begin{itemize}
\item
symplectic manifolds $(M^{2n},\omega)$ and $(\bar{M}^{2n-2k}, \bar{\omega})$;
\item
an orthogonal $S^k$-bundle $p\colon V \to \bar{M}$;
\item
a smooth, proper embedding $e\colon V\to M$ such that $e^*\omega=p^* \bar{\omega}$.
\end{itemize}
From these data, one immediately obtains a Lagrangian correspondence---the graph of $p$:
\begin{equation}\label{V-hat}
\widehat{V} := \{(x, px ): x\in V \} \subset (M , - \omega) \times (\bar{M} , \bar{\omega}). \end{equation}
Moreover, one can construct (a) an automorphism $\tau_V\in \aut(M,\omega)$, the {\bf fibred Dehn twist along $V$}, supported in a tubular neighbourhood of $V$, and determined up to Hamiltonian isotopies also supported in a tubular neighbourhood; and (b) a symplectic Morse--Bott fibration $(E_V,\pi,\Omega)$ over the disc $\bar{D}(r)$ such that $\pi^{-1}(r)\cong (M,\omega)$. The monodromy of $E_{V}$ around $\partial \bar{D}(r)$  is a fibred Dehn twist along $V$.

The construction for (a) is as follows. The submanifold $e(V)\subset M$ is coisotropic, for the annihilator of $e_*(T_{x} V) $ in $T_{e(x)}M$ is $\ker(D_x p)\subset e_*(T_x V)$ (the inclusion $\ker(D_x p))\subset e_*(T_{x} V) ^\omega$ is clear, and equality holds by a dimension-count).

Let $T_\lambda  = P_V \times_{\mathrm{O}(k+1)} (T^* S^k)_{\leq \lambda} $, where $p_V\colon P_V\to \bar{M}$ is the principal $\mathrm{O}(k+1)$-bundle of orthogonal frames. The zero-section gives an embedding $e'\colon V\to T_\lambda$, and with the symplectic form $\omega' = p_V^*\bar{\omega}+ d\langle \mu,\alpha\rangle + d\lambda_{\mathrm{can}}$, this is a coisotropic embedding with $e'^*\omega = e^*\omega$.

The coisotropic neighbourhood theorem tells us that, near $e(V)$, $\omega$ is determined up to symplectomorphism by $e^*\omega$. Hence $e$ extends to a symplectic embedding $\hat{e}\colon T_\lambda \to M$ for some $\lambda>0$ (unique up to Hamiltonian isotopies acting trivially on $e(V)$).

We define $\tau_V$ by $\tau_V (\hat{e}(x) ) = \hat{e} (\delta^{P_V}_h (x))$ for a function $h\in \mathcal{H}$. Since this is supported inside $\im(\hat{e})$, we may extend it trivially over $M$. Though $\tau_V$ depends on the choice of embedding $\hat{e}$, it  its Hamiltonian isotopy class does not, and we shall imprecisely refer to `the Dehn twist $\tau_V$'.

For (b), we again make use of the standardising embedding $\hat{e}$. Let $M_0$ be the closure of $M\setminus \im(\hat{e})$. We define
\[ E_V =  \big( M_0 \times \bar{D}(r)  \big) \cup_\Phi  E_P \]
where the gluing map $\Phi$ identifies the `horizontal boundaries' of the two pieces. In doing so it must marry both the projections to $\bar{D}(r)$ and the two--forms. (On $M_0\times \bar{D}(r)$ we use the trivial projection and the two--form which pulls back $\omega$.) A suitable map $\Phi$ is constructed by Seidel \cite[Lemma 1.10]{Sei} in the case $\bar{M}=\{\mathrm{pt.} \}$. It is $\mathrm{O}(k+1)$-equivariant, and can therefore be applied to associated bundles.

\subsection{Monodromy}\label{monodromy}
We now bring together the two threads of the discussion---vanishing cycles and fibred Dehn twists---by proving the following.

\begin{MonThm}\label{monod}
Let $(E,\pi,\Omega)$ be an elementary symplectic Morse--Bott fibration over $\bar{D}(r)$, with smooth fibre $M:=E_r$,
critical set $\bar{M} = E^{\crit}$, and vanishing cycle $\rho \colon V\to \bar{M}$.
Then the monodromy $\rho_E\in \aut(M, \Omega|M)$ is Hamiltonian--isotopic to the fibred Dehn twist $\tau_V$.
\end{MonThm}

Here we have to do a little more work in generalising from symplectic Lefschetz fibrations to symplectic Morse--Bott fibrations
(though the existence of Morse--Bott tubular neighbourhoods is already a useful preliminary step).

The theorem is a consequence of the following technical result.
\begin{Prop}\label{ass symp forms}
Let $(E,\pi,\Omega,J_0,j_0)$ be a symplectic Morse--Bott fibration with critical manifold $\bar{M}=E^\crit$, equipped with
holomorphic charts near $S^{\crit}$ and Morse--Bott tubular neighbourhood $ \tau \colon D_\epsilon N \to E$.
Identify $N$ with the associated bundle $P\times_{\mathrm{O}(k+1)}\C^{k+1}$, where $P\to \bar{M}$ is the
principal $\mathrm{O}(k+1)$-bundle of orthonormal frames of $N_{\R}$.
Then there is a family of two--forms $\{\Omega_t\}_{t\in [0,1]}$, such that
\begin{itemize}
\item $\Omega_0=\Omega$;
\item $\Omega_t$ tames $J_0$ for each $t$;
\item there exist one--forms $\alpha_t$ such that $d\Omega_t/dt = d\alpha_t$, with $\alpha_t | (E\setminus \im(\tau)) =0$ and $\iota^* \alpha_t =0$; and
\item $\tau^*\Omega_1$  is an associated symplectic form in a neighbourhood of the zero-section.
\end{itemize}
\end{Prop}

\begin{pf}[Proof of Theorem \ref{monod}]
When we deform $\Omega$ as in Prop. \ref{ass symp forms}, the monodromy around $\partial \bar{D}(0;r)$ does not change,
since $\Omega_t$ is constant over $\partial \bar{D}(0;r)$. The vanishing cycle at time $t$ is a coisotropic submanifold
\[  V_t \subset (M,\omega).   \]
By Lemma \ref{ham isotopy}, the isotopy of submanifolds $\{V_t\}$ is in fact Hamiltonian: it is generated by a path
$\phi_t \in \ham(M,\omega)$. Since $\tau_{\phi_t V_0}=\phi_t \circ \tau_{V_0} \circ \phi_t^{-1}$, the Hamiltonian isotopy
class of the associated Dehn twist is constant. We may therefore assume, by Prop. \ref{ass symp forms}, that $\bar{M}$
has a Morse--Bott tubular neighbourhood $N$ in which $\Omega$ is an associated symplectic form.
The argument is then virtually the same as in \cite[Proposition 1.15]{Sei}, so we shall be very brief.
It suffices, by an argument involving radial parallel transport over an annulus, to show that the monodromy of a
very small loop $\partial \bar{D}(\epsilon r)$ is Hamiltonian isotopic to the Dehn twist.
Over $\bar{D}(\epsilon r)$, one can identify the fibration with the union of a standard piece contained in $N$,
and a trivial piece (trivialised by symplectic parallel transport). The result then follows from the definition of the fibred Dehn twist.
\end{pf}

Proposition \ref{ass symp forms} will be deduced from a lemma:
\begin{Lem}\label{technical lemma}
Let $\bar{M}$ be a compact manifold, $F \to \bar{M}$ a real vector bundle of rank $r$ with Euclidean metric $g$, and $J$ an almost complex structure on the total space of $F\otimes \C$ such that (i) $J$ acts as scalar multiplication by $\ii$ on the fibres $F_x \otimes \C$, and (ii) the image of the zero-section $\iota\colon \bar{M}\to F\otimes \C$ is an almost complex submanifold. Let $\Omega$ be a symplectic form on the disc-subbundle $U = \{v\in F\otimes \C : g_{\C}(v,v) < R \} \subset F\otimes\C$, compatible with $J$.

Then there is another symplectic form $\Omega'$ on $U$, still taming $J$, equal to $\Omega$ near $\partial U$ and satisfying $\iota^* \Omega'=\iota^* \Omega$, but also invariant under unitary gauge transformations along $\im(\iota)$.
\end{Lem}
Here `unitary' is taken with respect to the hermitian metric $g_{\C}$, the hermitian extension of $g$ to the complexified bundle. Gauge invariance means that there is a constant $t\in \R$ such that, for any $x\in \bar{M}$, if $u_1,u_2\in \Th_{\iota(x)}(F\otimes\C) $,  $v_1,v_2\in \Tv_{\iota(x)}(F\otimes\C) $, then $\Omega'_{\iota(x)}(u_1 + v_1, u_2+v_2)= \Omega'_{\iota(x)}(u_1,u_2)+ t \imag g_{\C}(v_1,  v_2)$.
\begin{pf}
The almost complex structure $J$ gives rise to a $d^c$-operator on forms, $d^c = J \circ d$. We can write $\Omega$ as
\begin{equation}\label{decomp Omega}
 \Omega =  p^*\iota^*\Omega + d\alpha
 \end{equation}
where $p$ is the projection $F\otimes \C \to \bar{M}$. Moreover, we can write $\alpha = d^c \phi +   \beta$, where $\phi(z) = - h(z,z)$ for some hermitian metric $h$; $\iota^*d \beta=0$; and $d\beta| (F_x\otimes\C)$ vanishes at the origin. By adding a closed one--form to $\beta$, we may suppose that $\iota^*\beta=0$ as well.

We shall deal with the terms one by one. First, take a smooth, increasing  function $\chi \colon\R_{\geq 0}\to \R$, identically $0$ on $[0,1]$ and $1$ on $[2,\infty)$. Define $\chi_\delta \colon F\otimes\C \to \R $ for $\delta>0 $ by $\chi_\delta(v) = \chi ( | v |  /\delta ) $, where $| \cdot |$ is the hermitian metric obtained by complexifying $g$. Let
\[ \Omega_1=  p^*\iota^*\Omega  +  dd^c \phi  + d(\chi_\delta\beta), \]
so $\Omega_1$ is a closed two--form which differs from $\Omega$ only near the zero--section. We claim that $ \Omega_1$ tames $J$, provided that $\delta$ is small. Along the zero--section, $J$ is tamed by $p^*\iota^*\Omega+ dd^c\phi$; the same is therefore true inside a disc-bundle of some small radius $\delta'$. Take $\delta\leq \delta'/2$. By convexity, $J$ is tamed by the non-closed form $p^*\iota^*\Omega  +  dd^c \phi  + \chi_\delta \, d\beta$. This differs from $\Omega_1$ by $d\chi_\delta \wedge \beta$. Using the assumption that $\iota^*\beta=0$, one sees that as $\delta\to 0$, $d\chi_\delta \wedge \beta\to 0$ uniformly over the $2\delta$-tube. Hence, decreasing $\delta$ if necessary, we find that $J$ is tamed by $\Omega_1$.

We next modify the term $dd^c \phi$. For this, note that there exists a smooth convex function $f\colon \R \to \R$ such that $f(x)=0$ for $x\in[-1,1]$ and $f(x) = x^2 +c$ for $|x|>2$, such that, on $\C^r$ with its standard norm, $ - dd^c f( \|z\| )$ is a non-negative form. This form equals zero on $D^{2r}(0;1)$ and $4\omega_{\C^r}$ outside $D^{2r}(0;2)$. Let $h' = g_{\C}$ be the hermitian metric extending $g$, and introduce also a third hermitian metric $h''= h - \kappa h'$, where $\kappa>0$ is small enough that $h''$ is positive-definite.

Define functions $f_\epsilon=  \epsilon f\circ (\epsilon^{-1/2} h'')$ (here $h''$ abbreviates $z\mapsto h''(z,z)$). Then $ dd^c ( - \kappa h'  - f_\epsilon) $ is non-negative on the fibres, gauge-invariant near $\im(\iota)$, and equals $dd^c \phi$ outside the disc-bundle of radius $2\epsilon$. We shall prove that, for $\epsilon \ll \delta$ and $\kappa \ll 1$, the form
\[ \Omega_2 = \Omega_1 - dd^c\phi +  dd^c (-\kappa h' - f_\epsilon)  \]
(which differs from $\Omega_1$ inside the $2\epsilon$--tube only) still tames $J$.  Setting $\Omega' = \Omega_2$, we will then have a form with the required properties.

Consider tangent vectors $u$, of length 1 with respect to some metric (it doesn't matter which), attached to points inside the $2\epsilon$--tube.
We have $\Omega_2(u,Ju) = p^*\iota^*\Omega (u,Ju)+  dd^c ( - \kappa h'  - f_\epsilon)  (u,Ju) $.  This is certainly positive when $u$ is tangent to the fibre, so let us fix a horizontal distribution $H$ and assume that $u\in H$. The first term, $p^*\iota^*\Omega (u,Ju)$, is still positive when $\epsilon$ is less than some $\epsilon_0$, and bounded below by a constant $C>0$, independent of $\epsilon < \epsilon_0$.  On the other hand,  $|(dd^c h' )(u,Ju)| \leq C' $ for a similar constant $C'>0$, so taking $\kappa \leq  C/ 2C'$, we have  
\[p^*\iota^*\Omega (u,Ju) +  dd^c (-\kappa h') (u,Ju) \geq \frac{C}{2}>0. \]
The troublesome term in $\Omega_2(u,Ju)$ is the last one, $(dd^c  f_\epsilon)(u,Ju)$. We have
\[ d^c f_\epsilon =   q' (h'')  \,  d^c h'',  \]
where $q(t) =  \epsilon \, f(\epsilon^{-1/2} t^{1/2})$;  hence
\[ dd^c f_\epsilon =  q'' (h'') \, d h'' \wedge d^c h'' +  q' (h'') \, dd^c h''.  \]
Choose an $h''$--orthonormal frame $(e_i)$ for the real vector bundle underlying $F\otimes \C$ over a patch $U\subset \bar{M}$. This trivialises the bundle, identifying it with $U\times \R^{2r}$, and we may take $H$ to be the trivial horizontal  distribution. We then have $dh'' (u)=0$, so only the latter term matters.  We write the $U$--coordinates as $(x_j)$. We may as well take $u= \partial_{x_j}$, and we can then write  $Ju = \sum_i{ J_{ij} \partial_{x_i}} + \sum_k {Y_{kj} e_k} $. The terms in the resulting expression for $ dd^c f_\epsilon(u,Ju) $ then involve $Y_{kj}$ or $\partial Y_{kj}/{\partial x_l} $ as coefficients. But both $Y_{kj}$ and its horizontal derivatives go to zero along the zero--section, since the zero--section is almost complex. This implies that $(dd^c h'')(u,Ju)  $, evaluated at $(x,\xi)\in U\times \R^{2r}$, goes to zero as $\xi\to 0$.  Hence, if we take $\epsilon$ small enough, we can ensure that $|(dd^c  f_\epsilon)(u,Ju)| \leq \frac{C}{4}$ inside the $2\epsilon$--tube, whereupon $\Omega_2(u,Ju)\geq \frac{C}{4}>0$.
\end{pf}

\begin{pf}[Proof of \ref{ass symp forms}]
Replace the initial form $\Omega$ by an $\Omega'$ as in Lemma \ref{technical lemma}. (Note that $\Omega'$ is linearly homotopic to $\Omega$, since both are $J$--positive). The new form $\Omega'$ is gauge-invariant along the zero-section.

Let $\bar{\omega} = \iota^*\Omega'$, and let $t>0 $ be the unique constant such that $\Omega'| \Tv_{\iota(x)} N  = t\omega_{\C^r}$ for all $x\in \bar{M}$. Let $q\colon P\times \C^r \to N$ be the quotient map. Choose a connection form $\alpha$ on $P$, and introduce the  two--form $\eta = p^*\bar{\omega} + t(\omega_{\C^r}+ d\langle \mu, \alpha \rangle) \in \Omega^2(P\times \bar{M})$. Besides being $\mathrm{O}(r)$-invariant, $\eta$ tames $J$ in some neighbourhood of the zero-section.

We have $q^*\Omega' - \eta  = d\gamma$ for an invariant one--form $\gamma$ such that $d\gamma$ vanishes along $\im(\iota)$. Introduce cutoff functions $\chi_\delta$ as in the proof of the lemma, and consider the forms $ \eta + d(\chi_\delta \gamma)$. These are also invariant, and so descend to $N$. Choosing $\delta$ small enough that, within the $\delta$-tube, $\eta$ tames $J$, and so that the term $d\chi_\delta\wedge \gamma$ is very small, we find that these forms are also tame $J$. Thus we may set $\Omega_1 = \eta + d(\chi_\delta \gamma)$ and $\Omega_t = t\Omega_1+ (1-t) \Omega' $. \end{pf}

\section{Lagrangian correspondences}

\subsection{Symplectic Morse--Bott fibrations from families of curves}
\label{Lag corr}
In this section we will study some particular symplectic Morse--Bott fibrations and their vanishing cycles.
The key examples will be \emph{relative Hilbert schemes of $n$ points} on families of complex curves.

By considering the vanishing cycles of relative Hilbert schemes, we prove the following:
\begin{CorrThm}
Let $(\Sigma,j)$ be a Riemann surface, and $L\subset \Sigma$ an embedded circle.
Let $\bar{\Sigma}$ be the result of surgery on $L$, that is, the surface obtained by
excising a tubular neighbourhood $\mathrm{nd}(L)$ of $L$ and gluing in two discs.
Let $\bar{j}$ be a complex structure on $\bar{\Sigma}$ which agrees with $j$ outside $\mathrm{nd}(L)$.
Using $j$ and $\bar{j}$, we may consider the symmetric products of $\Sigma$ and of $\bar{\Sigma}$ as complex manifolds.
For any $n>0$, any $(s,t)\in \R_{>0}^2$, and any pair of K\"ahler forms $\omega \in \Omega^2(\sym^n(\Sigma))$
and $\bar{\omega}\in \Omega^2(\sym^{n-1}(\bar{\Sigma}))$ lying in cohomology classes
\[([\omega],[\bar{\omega}])= (s \eta_\Sigma + t\theta_\Sigma, s \eta_{\bar{\Sigma}} + t\theta_{\bar{\Sigma}} ) \]
(see below) there exists a Lagrangian submanifold
\[ \widehat{V}_L\subset \big (\sym^n(\Sigma)\times \sym^{n-1}(\bar{\Sigma}), \omega\oplus -\bar{\omega}\big) \]
such that
\begin{enumerate}
\item [(i)]
the first projection embeds $\widehat{V}_L$ into $\sym^n(\Sigma)$, and
\item [(ii)] the second projection is an $S^1$-bundle over $\sym^{n-1}(\bar{\Sigma})$.
\end{enumerate}
Moreover, we can construct such a $\widehat{V}_L$ canonically up to 
Hamiltonian isotopies through Lagrangians satisfying (i) and (ii). When $n=1$ (so $\sym^{n-1}(\bar{\Sigma})$ is a one-point set), $\widehat{V}_L \subset \Sigma$ is Hamiltonian isotopic to $L$.

In this construction, we can make $\widehat{V}_L$ vary smoothly with the input data $(j, \bar{j}, \omega, \bar{\omega})$. 
\end{CorrThm}
To explain the notation: for any closed Riemann surface $C$, there are two distinguished classes
in $H^2(\sym^n(C);\Z)$, both invariant under the action of the mapping class group of $C$.
They arise via the first Chern class $z$ of the universal divisor
\[ Z^{\mathrm{univ}} = \{ (x, D) : x\in \supp(D) \} \subset C\times \sym^n(C). \]
For $c\in H^*(C;\Z)$, let $c^{[1]} = \mathrm{pr}_{2!}(\mathrm{pr}_1^*c\cup  z)$.
Then $\eta_C = o_C^{[1]}$, where $o_C\in H^2(C;\Z)$ is the orientation class, and
\[\theta_C = \sum_i{ \alpha_i^{[1]}\cup \beta_i^{[1]}},  \]
where $\{ \alpha_i,\beta_i\}$ is a symplectic basis of $H^1(C;\Z)$.

The theorem is a consequence of (i) the general vanishing-cycle construction of the previous section,
together with (ii) an observation about the structure of the singular locus of the Hilbert scheme of $n$
points on a nodal curve.

Our chief concern will be elucidating the geometry of $\widehat{V}_L$. Its interpretation in terms
of points on $\Sigma$ and $\bar{\Sigma}$ is rather subtle---these are not `tautological' correspondences.

\begin{Defn}
By an {\bf elementary Lefschetz fibration} we will mean a triple $(E,\pi,J)$ consisting of a smooth
four--manifold (with boundary) $E$, a proper map $\pi\colon E\to \Delta$ to the closed unit disc $\Delta$,
and an almost complex structure $J$ on $E$ such that $D\pi \circ J = \ii D\pi$. We require that
\begin{enumerate}
\item[(i)]
There is precisely one point $c$ where $\pi$ fails to be submersive; it lies over $0\in \Delta$;
\item[(ii)]
$J$ is integrable in a neighbourhood of $c$;
\item[(iii)]
the complex Hessian form $\frac{1}{2}D^2_{c}\pi \colon T E\otimes_{\C} T E \to \C$ is non-degenerate.
By $D^2_c$ we mean the derivative viewed as a complex bilinear form, and non-degeneracy is over $\C$.
\end{enumerate}
A {\bf holomorphic elementary Lefschetz fibration} is one in which $J$ is globally integrable.
\end{Defn}

We have the following straightforward observation.
\begin{Lem}\label{elem lef fib}
(a) Given a pair $(\Sigma,L)$, an oriented surface with an embedded circle, there is a holomorphic elementary Lefschetz fibration $(E,\pi,J)$ and a diffeomorphism $\delta\colon \pi^{-1}(1) \to \Sigma$ such that $\delta^{-1}(L)$ is, topologically, a vanishing cycle for $E$. Moreover, the construction is canonical in the sense that to specify $(E,\pi,J; \delta)$ we have only to choose $(\Sigma,L)$ and a point in a contractible space.

(b) Let $(E,\pi,J; \delta)$ be as in (a). Let $\widetilde{E}_0 \to E_0$ be the normalisation of the singular fibre $E_0$.\footnote{Readers shy of algebraic geometry may wish to be reminded that the normalisation of a nodal curve is the  resolution of singularities in which one excises a neighbourhood $\{(a,b)\in \Delta\times \Delta : ab =0 \}$ of the node, and replaces it by two discs; and that this is an intrinsic (i.e. coordinate-independent) operation.} Then there is a canonical isotopy-class of diffeomorphisms $\bar{\delta} \colon \widetilde{E}_0\to \bar{\Sigma}$, where $\bar{\Sigma}$ is the surface obtained from $\Sigma$ by surgery on $L$.
\end{Lem}
\begin{pf}
(a) Let $q\colon \C^2\to \C$ be the map $(a,b)\mapsto a^2+b^2$, and let
\[ U =\{ x\in \C^2 :  |q(x)|\leq 1,  \| x \|^4 -  |q(x)|^2 \leq c \}.  \]
Let $Z=\{ (a,b) \in U: (a,b) \in  \sqrt{a^2 +b^2} \, \R^2 \}$. Then the map $ q\colon U\setminus Z \to \Delta$ is a holomorphic submersion, and its fibres are biholomorphic to $\Delta^* \amalg \Delta^*$; it is therefore a trivial holomorphic fibre bundle.

To build $(E,\pi,J)$ one glues the standard piece $U$ to a trivial fibration $\Delta \times (\Sigma\setminus A )\to \Delta$, where $A$ is an annular neighbourhood of $L$, via a holomorphic trivialisation of $U\setminus Z$.

The data needed to set up the gluing are (i) a positively oriented conformal structure $j$ on $\Sigma$; and (ii) the germ of a lift of $j$ to a Riemannian metric $g$ in a neighbourhood of $L$. These clearly form a contractible space. The metric $g$ determines a tubular neighbourhood $L\times [-\epsilon,\epsilon ]\hookrightarrow \Sigma$ of $L$. Its image is biholomorphic to a closed annulus, uniquely up to rotations of the annulus. But an embedding of the annulus (up to rotation) is exactly what is needed to set up the gluing.

(b) is an obvious consequence of the gluing construction.
\end{pf}

We now apply a moduli functor to the family $(E,\pi,J)$ to obtain a new family. Two examples to keep in mind are:
\begin{enumerate}
\item
The \emph{Hilbert scheme of $n$ points on $E$ relative to $S$}. Its fibre over $s\neq 0$ is $\sym^n(E_s)$.
\item
The \emph{Picard fibration of degree $n$}. Its fibre over $s\neq 0$ is the Picard variety $\Pic^n(E_s)$ parametrising holomorphic line bundles of degree $n$. The zero-fibre is the compactified Picard variety of $E_0$, parametrising torsion-free sheaves of rank 1 and degree $n$. (Here we require $\pi$ to be proper.)
\end{enumerate}
The latter two examples are usually constructed as GIT quotients of certain Hilbert schemes.

In these examples, the total space is \emph{non-singular}. These are smooth, projective varieties relative to the base.
The critical fibre has normal crossing singularities, and the structure of the normal crossing divisor can be related to
moduli spaces of objects on the normalisation $\widetilde{E}_0$. In the first two examples, the normal crossing divisor is
itself smooth, so the relative moduli space has a structure of symplectic Morse--Bott fibration.

There are also relative moduli spaces when the central fibre has more than one node, but these will typically not be globally smooth.

It is the Hilbert scheme example which will be developed in detail here, since this is the one which will lead to
Seiberg--Witten-like invariants for broken fibrations. The Picard fibration is briefly considered.

When $\pi$ is proper, one would like to add as a third example a suitably compactified moduli space of stable bundles of rank two and odd degree $d$ on $E$ relative to $S$. However, this is tricky. First, there is a question about whether to fix the determinant; in the algebraic geometry literature one typically does not. There is Gieseker's construction \cite{Gie}, in which  using bundles over semistable models for the nodal curve, in which the special fibre has non-smooth normal crossing singularities. Pandharipande \cite{Pan} takes a different approach, involving torsion--free sheaves (see also Nagaraj--Seshadri \cite{NS}). 

\begin{Rk}
In view of the complications in the algebraic geometry of moduli of stable bundles, it may be simpler to work with their gauge theoretic counterparts. Take  a compact oriented surface $\Sigma$, of genus $g\geq 2$,  with one boundary component. Consider the moduli space $M(\Sigma)$ of flat $\SU(2)$--connections which restrict to a fixed connection $A_0$  on the boundary, with $\hol_{\partial \Sigma}{A_0}=-1$, modulo gauge transformations fixing $A_0$. This is smooth of dimension $6g-6$, and the mapping class group of $(\Sigma,\partial \Sigma)$ acts on it.  

Callahan (unpublished thesis draft) showed that the action of a Dehn twist about a (non--contractible, non--boundary--parallel) separating circle $\gamma$ is a rank $1$ fibred Dehn twist. Seidel (also unpublished) showed that when $\gamma$ is non--separating, the action is by a rank 3 fibred twist. In this case the vanishing cycle corresponds to fixing the conjugacy class of  the holonomy  $\hol_\gamma$, and one can identify the reduced manifold with the corresponding moduli space $M(\bar{\Sigma})$ for the surgered surface. 

Thus, instead of building a Lagrangian correspondence via a symplectic Morse--Bott fibration (as we shall do for symmetric products), one can here reverse the process, first writing down the correspondence and then building a symplectic Morse--Bott fibration from it.
\end{Rk}

\subsection{Lagrangian correspondences via degeneration of symmetric products}

We construct Lagrangian correspondences between symmetric products, first in a special case where we can make things explicit, and then in general. The special case makes it clear that these correspondences \emph{do not} have any simple `tautological' interpretation in terms of points on the surfaces themselves.

\subsubsection{A genus--zero example}\label{genus 0}
There is a simple construction of a Lagrangian correspondence
\begin{equation}
\EuScript{L}_n  \subset \big( \sym^n(S^2) \times \sym^{n-1}(S^2\cup S^2), (-\omega_{\PS^n})\oplus \omega' \big). \end{equation}
Here $\omega_{\PS^n}$ is the Fubini-Study form on $\sym^n(S^2)=\PS^n$, and $\omega' = \bigcup_k{\omega_{\PS^{k-1}}\oplus\omega_{\PS^{n-k}}}$ on $\sym^{n-1}(S^2\cup S^2)=\bigcup_{k=1}^n{\PS^{k-1}\times \PS^{n-k}}$.

Let $a_{n,k}$ designate the diagonal action of $S^1$ on complex projective space $\PS^n $ with weights
\[ (\underbrace{ 1,\dots, 1}_{k} \underbrace {-1,\dots,-1 }_{n-k+1} ),  \quad  k\in \{ 1,\dots, n \}. \]
The action $a_{n,k}$ preserves the Fubini-Study form $\omega_{\PS^n}$, and is moreover a Hamiltonian $S^1$--action,
generated by the moment map $\zeta_{n,k} \colon \PS^n \to \R$ given by
\[  \zeta_{n,k} (z_0\colon\dots : z_n)=\frac{ ( |z_0|^2 + \dots + |z_{k-1} |^2) -( |z_k|^2 +\dots +|z_n|^2) }{2\pi\| z \|^{2} } .  \]
There is a natural identification of the reduced space $\zeta_{n,k}^{-1}(0)/S^1$ with
\begin{equation}\label{FS}
(\PS^{k-1} \times \PS^{n-k},\omega_{\PS^{k-1}}\oplus \omega_{\PS^{n-k}} ),
\end{equation}
To see that the identification is valid symplectically, note that there is an action of
$\mathrm{P}(\U(k) \times \U(n-k+1))\subset \mathrm{P}\U(n+1)$ on $\PS^n$ which commutes with the $S^1$--action.
The symplectic form on the reduced space must be invariant under this group, and is therefore a multiple
$c(\omega_{\PS^{k-1}}\oplus \omega_{\PS^{n-k}})$.

The quotient map $q_k\colon \zeta_{n,k}^{-1}(0) \to\zeta_{n,k}^{-1}(0)/\U(1) $ is identified with the unit circle bundle in the line
bundle $\mathcal{O}(1) \boxtimes \mathcal{O}(-1)$ (the tensor product of line bundles pulled back from the respective factors). The Duistermaat-Heckman formula \cite{DH} implies that, for $t\in (-1,1)$,
the cohomology class of the symplectic form on $\zeta_{n,k}^{-1}(t)/S^1$ varies linearly with $t$, with slope $(-1,1)$,
and from this one can read off that $c=1$.

It follows from the elementary properties of symplectic reduction that the graph
\[ \EuScript{L}_{n,k} =  \{ (x, q_k x) : \zeta_{n,k}(x) = 0  \} \subset ( \PS^n , - \omega_{\PS^n}) \times (\PS^{k-1} \times \PS^{n-k}, \omega_{\PS^{k-1}}\oplus\omega_{ \PS^{n-k}} ) \]
is Lagrangian. Set $\EuScript{L}_n = \bigcup_{k=1}^n{\EuScript{L}_{n,k}}$.

\subsubsection{Vanishing cycle interpretation}\label{VC interpretation}
The Lagrangian correspondences $\EuScript{L}_n$ in the previous example can be interpreted as vanishing cycles for a degeneration of $\PS^n$. We will give this degeneration as an explicit family of projective varieties $\EuScript{H}^{[n]}\to \C$; however, we should explain its geometric origin. We interpret $\PS^n$ as $\hilb^n(\PS^1)$, the Hilbert scheme of $n$ points (parametrising ideal sheaves $\mathcal{I}\subset \mathcal{O}_{\PS^1}$ of finite colength $n$). Given a degeneration of $\PS^1$ to a nodal curve, we may consider its \emph{relative} Hilbert scheme of $n$ points. Specifically, we consider the family $\widehat{\PS^1\times \C}\to \C$ obtained by blowing up $\PS^1\times \C$ at $((1:0), 0)$, and the associated family $\hilb_{\C}(\widehat{\PS^1\times \C})\to \C$. This is a degeneration of $\hilb^n(\PS^1)$. Explicit equations for the relative Hilbert scheme were given by Ran \cite{Ran}; the variety $\EuScript{H}^{[n]}$ is Ran's model for $\hilb^n_{\C}(\widehat{\PS^1\times \C})$.

Define a sequence of complex surfaces $X_n$, $n\geq 2$, as follows.
\begin{itemize}
\item $X_2=\PS^1\times \Delta$, where $\Delta$ is the closed unit disc.
\item For $n\geq 3$, 
 	\[ X_n \subset \PS^1_1 \times\dots \times \PS^1_{n-1}\times \Delta, \]
where $\PS^1_i$ is a copy of $\PS^1$. It is the subvariety cut out by the equations  
 \begin{equation}  
 b_i a_{i+1} = t a_i b_{i+1}, \quad i=1,\dots, n-2 
 \end{equation}
with $(a_i:b_i)$ homogeneous coordinates on $\PS^1_{i}$ and $t$ the coordinate on $\Delta$.
\end{itemize}
The projection onto the last factor, $ X_n \to \Delta$, has non-singular rational fibres over points in $\C^*$. The zero-fibre is, for $n\geq 2$, a chain of $n-1$ rational curves. One can obtain $X_{n+1}$ from $X_n$ by blowing up a point in the last curve in the chain. 

Now, define
\[ \EuScript{H}^{[n]} \subset \PS^n \times \PS^n \times X_n  \]
to be the subspace defined by the equations
\begin{align} 
& a_1 y_n=t b_1 y_0 ; \notag \\
&x_i b_i  y_0 = y_{n-i} a_i x_0 ,\quad i=1,\dots, n-1;   \\
& b_{n-1} x_n = t a_{n-1} x_0, \notag
\end{align}
where $(x_0: \dots : x_n)$ and $(y_0:\dots : y_n)$ are homogeneous coordinates on the two $\PS^n$ factors.

It is routine to verify the
\begin{Lem}
\begin{enumerate}
\item
$\EuScript{H}^{[n]}$ is a complex manifold of dimension $n+1$.
\item
Let $\pi\colon \EuScript{H}^{[n]}\to \Delta$ be the projection map. Then, for $t\in \C^*$, projection onto the first $\PS^n$-factor gives an isomorphism $\pi^{-1}(t) \cong \PS^n$.
\item    
The critical set $\EuScript{H}^{[n]}_{\mathrm{crit}} :=\crit(\pi) \subset \EuScript{H}^{[n]}$ is a complex submanifold of codimension two, contained in $\pi^{-1}(0)$. Moreover, $\EuScript{H}^{[n]}_{\mathrm{crit}}$ is naturally biholomorphic to $\bigcup_{k=1}^n{\PS^{k-1}\times \PS^{n- k}}$. 
\item
Near any point in $\crit(\pi)$, $\pi$ is holomorphically modelled on the map $(z_0,\dots,z_n)\mapsto z_0 z_1$. Hence $\pi^{-1}(0)$ has normal crossing singularities along $\crit(\pi)$.
\end{enumerate}
\end{Lem}
We give $\EuScript{H}^{[n]}$ the K\"ahler structure $\omega$ induced by the standard K\"ahler form on $\PS^n\times \PS^n\times (\PS^1)^{n-1}\times \C$. Then $(\EuScript{H}^n,\pi,\omega)$ is a symplectic Morse--Bott fibration over $\Delta$.
The vanishing cycle associated with the ray $[0,1]$ is a Lagrangian correspondence
\[ \widehat{V} \subset \PS^n  \times \left(\bigcup_{k=1}^n{\PS^{k-1}\times \PS^{n- k}}\right), \]
where the symplectic structure comes from a product of Fubini-Study forms as before. 

{\bf $S^1$--actions.}
Consider the diagonal $S^1$--action on $\PS^n\times \PS^n\times (\PS^1)^{n-1}\times \Delta$ with the following weights:
\begin{itemize}
\item
$  ((+1)^k ; (-1)^{n-k+1} )$ on the first $\PS^n$ (for some $k\in \{1,\dots , n\} $);
\item
$((-1)^k ; (+1)^{n-k+1} ) $ on the second $\PS^n$;
\item
$(-1,1)$ on each $\PS^1_i$;
\item 
$0$ on $\Delta$.
\end{itemize}
Denote this action by $A_{n,k}$. It is Hamiltonian, and its moment map is the sum of the moment maps of the factors. It is easy to check that $A_{n,k}$ leaves $\EuScript{H}^{[n]}$ invariant. The naturality of moment maps implies that the moment map $\mu_{n,k}\colon \EuScript{H}^{[n]}\to \R $ for the action on $\EuScript{H}^{[n]}$ is the restriction of the one on $\PS^n\times \PS^n\times (\PS^1)^{n-1}\times \Delta$.

Let $\rho_t\colon \pi^{-1}(1) \to \pi^{-1}(t)$ denote parallel transport in $\EuScript{H}^{[n]}$ along the ray $[t, 1]$. Since $A_{n,k}$ preserves the fibres of $\pi$, we have
\[\frac{d}{dt} (\mu\circ \rho_t ) =0.  \]
Introduce the `Lefschetz thimble',
\[ W = \{ \rho_t(x): x\in V, \, t\in[0,1] \}  \]
and also
\[  W' =  \bigcup_{k =1}^{n}{\mu_{n,k}^{-1}(0) \cap \pi^{-1}([0,1])}.\]
\begin{Lem} $ W = W'.$
\end{Lem}
\begin{pf}
We have $W\cap \pi^{-1}(0) = \EuScript{H}^{[n]}_{\mathrm{crit}}$. It is also true that $W'\cap \pi^{-1}(0) = \EuScript{H}^{[n]}_{\mathrm{crit}}$, as one can verify directly using the defining equations and the formula for the moment map. If $x\in \EuScript{H}^{[n]} $, with $\pi(x)= t\in (0,1]$, then $x\in W'$ if and only if $\rho_{[s,t]}^{-1} x\in W'$ for all $s\in (0,t]$, where $\rho_{[s,t] }$ is parallel transport over $[s,t]$. Hence $W = W'$. 
\end{pf}

We conclude from the lemma that there is a commutative diagram
\[ \xymatrix{%
& \bigcup_{k=1}^{n}{\zeta_{n,k}^{-1}(0)} \ar^{q_{n,k} }[r] \ar_{\rho_1}[dr]
 	&  \bigcup_{k=1}^{n}{\zeta_{n,k}^{-1}(0) /S^1 }\ar[d]^{\cong} \\
	& & \EuScript{H}^{[n]}_{\mathrm{crit}},
} \]
with $\rho_1$ the limiting parallel transport and $q_{n,k}$ the quotient map. Actually, one be a little more precise: the two spaces on the right can each be identified canonically and holomorphically with $\bigcup{\PS^{k-1}\times\PS^{n-k}}$, and under these identifications the vertical arrow becomes the identity map. Hence
\[ \widehat{V}  = \bigcup_{n=1}^{k} { \EuScript{L}_{n,k} }  \]
where both are considered as subspaces of
$ \PS^n \times \bigcup_{k=1}^{n}(\PS^{k-1}\times \PS^{n-k} ) .$

\subsubsection{The general construction}
{\bf Hilbert scheme of a nodal curve.}
An algebro-geometric interlude is required. 

We work in the category of complex analytic spaces. The Hilbert scheme of $n$ points on a complex curve $C$ is a complex space $\hilb^n(C)$ parametrising those ideal sheaves $\mathcal{I} \subset \mathcal{O}_C$ such that $\sum_{x\in X}{\dim_{\C}(\mathcal{O}_{C,x}/ \mathcal{I}_x)}$ is finite and equal to $n$. Its characteristic property is that there is a coherent ideal sheaf
\[\mathcal{I}_{\mathrm{univ}} \subset \mathcal{O}_{\hilb^n(C)\times C}\] 
on $\hilb^n(C)\times C$ such that the parametrisation sends $z\in \hilb^n(C)$ to $i_z^* \mathcal{I}_{\mathrm{univ}}$, where $i_z \colon C \owns x \mapsto (z,x)\in  \hilb^n(C)\times C $.

The cycle map, that is, the map $ \hilb^n(C)\to \sym^n(C)$ sending an ideal to its support, is an isomorphism if $C$ is non-singular.

One can consider, more generally, a family of curves $X \to S$. The relative Hilbert scheme, $\hilb^n_S(X)\to S$, is an analytic space with a sheaf $\mathcal{I}_{\mathrm{univ}}$ over $X\times_S \hilb^n_S(X)$ which identifies the fibre $\hilb^n_S(X)_s$ with the Hilbert scheme $\hilb^n(X_s)$, for each $s\in S$. There is a cycle map $\hilb^n_S(X)\to \sym^n_S(X)$.

\begin{Ex}
Let $U= \{(z_1,z_2)\in \C^2 : |z_1z_2|\leq 1\} $, and consider the family $\pi\colon U\to \Delta$, $(z_1,z_2)\mapsto z_1 z_2$. According to Ran, its Hilbert scheme $\hilb^n_{\Delta}(U)$ is given as follows. 
We assume $n\geq 2$. Recall the complex surfaces $X_n$ defined in section (\ref{VC interpretation}). The Hilbert scheme is the subspace
\[ \hilb^n_{\Delta}(U) \subset \C^n \times \C^n \times X_n \]
cut out by the equations
\begin{align*}
  	& a_1 y_n=t b_1 ;  \\
	& x_i b_i  = y_{n-i} a_i  ,\quad i=1,\dots, n-1;   \\
	& b_{n-1} x_n = t a_{n-1}, 
\end{align*}
where $(x_1,\dots,x_n;y_1,\dots y_n)$ are the coordinates on $\C^n\times \C^n$. (See \cite{Ran} or \cite{Per} for the details of the universal sheaf.)
\end{Ex}

One can give a direct construction of the Hilbert scheme of a general family of curves $\pi\colon X\to S$ over a curve $S$, where $\pi$ has non-degenerate quadratic critical points, by patching together fibre products of Hilbert schemes of $(z_1,z_2)\mapsto z_1$ and of $(z_1,z_2)\mapsto z_1 z_2$. The most important case is this:

\begin{Prop}\label{smoothness}
Suppose $(E,\pi,J)$ is an elementary Lefschetz fibration. Then the relative Hilbert scheme $\EuScript{E}^{[n]}:= \hilb^n(\pi)$ is a smooth manifold of dimension $2n+2$ equipped with a smooth map $\pi^{[n]} \colon \EuScript{E}^{[n]} \to \bar{D} $. Each fibre has the structure of complex analytic space. The critical locus $\crit(\pi^{[n]})$ is the singular set of $(\pi^{[n]})^{-1}(0)$, hence a complex space, and is itself smooth. 

The normal bundle $N= N_{\crit(\pi)/\EuScript{E}} \to \crit(\pi)$ has a structure of holomorphic vector bundle (induced by $J$), and the complex Hessian form on $N$ is non-degenerate.
\end{Prop}

The statement above is a little complicated because $J$ is not assumed integrable. If $(E,\pi,J)$ is a holomorphic fibration then $\EuScript{E}^{[n]}$ is itself a complex manifold, $\pi^{[n]}$ holomorphic, and for any point $x\in \crit(\pi^{[n]})$ there exist holomorphic charts centred on $x$ and on $\pi^{[n]}(x)$ in which $\pi^{[n]}$ takes the form $ (z_0,\dots, z_n)\mapsto z_0  z_1. $

\begin{Rk}
If $\pi\colon E\to \C$ has two non-degenerate critical points $c_1$ and $c_2$ in the same fibre then $\hilb^n(\pi)$ is singular for any $n\geq 2$. Consider, for example, the unique point in $\hilb^2(\pi)$ lying over $c_1 + c_2 \in \sym^2(E_0)$. This has a neighbourhood which is the fibre product of small neighbourhoods of the $c_i$ in $E$. Thus it is modelled on the singular quadric threefold
\[  \{ (a,b,c,d)\in \C^4: ab=cd \}. \]
\end{Rk}

\subsubsection{Global description of the Hilbert scheme of $n$ points on a nodal curve}\label{global}
\begin{figure}[t!]
\centering
\labellist
\pinlabel $\C^2$ at 200 50
\pinlabel $\C^2$ at 410 50
\pinlabel $\widetilde{\C^2}$ at 310 220
\pinlabel $\mathbb{P}^1$ at 310 100
\endlabellist
\includegraphics[scale=0.6]{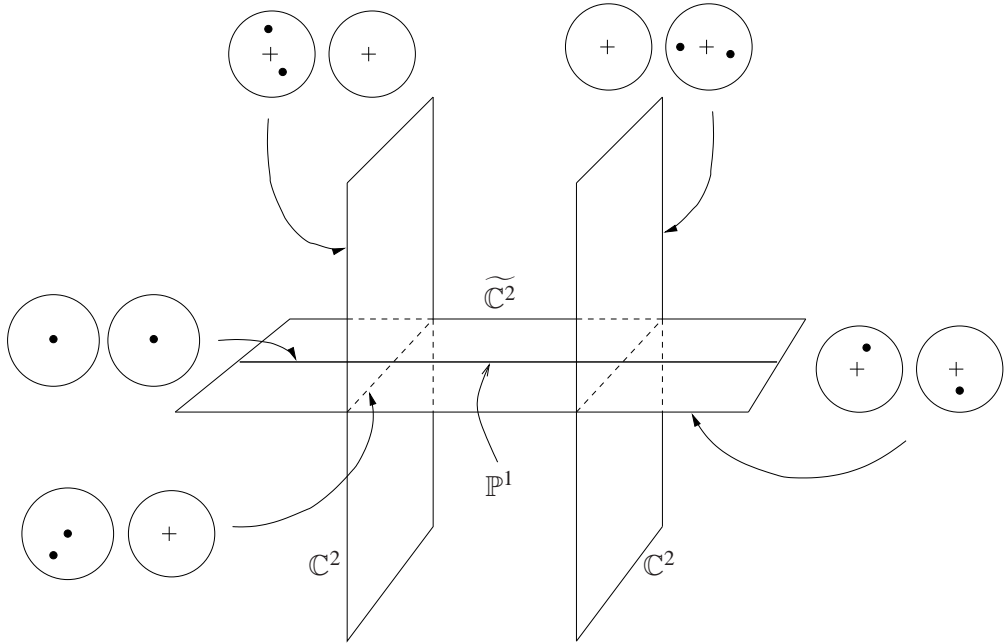}
\caption{\label{sym2}\emph{Schematic picture of $\hilb^2(\{zw=0\})$. It consists of two copies of $\C^2$ and one of $\C^2$ blown up at the origin. The exceptional $\PS^1$ is shown. The normal crossing divisor is isomorphic to $\C\cup\C$ (the normalisation of $\{zw=0\}$). The insertions indicate the behaviour of the cycle map: they show the support of the ideals parametrised by the Hilbert scheme (visualised on the normalisation of $\{zw=0\}$, with its distinguished points marked by crosses).}}
\end{figure}

The singular fibre $\EuScript{E}^{[n]}_0$ of the relative Hilbert scheme of an elementary Lefschetz fibration is the Hilbert scheme  $\hilb^n(C)$ of the nodal curve $C=E_0$. We can give a global description of such a Hilbert scheme as follows. See Figure \ref{sym2} for a picture of $\hilb^2(C)$ when $C = \{ zw =0 \}$.

Let $C$ be a complex curve with precisely one singular point---a node $c$. Thus $c$ has a neighbourhood $U$ isomorphic, as a local-ringed space, to
\[ U_\epsilon:= \{ (z_1,z_2\in \C^2:  z_1 z_2 =0,\, |z_1|^2+|z_2|^2 \leq \epsilon  ) \} \]
for some $\epsilon>0$. Let $\nu\colon \widetilde{C}\to C$ be the normalisation map (so $\widetilde{C}$ is obtained from $C$ by replacing $U_\epsilon$ by two discs of radius $\epsilon$), and let $(c^+,c^-)$ be a labelling of the two points in $\nu^{-1}(c)$. Now consider the embeddings
\[  \iota ^\pm \colon \sym^{n-1}(\widetilde{C})\to \sym^n(\widetilde{C}), \quad D\mapsto D+c^\pm. \]
Let $S= \im(\iota^+)\cap \im(\iota^-)$. Thus $S$ is an embedded copy of $\sym^{n-2}(\widetilde{C})$.

Define $\widetilde{\EuScript{H}}^{[n]}_C$ to be the blow-up of $\sym^n(\widetilde{C})$ along $S$. The embeddings $\iota^\pm$ lift uniquely to embeddings
\[  \widetilde{\iota}^\pm \colon \sym^{n-1}(\widetilde{C})\to \widetilde{\EuScript{H}}^{[n]}_C. \]
Moreover, $\im (\widetilde{\iota}^+)\cap \im(\widetilde{\iota}^- )=\emptyset$.

We now describe an analytic space $\EuScript{H}^{[n]}_C$. As a topological space, it is obtained from $\widetilde{\EuScript{H}}^{[n]}_C$ by gluing $\im(\widetilde{\iota}^+)$ to $\im(\widetilde{\iota}^+)$. That is,
\begin{equation}
\EuScript{H}^{[n]}_C  = \widetilde{\EuScript{H}}^{[n]}_C / \sim \end{equation}
where
\[  \widetilde{\iota}^+ (x) \sim \widetilde{\iota}^-(x) , \quad x\in \sym^{n-1}(\widetilde{C}).  \]
Its complex analytic structure is characterised by the property that the quotient map
\[ \widetilde{\EuScript{H}}^{[n]}_C \to \EuScript{H}^{[n]}_C\]
is holomorphic. This map then becomes the the normalisation map, and $\EuScript{H}^{[n]}_C$ becomes a complex space with normal crossing singularities.
\begin{Lem}
There is an ideal sheaf $\mathcal{I}_{\mathrm{univ}} $ on $ C\times \EuScript{H}^{[n]}_C$ which makes $\EuScript{H}^{[n]}$ into the Hilbert scheme $\hilb^n(C)$.
\end{Lem}

\begin{pf}
{\bf 1.} We begin locally near a node. The following easy algebraic lemma describes ideals in the local ring at a node:

\emph{Consider the $\C$-algebra $\mathcal{O}_{\C^2,0}$ of germs of holomorphic functions near the origin in $\C^2$. Let $R = \mathcal{O}_{\C^2,0}/\langle z_1z_2 \rangle$. Then a proper ideal $I\subset R$ has finite length $\dim_\C{R/I} = n$ if and only if it is of the form $I=I_{m,(a:b)}$,
where}
\begin{equation}\label{ideals}
I_{m,(a:b)} = \langle z_1^{m+1}, z_2^{n-m+1}, az_1^m+b z_2^{n-m} \rangle,\quad 1\leq m \leq n-1,\quad (a:b)\in \PS^1.
\end{equation}
The proof is left to the reader (it is, however, given in \cite[section 4.3.3]{Per}). The formula for $I_{m,(a:b)}$ defines an ideal sheaf over $\mathrm{Spec}(R)\times \bigcup_{m=1}^{n-1}{\PS^1}$ which makes $\bigcup_{m=1}^{n-1}{\PS^1}$ into the `local Hilbert scheme' $\hilb^n(R)$ of ideals in $R$.
 
{\bf 2.} To express this invariantly, think of $R$ as the local ring at a node $x\in C$. Let $C^+$, $C^-$ be the germs of the two sheets of $\widetilde{C}$, and $x^\pm \in C^\pm$ their distinguished points, mapping to the node $x$ by the normalisation map $n\colon \widetilde{C}\to C$. Then $\hilb^n(R)$ is given by pairs $(m,\lambda)$, where $m\in \{1,\dots n-1\}$ specifies a divisor $m [x^+]+ (n-m)[x^-]$ on $\widetilde{C}$, and $\lambda \in \mathbb{P} T_{(x^+,x^-)}(C^+\times C^-)$. 

The local ring $R= \mathcal{O}_{C,x}$ embeds as the subring of $\mathcal{O}_{C^+,x^+}\times \mathcal{O}_{C^-,x^-} $ of functions with $f(x^+)=f(x^-)$. The $m$th copy of $\mathbb{P} T_{(x^+,x^-)}(C^+\times C^-)$ parametrises ideals by the map $(m,\lambda)\mapsto I_{m,\lambda}$, where $ I_{m,\lambda}$ is the ideal of functions $f$ vanishing to order $m$ at $x^-$, to order $n-m$ at $x^+$, and such that the leading order terms of $f|C^+$ and $f|C^-$ are in the ratio $\lambda$.

{\bf 3.} Now we turn to the global structure of ideal sheaves over $C$. The local classification of ideals implies that, to specify an ideal sheaf of colength $n$ over $C$ is precisely to give (i) an effective divisor $D\in \sym^n(C)$; and (ii) if $\{x^+,x^-\} \subset D $, a point $\lambda \in \mathbb{P} T_{(x^+,x^-)}(C^+\times C^-)$. This is simply the statement that there is a bijection from $\EuScript{H}^{[n]}$ to the set of colength $n$ ideal sheaves, $[D,\lambda]\mapsto \mathcal{I}_{D,\lambda}$.

We finally exhibit a universal ideal sheaf $\mathcal{I}_{\mathrm{univ}}$, inducing the bijection $[D,\lambda]\mapsto \mathcal{I}_{D,\lambda}$. It is perhaps clearer to give the corresponding subscheme $\mathcal{Z}_{\mathrm{univ}}$: take the universal divisor
$ \Delta_{\mathrm{univ} }$ in $\widetilde{C}\times\sym^n(\widetilde{C})$, and let $\widetilde{ \Delta}_{\mathrm{univ} }$ be its proper transform in $\widetilde{C}\times\widetilde{\EuScript{H}}^{[n]}_C$. Then $\mathcal{Z}_{\mathrm{univ}}$ is the push-forward of $\widetilde{ \Delta}_{\mathrm{univ} }$ in $C\times \EuScript{H}^{[n]}_C$.

It is enough to check universality of $\mathcal{Z}_{\mathrm{univ}}$ locally on $C$. It is certainly universal near regular points. Moreover, $\mathcal{Z}_{\mathrm{univ}}$, restricted to $C\times \{ [D,\lambda] \} $, is cut out near the singular point $x$ by the ideal $I_{m,\lambda}$ as in (\ref{ideals}). This gives the result.
\end{pf}

The following corollary, an immediate consequence of the identification of $\EuScript{H}^{[n]}_C$ as the Hilbert scheme $\hilb^n(C)$, is crucial for this paper.
\begin{Cor}\label{singular set}
The singular set $\mathrm{sing}(\hilb^n(C))$, i.e. the normal crossing divisor in $\hilb^n(C)$, is naturally identified with $\sym^{n-1} (\widetilde{C}) $. 
\end{Cor}

The normal cone to the divisor $ \mathrm{sing}(\hilb^n(C))$ is the union of two normal line bundles, $N_+$ and $N_-$. These are identified with the normal bundles to $\im(\widetilde{\iota}^+)$ and $\im(\widetilde{\iota}^-)$ in the normalisation $\widetilde{\EuScript{H}}^{[n]}_C$. 

\begin{Prop}\label{normal lines}
Considered as holomorphic line bundles over $\sym^{n-1}(\widetilde{C})$, we have
\begin{align*}
& N_+ = \mathcal{O}(\delta^{n-1}_+ - \delta^{n-1}_-), \\
& N_- = \mathcal{O}(\delta^{n-1}_- - \delta^{n-1}_+).
\end{align*}
Here $\delta^{n-1}_\pm$ is the hypersurface $\{D\in \sym^{n-1}(\widetilde{C}): x^\pm\in \supp(D)\}$.
\end{Prop}
\begin{pf}
The situation is that we have a pair $(Y_1, Y_2)$ of codimension $1$ complex submanifolds in a complex manifold $X$, with transverse intersection $Z=Y_1\cap Y_2$. We want to understand the normal bundles $N_{\widetilde{Y_i} / \widetilde{X } }$ to the proper transforms $\widetilde{Y_i}$ of $Y_i$ in $\widetilde{X}$, the blow-up  of $X$ along $Z$. But the blow-down map $\widetilde{X}\to X$ identifies $\widetilde{Y_i}$ with $Y_i $, and $N_{\widetilde{Y_i} / \widetilde{X } } $ with $N_{Y_i/X}\otimes  \EuScript{O}_{Y_i}(-Z) $.

Applying this to $Y_1=\delta^n_-$, $Y_2= \delta^n_+$ in $X= \sym^n(\widetilde{C})$, we obtain
\begin{align*}
N_{\widetilde{Y_1} / \widetilde{X}}
& \cong N_{\delta_-^n /\sym^n(\widetilde{C})}\otimes \mathcal{O}_{\delta_-^n}(- \delta_+^{n-1})\\
& \cong N_{\delta_-^n /\sym^n(\widetilde{C})}\otimes \mathcal{O}_{\delta_-^n}(- \delta_-^{n-1})\otimes \mathcal{O}_{\delta_-^n}(\delta^{n-1}_-  - \delta^{n-1}_+)\\
& \cong N_{\delta_-^n /\sym^n(\widetilde{C})}\otimes  (\mathcal{O}_X(- \delta_-^n))|{\delta^n_-}\otimes \mathcal{O}_{\delta_-^n}(\delta^{n-1}_-  - \delta^{n-1}_+)\\
&\cong  \mathcal{O}_{\delta_-^n}(\delta^{n-1}_-  - \delta^{n-1}_+),
\end{align*}
where the last isomorphism uses the adjunction formula. Compute $N_{\widetilde{Y_2} / \widetilde{X}}$ in the same fashion.
\end{pf}

\begin{Rk}
(i) The proposition shows that $N_+$ is dual to $N_-$. This is no accident: thinking of $\hilb^n(C)$ as the zero-fibre of a family $\EuScript{E}^{[n]}$ over the disc, the Hessian form on the  normal bundle to $\mathrm{Sing}(\hilb^n(C))$ in $\EuScript{E}^{[n]}$ gives a perfect pairing of $N_+$ with $N_-$.

(ii) When $\widetilde{C}$ is connected, we have $c_1(N_\pm)=0$. However, when $\widetilde{C}$ is disconnected, $N_\pm$ can be topologically non-trivial.
\end{Rk}

\subsubsection{Structure of the vanishing cycles} We now have an explicit (if opaque) local construction of the relative Hilbert scheme $\pi\colon \EuScript{E}^{[n]}\to \Delta$ of a holomorphic elementary Lefschetz fibration $E\to \Delta$, and a global description of its special fibre $\EuScript{E}^{[n]}_0$. With these in hand, we can now set out the following theorem.
\newpage
\begin{Thm}[Structure theorem]\label{structure}
\begin{enumerate}
\item
For any K\"ahler form $\Omega$ on $\EuScript{E}^{[n]}$, $ (\EuScript{E}^{[n]},\pi,\Omega)$ is a symplectic Morse--Bott fibration. 
\item  
There is a natural holomorphic identification of the critical set $\crit(\pi)$ with $\sym^{n-1}(\widetilde{E_0})$.

\item
The vanishing cycle
\[ \widehat{V}  \subset \sym^n(E_1) \times \sym^{n-1}(\widetilde{E_0})  \] 
associated with the ray $[0,1]\subset \Delta$ is Lagrangian with respect to $\omega \oplus -\bar{\omega}$, where $\omega$ (resp. $\bar{\omega}$) is the restriction of $\Omega$ to $ \sym^n(E_1)$ (resp. $ \sym^{n-1}(\widetilde{E}_0) $). 

\item
The projection $ \widehat{V}  \to  \sym^n(E_1) $ is a smooth embedding.
\item
The projection $ \widehat{V} \to \sym^{n-1}(\widetilde{E}_0) $ is an $S^1$-bundle, isomorphic to the unit circle bundle in $ \mathcal{O}(\delta^{n-1}_+ - \delta^{n-1}_-)$. In particular, when $\widetilde{E}_0$ is connected, $\widehat{V}$ is a trivial $S^1$-bundle over $\sym^{n-1}(\widetilde{E}_0)$.
\end{enumerate}
\end{Thm}

\begin{pf}
For the most part, this is a restatement of what we have already established. Parts (1)--(3) have been proved (Proposition \ref{smoothness}, Corollary \ref{singular set}, and Section \ref{vc}), and (4) is obvious.
  
As to (5), this also goes along the lines of the discussion in Section \ref{vc}. The normal bundle $N= N_{\crit(\pi)/\EuScript{E}^{[n]}}$ carries a non-degenerate complex Hessian form $H$, making it an $\mathrm{O}(2,\C)$-bundle. There is a totally real subbundle $N_{\R} = \{ v\in N : H(v,v)\in \R \}$, reducing the structure to $\mathrm{O}(2,\R)$. The projection $ \widehat{V} \to \sym^{n-1}(\widetilde{E}_0) $ is isomorphic to the unit circle bundle in $N_{\R}$. But $N$ is isomorphic to a sum of line bundles $L\oplus L^\vee$, so that the Hessian corresponds to the canonical pairing on $L\oplus L^\vee$.  Projection onto the $L$-summand gives an isomorphism $N_{\R}\cong L$. By Proposition \ref{normal lines}, we have $L\cong \mathcal{O}(\delta^{n-1}_+ - \delta^{n-1}_-) $, so we are done. 
\end{pf}

\subsection{The Picard fibration} Canonically associated with an elementary Lefschetz fibration $E\to \Delta$, with a chosen section $\sigma$, is a proper family 
\[ \Ppic^n_\Delta(E)\to \Delta  \quad (n\in \Z), \] 
the degree $n$ Picard fibration (see \cite{Cap}, for example). A (closed) point of $\Ppic^n_\Delta(E)$ represents a point $s\in \Delta$ and a torsion-free coherent sheaf on $E_s$, of rank one and first Chern class $n \in H^2(E_s;\Z)=\Z$. For explicit descriptions, see \cite{DS,Per}. The key points about its structure are as follows: 
\begin{itemize}
\item $\Ppic^n_\Delta(E)$ is smooth.
\item For $s\neq 0$, the fibre $\Ppic^n(E_s)$ is canonically identified with the Picard variety $\Pic^n(E_s)$.
\item The critical fibre $\Ppic^n(E_0)$ has normal crossing singularities. If $E_0$ is irreducible, the singular locus $\mathrm{Sing}(\Ppic^n(E_0))$ is naturally identified with $\Pic^{n-1}(\widetilde{E}_0)$. If $E_0$ is reducible with one node, and $\widetilde{E}_0$ has components $C_1$ and $C_2$, then  $\mathrm{Sing}(\Ppic^n(E_0))=\bigcup_{i=0}^n{\Pic^i (C_1)\times \Pic^{n-i}(C_2)}$.
\item
There is a holomorphic Abel--Jacobi map $u\colon \hilb^n_\Delta(E) \to \Ppic^n_\Delta(E)$. Indeed, if $\EuScript{L}$ is a rank one torsion-free sheaf on a nodal curve $C$, then every $s\in \mathbb{P}H^0(\EuScript{L}) $ defines a closed subscheme $(s)$ of $C$, and conversely, every colength $n$ subscheme arises this way. If $n>2g-2$, where $g$ is the genus of the regular fibres, then the fibres have constant dimension, and hence this is a projective vector bundle.
\item 
There is a holomorphic line bundle $\Theta \to \Ppic^n_\Delta(E)$ which is ample relative to $S$, and which restricts over regular fibres to the standard theta line-bundle over the Picard torus.
\end{itemize}
The polarisation in the last point gives rise to a K\"ahler structure on $\Ppic^n_\Delta(E)$. With this it becomes a symplectic Morse--Bott fibration over the disc, of rank one. It has a Lagrangian vanishing cycle
\[   \widehat{V}_{\Ppic} \subset   \Pic^n(\Sigma) \times  \Pic^{n-1}(\bar{\Sigma}), \]
where we have used the identifications $\delta\colon E_1 \to \Sigma$ and $\bar{\delta}\colon \widetilde{E}_0 \to \bar{\Sigma}$.

At a topological level, the Picard fibration of an elementary Lefschetz fibration $E$ is easily understood, providing that the vanishing cycle $L\subset E_1$ is non-separating. The regular fibres of
$\Ppic_{\Delta}^n(E)$ are diffeomorphic to $\T^{2g}$, and the family is the fibre product of a trivial family $\T^{2g-2}\times \Delta \to \Delta$ and a genus 1 elementary Lefschetz fibration $E'$. The vanishing cycle for such a fibre product (equipped with a product symplectic form) is
\[  L' \times \mathrm{diag}_{\T^{2g-2}}  \subset \T^2 \times \T^{2g-2}\times \T^{2g-2} =  \T^{2g} \times  \T^{2g-2} ,\]
where $L'$ is the vanishing cycle for $E'$.

\subsection{Lagrangian correspondences between symmetric products: examples}\label{exx}
\subsubsection{Genus $0$}
We studied the case $\Sigma=S^2$ in subsection \ref{genus 0}.
The characterisation of $\widehat{V}$ as an $S^1$-bundle given in point 5 of the structure theorem \ref{structure} is consistent with our explicit picture.

\subsubsection{$\sym^2$ in genus $2$} Let $E\to \Delta$ be a degeneration of a genus 2 curve along a non-separating loop $L$.
As discussed above, the compactified Picard family $\Ppic^n_{\Delta}(E)$ is, topologically, the
fibre product of a trivial $\T^2$-bundle and a genus 1 Lefschetz fibration $E'$ with vanishing cycle $L'$.
The vanishing cycle for $\Ppic^n_{\Delta}(E)$ (with respect to a product symplectic form) is
\[ L' \times \mathrm{diag}_{\T^2}  \subset \T^2 \times \T^2\times \T^2 =  \T^4 \times  \T^2 .  \]
The Abel--Jacobi map expresses $\sym^2(\Sigma)$ as the blow-up at a point of the 4-torus $\Pic^2(\Sigma)$. Likewise, the relative Hilbert scheme $\hilb^2_{\Delta}(E)$ is the blow-up of the compactified Picard family $\Ppic^2_{\Delta}(E)$ along the section $\sigma$ (where $\sigma(s)$  is the point corresponding to $K_{E_s}$). Note that $\sigma(0)$ is a regular point. We can choose a symplectic form on $\Ppic^n_{\Delta}(E)$ so that the Lefschetz thimble in the Picard family is disjoint from $\im(\sigma)$. We obtain a symplectic form on the blow-up in the standard way. We then find that the vanishing cycle $\widehat{V}$ for the relative Hilbert scheme is just the proper transform of that for the Picard family:
\[ \widehat{V} = L' \times \mathrm{diag}_{\T^2} \hookrightarrow  \widetilde{\T^4} \times  \T^2 .  \]

\subsubsection{Projective bundle range, $n>2g-2$} In this numerical range, there is a holomorphic vector bundle $p\colon U \to \Pic^n(\Sigma)$, of rank $n-g+1$, such that the Abel--Jacobi map identifies $\sym^n(\Sigma)$ with $\PS U$. We can choose K\"ahler forms adapted to this projective bundle.\footnote{Caution: our procedure here is \emph{not} quite the same as the one used to form symplectic associated bundles in Section 2.}

Recall (from \cite[p. 209]{MS}, for example) that in a principal $G$-bundle $P\to B$ with connection $\alpha$, one has a vertical cotangent bundle $i_\alpha \colon T^*_{vert} P \hookrightarrow T^* P$, which carries the closed two--form $i_\alpha^* \omega_{can}$. If $(F,\sigma)$ is a Hamiltonian $G$-space then $ P\times_G F$ can be obtained by symplectic reduction of $T^*_{vert}P \times F$, and therefore carries a closed two--form $\tau=\tau(\sigma,\alpha)$ inherited from  $i_\alpha^*\omega_{can} \oplus \sigma$. This form restricts to $\sigma$ on the fibres, and the connection on $P\times_G F\to B$ it induces is the one associated with $\alpha$.

Let us apply this to the principal $\U(n-g+1)$-bundle $P$ of unitary frames for $U$, and to $(F,\sigma)=(\PS^{n-g},\omega_{\mathrm{FS}})$. We obtain a closed two--form $\tau = \tau(\zeta,\alpha)$ on the associated bundle $\PS U = P\times_{\U(n-g+1)}\PS^{n-g}$, which restricts to $\omega_{\mathrm{FS}}$ on the fibres of $\PS U$. Then $ \tau +  c p^*\zeta $ is a symplectic form when $c\gg 0$. (It seems likely that, when $\alpha$ is the Chern connection, this form is of type $(1,1)$, but we have not verified this.)

The relative Hilbert scheme $\hilb^n_\Delta(E)$ is a projective bundle over $\Ppic^n_\Delta(E)$, and consequently we can apply the same procedure to obtain a symplectic form on $\hilb^n_\Delta(E)$.

One then easily checks that symplectic parallel transport in $\hilb^n_\Delta(E)$ is the horizontal lift of symplectic parallel transport in $\Ppic^n_\Delta(E)$ via the connection $\alpha$. It follows that the coisotropic vanishing cycle $V_L \subset \sym^n(\Sigma)$ is the Abel--Jacobi preimage of the vanishing cycle $V_{\Ppic}$ for the Picard fibration, and that the limiting parallel transport map $\rho\colon V_L \to \sym^{n-1}(\bar{\Sigma})$ is the lift of limiting parallel transport from the Picard fibration, $V_{\Ppic}\to (\Ppic^n_\Delta(E))^\crit$.

\subsubsection{Cohomological correspondences}\label{macdonald}
Recall the structure of the cohomology ring of the symmetric products of a compact Riemann surface $C$, as first established by MacDonald \cite{Mac}. If $C$ is connected, we associate with it a sequence of graded rings
\begin{equation}\label{algebra A}
 \mathbb{S}(C,n)  = \bigoplus_{i=0}^n {  \Z[U]/(U^{i+1})\otimes_{\Z} \Lambda^{n-i}H^1(C;\Z) },\quad \deg U= 2.  \end{equation}
(Here $\mathbb{S}(C,0)=\Z$.) One should think of $U$ as the generator of $H^2(C;\Z)$. If $C$ is disconnected, with components $C_1,\dots, C_m$, we put
\[ \mathbb{S}(C,n) =  \bigoplus_{ \stackrel{n_j\geq 0} {n_1+\dots + n_m=n}}{\mathbb{S}(C_1,n_1)\otimes\dots\otimes \mathbb{S}(C_m,n_m) }  \]
There is a ring isomorphism
\begin{equation} a_C\colon  \mathbb{S}(C,n) \to H^*(\sym^n(C);\Z).  \end{equation}
To specify $a_C$, it is enough to do so on monomials $1\otimes \dots \otimes U_i \otimes \dots \otimes 1 $ and $1\otimes \dots \otimes \lambda_i \otimes \dots \otimes 1 $, where $\lambda_i\in H^1(C_i;\Z)$, since such elements generate the ring $\mathbb{S}(C,n)$. The first Chern class $z$ of the universal divisor in $C\times\sym^n(C)$ gives rise to operations
\begin{equation}\label{operations}
  H^*(C;\Z) \to H^{*+2k-2}(\sym^n(C);\Z) ,\quad c\mapsto c^{[k]} := 
 \mathrm{pr}_{2!} \big( z^k \cdot \mathrm{pr}_1^* c  \big). \end{equation}
(We only need the operation $c\mapsto c^{[1]}$ here.) We have 
\[ a_C(1\otimes \dots \otimes U_i \otimes \dots \otimes 1) =o_{C_i}^{[1]},\quad
a_C(1\otimes \dots \otimes \lambda_i \otimes \dots \otimes 1) =\lambda_i^{[1]} ,\]
where $o_{C_i}$ is Poincar\'e dual to the point class in $C_i$.

There are two distinguished classes in $H^2(\sym^n(C);\Z)$ (when $n=1$ they are linearly dependent). The first is the class
\begin{equation} 
\eta_C = a_C( o_C ) = o_C^{[1]},\end{equation}
where $o_C\in H^2(C;\Z)$ is the orientation class. This is dual to the divisor $x+\sym^{n-1}(C)$. The second is 
\begin{equation}   \theta_C = a_C (\sum_i {\alpha_i \wedge \beta_i }), 
\end{equation}
where $\{ \alpha_i,\beta_j \}$ is a symplectic basis for $H^1(C;\Z)$ (so $\sum_i {\alpha_i \wedge \beta_i }$ is the intersection form on $H_1(C;\Z)$). These two classes are invariant under the action of the mapping class group of $C$, and when $C$ is connected they generate the invariant subring.

The vanishing cycle
\[ \widehat{V} \subset \sym^n(E_1)\times \sym^{n-1}(\widetilde{E}_0),   \]
associated with $\hilb^n_{\Delta}(E)$ is Lagrangian with respect to a form of shape $\omega\oplus -\bar{\omega}$. The symplectic forms $\omega$ and $\bar{\omega}$ cannot be chosen arbitrarily: they must be the restrictions of a closed two--form on $\hilb^n_{\Delta}(E)$, and this entails a relation between their cohomology classes. The pair $([\omega],[\bar{\omega}])$ is constrained to lie on a \emph{correspondence} $\mathrm{corr}^2_n\otimes_{\Z}\R$, where
\[  \mathrm{corr}_n^{2} \subset  H^2(\sym^n(E_1);\Z)\times H^2(\sym^{n-1}(\widetilde{E}_0);\Z) \]
is the space of common restrictions of classes in $ H^2(\hilb^n_\Delta(E);\Z)$.
\begin{Prop}\label{coh corr}
Let 
\[C^k = \{ (j_1^*c,n^*j_0^*c ): c\in H^k(E;\Z) \} \subset H^k(E_1;\Z) \times H^k(\widetilde{E}_0;\Z)\]
where $j_s\colon E_s\to E$ is the inclusion of the fibre over $s\in \Delta$, and $n\colon \widetilde{E}_0\to E_0$ the normalisation map. Then $\mathrm{corr}_n^2$ is spanned by the classes 
\begin{itemize}
\item
$(x^{[1]}\cup y^{[1]},\bar{x}^{[1]}\cup \bar{y}^{[1]})$, where $(x, \bar{x})$ and $(y,\bar{y})$ are in $ C^1$;
\item 
$(z^{[1]},\bar{z}^{[1]})$, where $(z,\bar{z})\in C^2$;
\item
if $E_0$ is irreducible, the class $(\theta_{E_1}, \theta_{\widetilde{E}_0})$. 
\end{itemize}
\end{Prop}
In an appendix to this paper, we describe an embedding $\mathbb{S}(E_0,n) \hookrightarrow H^*(\hilb^n_{\Delta}(E);\Z)$, and show that the quotient is isomorphic to $\mathbb{S}(\widetilde{E},n-2)[-2]$. The proof of Proposition \ref{coh corr} is given there, as a corollary of the calculation of $H^2$. 

We now identify a convenient set of K\"ahler classes on $\hilb^n_\Delta(E)$.\footnote{There appears to be a small error concerning this point in the paper of Donaldson and Smith \cite{DK}. In Theorem 3.6, the line bundle $\Lambda^r(\pi)$ is not, as claimed, relatively ample. However, if one twists it by the Abel--Jacobi--pullback of the theta line bundle $\Theta$ over the Picard fibration, the resulting line bundle \emph{is} relatively ample.} Let $\eta_{\Delta}$ be the class which restricts to $\hilb^n(E_0)$ as $\eta=o_{E_0}^{[1]}$, and let $\theta_\Delta$ be the pullback by the Abel--Jacobi map of $c_1(\Theta_\Delta)$, where $\Theta_\Delta\to \Ppic^n_\Delta(E)$ is the theta line-bundle.

\begin{Lem}\label{Kahler classes}
For any $s,t>0 $, the class $s\eta_{\Delta}+t \theta_{\Delta}$ is represented by a K\"ahler form.
\end{Lem}
Again, the proof fits in with the results on the cohomology of the Hilbert scheme, and so is consigned to the appendix.

\subsection{Proof of Theorem A}
The proof is a matter of assembly: all the pieces have already been constructed. Recall that the input data comprises the pair $(\Sigma,L)$, together with complex structures $j$ on $\Sigma$ and $\bar{j}$ on the surgered surface $\Sigma_L$; and K\"ahler forms $\omega\in \Omega^2(\sym^n(\Sigma))$ and $\bar{\omega} \in \Omega^2(\sym^{n-1}(\bar{\Sigma}))$. 

We can form a holomorphic elementary Lefschetz fibration $E\to \Delta$ with $(\Sigma,j)$ as smooth fibre $E_1$ and topological vanishing cycle $L$, as in Lemma \ref{elem lef fib}. Recall that the construction is essentially canonical. Both $\bar{\Sigma}$ and $\widetilde{E}_0$ are built from $\Sigma$ by excising a neighbourhood of $L$ and gluing in a pair of discs. Because of this, we can choose a diffeomorphism $\bar{\delta} \widetilde{E}_0 \cong \bar{\Sigma}$ which is the identity outside the neighbourhood, and such that $\bar{\delta}^*\bar{j}$ coincides with the complex structure of $ \widetilde{E}_0$. Indeed, this requirement pins down $\bar{\delta}$ precisely, because the open disc has no compactly supported holomorphic automorphisms.

Form the relative Hilbert scheme $\hilb^n_\Delta(E)\to \Delta$, and choose a K\"ahler form
\[ \Omega\in \Omega^2( \hilb^n_\Delta(E)),\quad [\Omega] = s\eta_\Delta+ t\theta_\Delta. \] 
Using $\Omega$ we construct the vanishing cycle $\widehat{V}_\Omega\subset \sym^n(\Sigma;j)\times \sym^{n-1}(\bar{\Sigma};\bar{j})$. This is not quite what is wanted, because the restrictions $\omega_1 =  \Omega|\sym^n(\Sigma)$ and $\bar{\omega}_1= \bar{\delta}^{-1*}(\Omega|{\sym^{n-1}(\widetilde{E}_0)})$ are not the given forms $\omega$ and $\bar{\omega}$. We can remedy this using the convexity of the space of K\"ahler forms. Set 
\[ \omega_t = t \omega_1 + (1-t)\omega , \quad \bar{\omega}_t = t\bar{\omega}_1  + (1-t) \bar{\omega}.\] 
Since the $\omega_t$ (resp. $\bar{\omega}_t$) are cohomologous K\"ahler forms, they give rise to symplectic flows 
\begin{align}\label{flow1} 
&\phi_t\colon (\sym^n(\Sigma),\omega_0)\to (\sym^n(\Sigma),\omega_t),\\
\label{flow2}
&\bar{\phi_t}\colon (\sym^{n-1}(\bar{\Sigma}),\bar{\omega}_0 )\to   (\sym^{n-1}(\widetilde{E}_0),\bar{\omega}_t ).
\end{align}
Put
\[ \widehat{V}_L =  (\phi_1^{-1} \times \bar{\phi}_1^{-1}) ( \widehat{V}_\Omega). \]
This is a Lagrangian submanifold of $ \sym^n(\Sigma;j)\times \sym^{n-1}(\bar{\Sigma};\bar{j})$ with respect to $-\omega\oplus \bar{\omega}$. 

Since the space of auxiliary choices involved in the construction is path connected (indeed, contractible), the Lagrangian isotopy class of $\widehat{V}_L$ is canonically determined. We must check that so too is the Hamiltonian isotopy class. To simplify the notation, we will assume that the data are all fixed apart from the K\"ahler form, which varies in a path $\Omega_t$ and prove Hamiltonian isotopy of the resulting vanishing cycles $\widehat{V}_{L,t}$. The remaining ambiguity (i.e. the precise choice of $(E,\pi,J)$) can be handled in just the same way.
 
The $\Omega_t$ are cohomologous forms. By Lemma \ref{globally coisotropic}, they extend to a closed two--form $\Xi$ on $\hilb^n_\Delta(E)\times [0,1]$ with the property that the union of the coisotropic vanishing cycles, 
\[ \EuScript{V}=  \bigcup_t{V_{\Omega_t}}\subset  \sym^n(\Sigma)\times [0,1],  \]
is \emph{globally coisotropic} with respect to $\zeta:=\Xi | (\sym^n(\Sigma)\times [0,1])$. Then $\zeta| \EuScript{V}$ is the pullback of a closed form $\bar{\zeta}\in \Omega^2(\sym^{n-1}(\widetilde{E}_0)\times [0,1])$ (a closed extension of the family of forms $\Omega_t | \sym^{n-1}(\widetilde{E}_0)$). Using the maps (\ref{flow1}, \ref{flow2}), define
\begin{align*} 
&\Phi \colon \sym^n(\Sigma)\times[0,1] \to \sym^n(\Sigma)\times[0,1] , && (x,t)\mapsto (\phi_t(x),t) ,\\
& \bar{\Phi}\colon \sym^{n-1}(\bar{\Sigma}) \times[0,1] \to   \sym^{n-1}(\widetilde{E}_0)\times[0,1],&& (\bar{x},t)\mapsto (\bar{\phi}_t(\bar{x}),t).
\end{align*}
Then $\Phi^{-1}(\EuScript{V})$ is globally coisotropic with respect to $(\Phi^{-1})^* \zeta$. Using the bundle isomorphism $\bar{\Phi}$, we obtain a sub-bundle
\[\bigcup_{t\in [0,1]} \widehat{V}_{L,t}  \subset  (\sym^n(\Sigma)\times[0,1]) \times_{[0,1]} (\sym^{n-1}(\bar{\Sigma})\times [0,1])  \]
which is globally isotropic with respect to $\zeta \oplus -\bar{\Phi}^*\bar{\zeta}$. This shows that $\widehat{V}_{L,0} $ is Hamiltonian--isotopic to $\widehat{V}_{L,1}$, and so concludes the proof.
 
\subsection{A partial model}\label{van cycles}
It is possible to write down explicit hypersurfaces in $\sym^n(\Sigma)$ (for instance as level sets of $S^1$-valued functions) which represent the smooth isotopy class of the vanishing hypersurface $V\subset \sym^n(\Sigma)$. However, proving the validity of such a model would entail a lengthy diversion, so we will content ourselves with something much weaker.

Let $A\subset \Sigma$ be the image of an embedding of an annulus, 
\[\iota\colon \{z\in \C: \epsilon \leq |z| \leq \epsilon^{-1}\} \hookrightarrow \Sigma,\] 
with $L = i (\{ |z|=1\})$ its unit circle. Set 
\[ A' = \iota(\{ z: 2\epsilon < |z| < (2\epsilon)^{-1} \}) \subset A, \]
a smaller, open annulus. Then $\sym^n(\Sigma)$ contains open subsets 
\[U_i = \sym^i(A')\times \sym^{n-i}(\Sigma\setminus A),\quad i=0,\dots, n.\] 
\begin{Lem}\label{top model}
The pair $(V,\rho)$ (hypersurface in $\sym^n(\Sigma)$, $S^1$-bundle projection to $\sym^{n-1}(\bar{\Sigma})$) is isotopic to a pair $(V',\rho')$ such that 
\begin{itemize}
\item  $V' \cap U_0 =\emptyset$; 
\item $ V' \cap U_1 = L \times \sym^{n-1}(\Sigma\setminus A)$; and
\item $\rho' | ( V' \cap U_1)$ is given by the obvious projection to $\sym^{n-1}(\Sigma\setminus A)\subset \sym^{n-1}(\bar{\Sigma})$.
\end{itemize}
\end{Lem}

\begin{pf}
As pointed out below Definition (\ref{good forms}), one does not need \emph{closed} two--forms to construct topological vanishing cycles; `good' two--forms are adequate (for the almost complex required by the definition we take our usual integrable one). Moreover, any two good two--forms have smoothly isotopic vanishing cycles.

Let $E\to \Delta$ be the  Lefschetz fibration associated with $(\Sigma,L)$ (so $E_1=\Sigma$). The space $E$ is the union of a trivial fibration $E' =  \Delta\times (\Sigma\setminus \mathrm{int}(A))$ and $E_0 = E\setminus \mathrm{int}(E')$. Let $E_0'\subset E_0$ be a slightly smaller open set: say $E_0' = E\setminus (\Delta\times (\Sigma\setminus A'))$.

In the relative Hilbert scheme of $E$ there are open subsets
\begin{align*}
 W_0 & =\hilb^n_{\Delta}(E') = \Delta\times \sym^n(\Sigma\setminus \mathrm{int}(A)),\\
 W_1 & = E'_0 \times_{\Delta}\hilb^{n-1}_{\Delta}(E') = E'_0 \times \sym^{n-1}(\Sigma\setminus \mathrm{int}(A)).  
 \end{align*}
Take a closed two--form form on $W_0$ which is the pullback of a K\"ahler form on $\sym^n(\Sigma\setminus \mathrm{int}(A)))$. Take a closed two--form form on $W_1$ which is the product of a K\"ahler form on $\sym^{n-1}(\Sigma\setminus \mathrm{int}(A))$ and a symplectic form on $E'_0$ which has $L$ as its vanishing cycle.

Extend these forms to a good two--form on $\hilb^n_{\Delta}(E)$. Then the parallel transport over the ray $[0,1]$ preserves $W_0$ and $W_1$, and the conclusions of the lemma are immediate. \end{pf}

The remaining sets $U_i$ extend to subsets $W_i$ of the Hilbert scheme in a similar way. Like $W_1$, the higher $W_i$ are fibre products, and we can choose K\"ahler forms on them which respect the fibre product structure. Extending these to globally--defined good two--forms, and using  the path--connectedness of the space of good two--forms, we obtain the following lemma. 
\begin{Lem}\label{good forms remark}
We can isotope $V$ to another hypersurface $V'$ such that, for each $i$, $V'\cap (\sym^i(A) \times \sym^{n-i} (\Sigma\setminus A))  =  V_i'\times \sym^{n-i}(\Sigma\setminus A)$ for hypersurfaces $V'_i \subset \sym^i(A)$. In other words, for an $n$-tuple $\mathbf{x}\in\sym^n(\Sigma)$, membership of $V'$ depends only on those points of $\mathbf{x}$ which lie in $A$. Moreover, the $V'_i$ can be taken to be `universal', i.e. independent of the topology of $\Sigma$. 
\end{Lem}
Now let $V$ be the image of $\widehat{V}$ in $\sym^n(\Sigma)$. We compute the fundamental class $[V]\in H_{2n-1}(\sym^n(E_1) ;\Z)$. Note that there is an isomorphism 
\[\mu\colon  H_1(\Sigma;\Z)\to H^1(\sym^n(\Sigma);\Z),\quad  h \mapsto \mu(h) = c / h\]
where $c\in H^2(\Sigma \times \sym^n(\Sigma);\Z)$ is dual to the universal divisor. (We have $\mu(h) = (\mathrm{PD} (h) )^{[1]}.)$
\begin{Lem}\label{fund class}
Given an orientation of $L$, there is an orientation of $V$ such that $[V]$ is Poincar\'e dual to $\mu([L])$. 
\end{Lem}
\begin{pf}
If $\phi\in \diff(\Sigma)$ acts as the identity in a neighbourhood of $L$ then, by the naturality of the Hilbert scheme construction, $\sym^n(\phi)$ preserves $[V]$. One such diffeomorphism is the Dehn twist along a circle parallel to $L$. Now, $(\sym^n(\phi))^* \cdot \mu(h)= \mu(\phi_*h)$, so we must have $[V]=\mu(h)$ where $h$ is invariant under the stabiliser of $[L]$ in $\Sp H^1(\Sigma;\Z)$. It follows that $h$ is a multiple of $[L]$. If $[L]=0$ we are done. If not, we need to see that it is a unit multiple. We can test this by computing the intersection number of $V$ with the 1-cycle $\Lambda  = \{[x,\dots,x]: x\in L' \}$, where $L' \subset\Sigma$ is a circle which intersects $L$ transversely in a single point. Deform $V$ to $V'$ as in Remark \ref{good forms remark}. The remark implies that $V'$ intersects $\Lambda$ only in $\sym^n(A)$, and that the intersection number is independent of the topology of $\Sigma$. But since the intersection number is $\pm 1$ in our genus 1 and 2 examples, the general case follows. 
\end{pf}

\subsubsection{Iterated vanishing cycles and Heegaard tori} As we have seen, a circle $L\subset \Sigma$ gives rise to a Lagrangian correspondence
\[ \widehat{V}_L = \widehat{V}_L^n \subset \big( \sym^n(\Sigma)\times \sym^{n-1}(s_L \Sigma),  -\omega\oplus \bar{\omega}\big). \]
where $s_L\Sigma$ denotes the result of surgery along $L$. If $L'$ is another circle, disjoint from $L$, then $L'$ is still visible in $s_L \Sigma$, and it too gives rise to a vanishing cycle 
\[ \widehat{V}_{L'}^{n-1} \subset \sym^{n-1}(s_L\Sigma)\times \sym^{n-2}(s_{L'} s_L \Sigma).  \]
We can then compose these correspondences to obtain
\[  \widehat{V}_{L'}^{n-1}\circ\widehat{V}_L^{n} \subset\sym^n(\Sigma)\times \sym^{n-1}(s_{L'} s_L \Sigma). \]
(Recall that, given smooth correspondences $C_{12} \subset M_1\times M_2$ and $C_{23}\subset M_2\times M_3$, their composition $C_{23}\circ C_{12}\subset M_1\times M_3$ is the set of pairs $(x_1,x_3)$ such that there exists $x_2\in M_2$ with $(x_1,x_2)\in C_{12}$ and $(x_2,x_3)\in C_{23}$).
Thus $\widehat{V}_{L'}^{n-1}\circ\widehat{V}_L^{n}$ is an $(S^1\times S^1)$-bundle over $\sym^{n-2}(s_{L'}s_L\Sigma)$. For it to be Lagrangian, our symplectic form on $\sym^{n-1}(s_L\Sigma)\times \sym^{n-2}(s_{L'} s_L \Sigma) $ should be of shape $-\bar{\omega}\oplus\bar{\bar{\omega}}$, where $\bar{\omega}$ is the same form as was used to form the first vanishing cycle. This is possible, because one can extend $\bar{\omega}$ to a K\"ahler form on the relative Hilbert scheme which has $\sym^{n-1}(s_L\Sigma)$ as its smooth fibre.

It is straightforward to prove, using good two--forms, that $\widehat{V}_{L'}^{n-1}\circ\widehat{V}_L^{n}$ is smoothly isotopic to $\widehat{V}_{L}^{n-1}\circ\widehat{V}_{L'}^{n}$.
\begin{Conj}\label{commutativity}
$\widehat{V}_{L'}^{n-1}\circ\widehat{V}_L^{n}$ is Hamiltonian isotopic to  $\widehat{V}_{L}^{n-1}\circ\widehat{V}_{L'}^{n}$.
\end{Conj}

If $L_1,\dots L_n$ are disjoint circles in $\Sigma$, there is an $n$-times iterated correspondence \[T \subset \sym^n(\Sigma). \]
The second factor does not appear here because it is $\sym^0(s_{L_n}\circ\dots \circ s_{L_1}\Sigma)$, i.e. a point.
\begin{Lem}
The iterated vanishing cycle $T\subset \sym^n(\Sigma)$ is smoothly isotopic to $L_1\times \dots\times L_n$.  
\end{Lem}
\begin{pf}
We prove by induction on $n$ that we can find a sequence of good two--forms on successive relative Hilbert schemes so that the iterated topological vanishing cycle $T$ is equal to  $L_1\times \dots L_n $. This is certainly true when $n=1$, so we suppose that $n>1$ and that the iterated topological vanishing cycle
\[T' := \widehat{V}_{L_2}^{n-1} \circ\dots \circ \widehat{V}_{L_n}^{1} \subset \sym^{n-1}(\Sigma_1)\]
is $L_2\times \dots \times L_n$. Let $E_{L_1}$ be the Lefschetz fibration in which $\Sigma$ degenerates to $\Sigma_1$ along $L_1$. Choose a good two--form on $\hilb^n_\Delta(E_{L_1})$ as in Lemma \ref{top model}. The projection to $\sym^n(\Sigma)$ of its vanishing cycle $\widehat{V}_{L_1}^n$ then intersects $T'$ precisely along $L_1\times\dots \times L_n$.
\end{pf}
In particular, suppose that $L_1\dots, L_g \subset \Sigma$ are disjoint and linearly independent in homology, where $g$ is the genus of $\Sigma$. Then $T$ is isotopic to the \emph{Heegaard torus} $L_1\times \dots \times L_g$. These tori are the cornerstone of Heegaard Floer homology. 

The link with Heegaard Floer homology will be developed in a future paper, which will contain a symplectic refinement of this topological observation, and also address the symplectic interpretation of handle-sliding as an operation on the Heegaard tori.

\subsubsection{A homotopical lemma}
The following lemma will help us control boundary bubbling when we consider pseudoholomorphic curves.
\begin{Lem}\label{pi2 surjects}
Assume that $n\geq 2$. Fix a basepoint $x\in \widehat{V}$, and consider the natural homomorphism 
\[  \theta_x\colon  \pi_2 \big(\sym^n(\Sigma)\times\sym^{n-1}(\bar{\Sigma});x \big) \to \pi_2 \big(\sym^n(\Sigma)\times\sym^{n-1}(\bar{\Sigma} ), \widehat{V};x\big).\] 
When $\bar{\Sigma}$ is connected, $\theta_x$ is surjective.  When $\bar{\Sigma}$ is disconnected, $\coker (\theta_x)$ is a non-trivial cyclic group.
\end{Lem}
\begin{pf}
Write $\Sigma_{[n]}=\sym^n(\Sigma)$, $\bar{\Sigma}_{[n-1]}=\sym^{n-1}{\bar{\Sigma}}$. 
By the exact sequence of homotopy groups for the pair $(\Sigma_{[n]}\times \bar{\Sigma}_{[n-1]}, \widehat{V})$, we have
\[ \coker(\theta_x)\cong \ker\big(\pi_1(\widehat{V};x)\to \pi_1(\Sigma_{[n]}\times \bar{\Sigma}_{[n-1]};x)\big).\] 
Suppose that $\bar{\Sigma}$ is connected. Then $\widehat{V}$ is a trivial $S^1$-bundle over $\bar{\Sigma}_{[n-1]}$ by point (5) of Theorem \ref{structure}, and so has fundamental group $\Z\times H_1(\bar{\Sigma})$ when $n>2$, and $\Z\times \pi_1(\bar{\Sigma})$ when $n=2$. Any class $h$ in its kernel must clearly be a multiple of the fibre of the $S^1$--bundle. But the fibre--class is non-trivial in $H_1(\Sigma_{[n]})=\pi_1(\Sigma_{[n]})$ (it is the class defined by the vanishing circle $L\subset \Sigma$) so the result follows. (There is an alternative argument using the Abel--Jacobi map---this was the approach taken in \cite{Per}).) 

Now suppose $\bar{\Sigma}$ is disconnected. As in the connected case, any class $h$ in the image of  $\pi_2 (\Sigma_{[n]}\times \bar{\Sigma}_{[n-1]},\widehat{V} ; x)\to \pi_1(\widehat{V};x)$ must  be a multiple of the fibre $F$ of the $S^1$--bundle which passes through $x$.  Thus $\coker(\theta_x)$ is cyclic. It must be non-trivial because $F_x$ is in the kernel of the Hurewicz map to $\pi_1(\widehat{V})\to H_1(\widehat{V})$, and hence in the kernel of any homomorphism from $\pi_1(\widehat{V})$ to an abelian group.
\end{pf}

\section{Lagrangian matching conditions from broken fibrations}
\label{LMC from BF}
In this section we explain how the Lagrangian correspondences between symmetric products, constructed and analysed in the previous section, can be cast as `Lagrangian matching conditions' associated with broken fibrations.

\emph{Start with an elementary broken fibration $(X_{br},\pi_{br} )$ over an annulus $A=\{ z\in \C : 1/2\leq |z|\leq2\} $, with connected critical set $Z$ mapping diffeomorphically to $\{|z|= 1 \}$. Let $Y = \pi_{br}^{-1}(\{|z|=1/2\})$ and $\bar{Y}= \pi_{br}^{-1}(\{|z|=2\})$. We suppose that the fibre $\Sigma = \pi_{br}^{-1}(1/2)$ is connected of genus $g$, and that $\pi_{br}^{-1}(2)$ is either connected of genus $g-1$, or else disconnected with components of genera $g_1$ and $g-g_1$.}

We supplement the fibration with the following data, which are all determined up to deformation by $(X_{br},\pi_{br} )$:
\begin{itemize}
\item
A near-symplectic form $\omega$, positive on the fibres of $\pi_{br}$ at regular points, and vanishing along $Z$ (see \cite{ADK}). Once this is chosen, $\bar{Y}$ is identified with the mapping torus $\torus(\phi)$ of its symplectic monodromy $\phi\in \aut(\Sigma,\omega|\Sigma)$. Likewise, $\bar{Y} =\torus(\bar{\phi})$ for some $\bar{\phi}\in \aut(\bar{\Sigma},\omega|\bar{\Sigma})$.
\item
A surface $Q\subset \bar{Y}$ (a torus or Klein bottle) which shrinks to $Z$. Recall that we refer to this as the attaching surface of the fibration. One can choose $\omega$ so that $Q$ is isotropic, though this is rarely necessary for our puropses.
\item
A complex structure $j$ on $\Tv Y$, compatible with $\omega$. There is also a complex structure $\bar{j}$ on $\Tv \bar{Y}$, but usually we choose this in a particular way, to be explained presently.
\end{itemize}
We shall write $(Y^{[n]},\pi^{[n]})$ (resp. $\bar{Y}^{[n-1]},\bar{\pi}^{[n-1]}$) for the relative symmetric product $\sym^n_{S^1}(Y)\to S^1$ (resp. $\sym^{n-1}_{S^1}(\bar{Y})\to S^1$). These are to be considered as differentiable families of complex manifolds.

To state the following theorem---a parametric version of Theorem A---we need to specify some cohomology classes on the relative symmetric products. 

There are operations 
\[ H^*(Y;\Z) \to H^{*+2k-2}(Y^{[n]} ;\Z),\quad c\mapsto c^{[k]}  \]
defined using the universal divisor as in (\ref{operations}), and one obvious family of classes to consider are those of form $w^{[1]}$, where $w\in H^2(Y;\Z)$. Another useful cohomology class is the first Chern class of the vertical tangent bundle.
\begin{MatchThm}[Lagrangian matching condition]
Let $\Omega \in \Omega^2(Y^{[n]})$ and $\bar{\Omega} \in \Omega^2(\bar{Y}^{[n-1]})$ be closed, fibrewise-K\"ahler two--forms representing a pair of cohomology classes of form 
\[  [\Omega]  =  w^{[1]} - \lambda c_1(\Tv Y^{[n]}),\quad   [\bar{\Omega}]  =  \bar{w}^{[1]} - \lambda c_1(\Tv \bar{Y}^{[n-1]}), \quad \lambda>0, \]
where $w$ and $\bar{w}$ are common restrictions of a class $W\in H^2(X_{br};\R)$. Then, inside the fibre product 
\[ (Y^{[n]}\times_{S^1}\bar{Y}^{[n-1]},\pi^{[n]} \times_{S^1} \bar{\pi}^{[n-1]} , (-\Omega) \oplus \bar{\Omega}) \]
there is a canonical isotopy class of isotropic sub-bundles $\EuScript{Q}$ with the property that the projection $\EuScript{Q} \to Y^{[n]}$ is an embedding and the projection $\EuScript{Q} \to \bar{Y}^{[n-1]}$ is an $S^1$-bundle.
\end{MatchThm}
(Note that there is an implicit restriction on $\lambda$ here: if it is too large, the classes will cease to admit K\"ahler forms.)  
\begin{Rk}
The zeroth symmetric product of a space is a point, so when $n=1$, $\EuScript{Q}$ lies inside $Y$ itself. In this case $\EuScript{Q}$ is isotopic to $Q$ itself.
\end{Rk}

The situation differs from Theorem A only in that we are now considering not a single surface $\Sigma$, but a family of surfaces parametrised by $S^1$. The proof runs as follows: (1) we construct an $S^1$-family of elementary Lefschetz fibrations $E_t$; (2) we form the associated $S^1$-family of relative Hilbert schemes, and endow it with a global closed two--form; (3) we form the $S^1$-family of vanishing cycles, arising from symplectic parallel transport into the critical locus, along an $S^1$-family of rays; (4) we fine-tune the global two--form so that the family of vanishing cycles becomes globally coisotropic. The only other point to take care of is that we have specified appropriate cohomology classes. By MacDonald's formula 14.5 in \cite{Mac},
\[c_1( T\sym^n(\Sigma))  =  (n+1-g)\eta_{\Sigma} - \theta_{\Sigma}.  \]
\begin{pf}
\emph{An $S^1$-family of Lefschetz fibrations.} Given a Riemann surface $\Sigma$ with an embedded circle $L$, we can construct an elementary Lefschetz fibration $(E,\pi, J)$ over the closed unit disc, with smooth fibre $E_1 = \Sigma$ and vanishing cycle $L$. It is canonical up to deformation. If we are additionally supplied with $\phi\in \diff(\Sigma)$ leaving $L$ invariant, then we can extend $\phi$ to $\Phi\in \diff(E)$, again canonically up to deformation. This is most clearly seen by directly building the mapping torus of $\Phi$ as a five-dimensional manifold  $\mathbf{E}$ fibred over $S^1$. This is simply a parametric version of Lemma \ref{elem lef fib}. Each fibre $\Sigma_t$ of the mapping torus $\torus(\phi)$ contains a circle $Q_t$, and we may build a Lefschetz fibration $(E_t,\pi_t,J_t)$ with smooth fibre $\Sigma_t= \pi_t^{-1}(1)$ and topological vanishing cycle $Q_t$. Because the construction varies smoothly with parameters, these fit together to form a manifold $\mathbf{E}$ equipped with a smooth map $\mathbf{E}\to S^1\times \Delta$. The composite $\mathbf{E}\to S^1\times \Delta \to S^1$ has compact fibres, and so $\mathbf{E}$ may be identified with the mapping torus $\torus(\Phi)$ of a self-diffeomorphism $\Phi$ of $E$, extending $\phi$.

\emph{The $S^1$-family of relative Hilbert schemes.} We now form the set
\[\EuScript{E}= \bigcup_{t\in S^1}{\hilb^n_{\Delta} (E_t)} . \]
This carries a natural topology and differentiable structure, extending the standard ones on the $\hilb^n_\Delta(E_t)$: it is characterised by the map $\EuScript{E}\to S^1$ being smooth.

The normalisations of the critical fibres of the $E_t$ fit together to form a three--manifold fibred over $S^1$, with a complex structure on its vertical tangent bundle. We may identify this with $\bar{Y}$. This identification induces a map $j_{crit}\colon \bar{Y} \to \mathbf{E}$. We have
\[ \sym^{n-1}_{S^1} (\bar{Y})  =  \bigcup_{t\in S^1}{\mathrm{Sing} ( \hilb^n(\pi_t^{-1}(0) )   } = \EuScript{E}^{\crit}.\] 

\emph{Cohomology of $\EuScript{E}$.} We need to verify that the pairs of classes $(w^{[1]},\bar{w}^{[1]})$ and $(c_1(\Tv Y^{[n]}) , c_1(\Tv \bar{Y}^{[n-1]}))$ arise as common restrictions of classes in $\EuScript{E}$. In the case of the Chern classes, this is almost immediate: the vertical tangent bundle  $\Tv \EuScript{E}$ (vertical with respect to the projection $\EuScript{E}\to S^1$) comes equipped with a complex structure, and $c_1(\Tv \EuScript{E})$ restricts to $Y^{[n]}$ as $c_1(\Tv Y^{[n]})$. It restricts to $\EuScript{E}^\crit$ as $ \Tv \bar{Y}^{[n-1]} \oplus  N^+ \oplus N^-$ for a pair of complex line bundles $N^+$ and $N^-$. However, $N^+\cong (N^-)^\vee$ because of the non-degeneracy of the Hessian pairing, hence $c_1( \Tv \bar{Y}^{[n-1]} \oplus  N^+ \oplus N^-)= c_1(  \Tv \bar{Y}^{[n-1]} )$.

Now consider a pair $(w^{[1]},\bar{w}^{[1]})$, where $(w,\bar{w})$ is the pair of restrictions of $W\in H^2(X_{br})$). There is a unique class $W' \in H^2(\mathbf{E})$ such that $j_{crit}^* W' = w | \bar{Y}$ and $W'| Y = w|Y$. To see this, observe that we can certainly find a family of classes $W'_t \in H^2( E_t) $, $t\in S^1$, which agree with the restrictions of $w$ to $Y_t$ and $\bar{Y}_t$. The possible ways of extending the family $\{W_t'\}$ to a single class $W'$ are parametrised by a quotient of $H^1(E_t)$, so there will precisely one which extends $W$. There is a natural operation $ H^2(\mathbf{E})\to H^{2}(\EuScript{E})$, $ c\mapsto c^{[1]}$, and applying it to $W'$ we get the desired result.

\emph{Adjusting the two--forms.} We can now certainly find a global closed two--form $\mathbf{\Omega}$ on $\EuScript{E}$ which restricts to the classes $[\Omega]$ and $[\bar{\Omega}]$. Moreover, we can arrange (cf. Lemma \ref{Kahler classes}) that is K\"ahler on each $\hilb^n_\Delta(E_t)$. We must perturb it to ensure that the family of coisotropic vanishing cycles $\EuScript{V} = \bigcup_{t\in S^1}{V_{Q_t}}\subset Y^{[n]}$ is globally coisotropic. As we have seen in Section 2 (Lemma \ref{globally coisotropic}), and again in Section 3 (proof of Theorem A, this is achievable.

The result now follows.
\end{pf}

\subsection{Topology of the Lagrangian matching condition}
\label{Lag topology}
It is useful to know a little about the topology of $\EuScript{Q}$. We have already studied the topology of the vanishing cycle $\widehat{V}$---the fibre of the bundle map $\EuScript{Q}\to S^1$---so what we are interested in here is the behaviour of the $S^1$-family.

As the result of surgery, the surface $\bar{\Sigma}$ contains two distinguished points $\{x_+,x_-\}$. In $\bar{Y}$, the pairs of distinguished points on the fibres trace out a 1-manifold $\Gamma$, well-defined up to isotopy---a two-fold cover of $S^1$. If $Q$ is orientable, the covering is trivial, and $\Gamma$ is the disjoint union of two sections of $Y\to S^1$ which we write as $\Gamma^+$ and $\Gamma^-$. If $Q$ is non-orientable, $\Gamma$ is non-trivial.

\begin{Rk}\label{normal form for VC}
As in Section (\ref{van cycles}), $\EuScript{Q}$ is isotopic, through Lagrangian-subbundles, to another such bundle $\EuScript{Q}'$ which has the following property: let $\nu$ be a tubular neighbourhood of $Q$
in $Y$. Consider the open set $\EuScript{U}\subset Y^{[n]}$ of $n$-tuples $[x_1,\dots,x_n]\in \Sigma_t$ such that $x_1$ lies in $\nu$ and $x_2,\dots x_n$ lie outside $\nu$. Then, inside $\EuScript{U}$, $\EuScript{Q}'$ is cut out as $Q\times_{S^1}\sym^{n-1}_{S^1}(Y\setminus \nu)$. Moreover, the bundle map to $\bar{Y}^{[n-1]}$ is the projection to $\sym^{n-1}_{S^1}(Y\setminus \nu)$. 

This can be rephrased as follows. Let $\bar{\nu}$ be a tubular neighbourhood of $\Gamma$ in $\bar{Y}$. Then the $S^1$-bundles $Q\times_{S^1}\bar{Y}^{[n-1]} \to\bar{Y}^{[n-1]} $ and $\EuScript{Q}' \to \bar{Y}^{[n-1]}$ become isomorphic on restriction to $\sym^{n-1}_{S^1}(\bar{Y}\setminus \bar{\nu})\subset\bar{Y}^{[n-1]}$.
\end{Rk}

For the following lemma, we should clarify our notation. We have submersions
\[  \EuScript{Q}\stackrel{\rho}{\longrightarrow} \bar{Y}^{[n-1]}\stackrel{\pi^{[n-1]}} {\longrightarrow} S^1 . \]
The map $\rho $ is an $S^1$-bundle, and its structure group is naturally reduced to $\mathrm{O}(2)$-bundle. We write $\Tv \EuScript{Q}$ to mean $\ker D(\pi^{[n-1]}\circ \rho) \subset T\EuScript{Q}$.

\begin{Lem}\label{Q as line bundle}
Consider $\EuScript{Q}\to\bar{Y}^{[n-1]}$ as an $\mathrm{O}(2)$-bundle. 
\begin{enumerate}
\item
When $Q$ is orientable (a torus), we have $w_1(\EuScript{Q} )=0$. 
\item
In general, 
$w_1(\EuScript{Q} )$ restricts trivially to $H^1(\sym^{n-1}(\bar{\Sigma});\Z/2)$.
\item
Let $\delta_\Gamma\subset\EuScript{Q}$ be the codimension-two sub-fibre bundle of $(n-1)$-tuples containing a point of $\Gamma$. Then $[\delta_\Gamma]$ is Poincar\'e dual to $ w_2(\EuScript{Q})$. If $Q\to S^1$ is orientable then $\EuScript{Q}$ reduces to an $\SO(2)$-bundle, and $c_1(\EuScript{Q}) $ is Poincar\'e dual to $[\delta_{\Gamma^+}] - [\delta_{\Gamma^-}]$.
\item
$w_1(\Tv \EuScript{Q})=  \rho^*w_1(\EuScript{Q})$. 
\item
$w_2(\Tv \EuScript{Q})= \rho^*w_2(\EuScript{Q}) + \rho^*w_2(\Tv Y^{[n]})$.
\end{enumerate}
\end{Lem}
\begin{pf}
(1) and (2) are clear. 

(3) The class $ w_2(\EuScript{Q}) - \mathrm{P.D.}[\delta_\Gamma] $ restricts trivially to the fibres, and is therefore equal to $t c$ for some $c\in H^1(\sym^{n-1}(\bar{\Sigma}))$. Since the $\mu$ map on $H_1$ is an isomorphism, we may write $c= \mu(h)$ for some $h\in H_1(\bar{\Sigma};\Z/2)$. But $h$ must be invariant under diffeomorphisms of $\bar{\Sigma}$ trivial near $\{x_+,x_-\}$, hence $h=0$. In the oriented case,  $c_1$ can be handled similarly.

(4)-(5) The vertical tangent bundle $\Tv \EuScript{Q}$ fits into a short exact sequence 
\[ 0 \to \rho^* \EuScript{Q} \to  \Tv \EuScript{Q} \to \rho^* \Tv Y^{[n]}\to 0. \]
But $w_1( \Tv Y^{[n]}) =0$, so the Whitney formula gives the result.
\end{pf}

The last lemma dealt with $\EuScript{Q}$ as a bundle over $\bar{Y}^{[n-1]}$. The next deals with its fundamental class as a submanifold of $Y^{[n]}$.
\begin{Lem} 
Let $\delta_Q$ be the sub-fibre bundle of $Y^{[n]}\to S^1$ of $n$-tuples which contain a point of $Q$. Then $[\EuScript{Q}] = [\delta_Q]\in H_{2n}(Y^{[n]})$. Here we use $\Z$-coefficients when $Q$ is a torus, $\Z/2$-coefficients otherwise.
\end{Lem}
\begin{pf}
By Lemma \ref{fund class}, the difference $[\EuScript{Q}] - [\delta_Q]$ gives zero when intersected with the fibre $\sym^n(\Sigma)$. It is therefore a class supported on $\sym^n(\Sigma)$. But both $[\EuScript{Q}] $ and $[\delta_Q]$ go to zero inside the $S^1$-family of Hilbert schemes $\EuScript{E}$, whereas fundamental classes of components of $\sym^n(\Sigma)$ do not, hence $[\EuScript{Q}] - [\delta_Q]=0$. \end{pf}

\appendix
\section{Cohomology of the relative Hilbert scheme}

In this appendix we describe the additive structure of the integral cohomology of $\hilb^n_\Delta(E)$, where $E\to\Delta$ is a
connected elementary Lefschetz fibration. This is more than is really needed here---it is $H^2$ which is actually used.
One reason to investigate the full cohomology is that there is an interesting homomorphism $H^*(\hilb^n_\Delta(E))\to HF_*(\tau)$
to the Floer homology of the symplectic monodromy of the Hilbert scheme (a fibred Dehn twist). Notice that $\hilb^n_\Delta(E)$ deformation-retracts to its central fibre $\hilb^n(E_0)$, which therefore has the same cohomology; we shall pass freely from one space to the other.

Some notation is required:

(i) We have inclusion maps
\[ \xymatrix{
 \sym^n(E_1) \ar[r]^{i_1} & \hilb^n_{\Delta}(E)  &  \sym^{n-1}(\widetilde{E}_0)\ar[l]_{i_{\crit}}.
  } \]

(ii) We write $j_s\colon E_s\to E$ for the inclusion of the fibre $E_s$ of $E$.
Also write $n\colon \widetilde{E}_0\to E_0$ for the normalisation map, and $j_{\mathrm{crit}} =  n\circ j_0$.

(iii) The relative Hilbert scheme has a universal divisor
\[ Z^{\univ}_\Delta\subset \hilb^n_\Delta(E) \times_{\Delta} E.  \]
This has a dual class $z \in H^2(\hilb^n_\Delta(E)\times_\Delta E )$. There are resulting operations
\[  H^*(E) \to H^{*+2k-2}(\hilb^n_\Delta(E)) ,\quad c\mapsto c^{[k]} :=
 p_{1!} \big( z^k \cdot p_2^* c  \big) , \]
($\Z$ coefficients) where $p_1$ and $p_2$ are the projections from the fibre product onto $ \hilb^n_\Delta(E) $ and $E$ respectively.
(We will only be concerned with the first of these operations, $c\mapsto c^{[1]}$.) There are similar operations
\begin{align*}
& H^*(E_1)\to  H^{*+2k-2}(\sym^n(E_1)) ,\\
& H^*(\widetilde{E}_0)\to H^{*+2k-2}(\sym^{n-1}(\widetilde{E}_0)),
\end{align*}
also denoted $c\mapsto c^{[k]}$ and defined using universal divisors in exactly the same way. These are related as follows:
\begin{Lem}\label{coh lem}
For any $c\in H^*(E)$, we have $i_1^* c^{[k]} = (j_1^*c)^{[k]}$, for any $k\geq 1$.
If $c$ has positive degree then $i_{\crit}^* c^{[1]} = (j_{\mathrm{crit}} ^*c)^{[1]}$.
\end{Lem}

\begin{pf}
The first assertion follows from a comparison of the universal ideal sheaves (divisors). These divisors are
 \begin{align*}
  & Z^{\mathrm{univ}}_{\Delta} \subset E\times_\Delta \hilb^n_\Delta(E),  \\
  & Z^{\mathrm{univ}}_{E_1} \subset E_1 \times \sym^n(E_1).
 \end{align*}
Clearly, $Z^{\mathrm{univ}}_{\Delta}$ pulls back to $Z^{\mathrm{univ}}_{E_1}$, and this immediately gives the result.

The second assertion can also be proved by studying the universal divisor
$Z^{\mathrm{univ}}_{\widetilde{E}_0} \subset \widetilde{E}_0 \times \sym^{n-1}(\widetilde{E}_0)$,
but since this leads to rather indigestible formulae, we shall instead prove it using  Poincar\'e-Lefschetz duality.
Any class $c\in H^k(E)$, $k=1,2$, is Poincar\'e-Lefschetz dual to a cycle $\zeta_c$,
representing a class $[\zeta_c] \in H^{4-k}(E,\partial E)$. We can take $\zeta_c$ to be a smooth cycle which maps
submersively to $\Delta$ (in particular, $\zeta_c\cap E^\crit =\emptyset$). Then $c^{[1]}$ is dual to the cycle $\zeta_c^{[1]}$ of points $\mathrm{x}\in \hilb^n_\Delta(E)$ with $\supp(\mathrm{x})\cap \zeta_c \neq \emptyset$. So $\zeta_c^{[1]}\cap \crit(\hilb^n(C))$ is the set of points $\mathbf{x}\in \sym^{n-1}(\widetilde{E}_0)$ which hit $\zeta_c\cap E_0$. But this means that
$\zeta_c^{[1]}\cap \crit(\hilb^n(C))$ is Poincar\'e dual (in $\sym^{n-1}(\widetilde{E}_0)$) to $(j_{\mathrm{crit}} ^*c)^{[1]}$.
\end{pf}

Recall from subsection \ref{macdonald} the algebras $\mathbb{S}(C,n)$ associated with a curve $C$. Each such algebra comes with a homomorphism $a_C\colon \mathbb{S}(C,n) \to H^*(\sym^n(C))$. When followed by the map on cohomology induced by the cycle map $\hilb^n(C)\to \sym^n(C)$, this gives a homomorphism  $ \mathbb{S}(C,n) \to H^*(\hilb^n(C))$ which we continue to denote by $a_C$.

\begin{Thm}\label{hilb coh} The integral cohomology of the Hilbert scheme of $n$ points on a compact one-nodal curve $C=E_0$
fits into a natural short exact sequence
\[  \xymatrix{
0 \ar[r] & \mathbb{S}(C,n) \ar[r]^-{a_C}  & H^*(\hilb^n(C)) \ar[r]^{b} &\mathbb{S}(\widetilde{C},n-2)[-2]  \ar[r] & 0 .
} \]
\end{Thm}
\begin{pf}
Let $ \widetilde{\EuScript{H}}^{[n]}_C$ be the normalisation of the Hilbert scheme, so there is a
quotient map $\nu\colon \widetilde{\EuScript{H}}^{[n]}_C\to \hilb^n(C)$. We established in Section \ref{global} that $\widetilde{\EuScript{H}}^{[n]}$ is the blow-up of $\sym^n(\widetilde{C})$ along $\sym^{n-2}(\widetilde{C})$, and that $\nu$ folds together two embeddings $\tilde{i}_\pm\colon \sym^{n-1}(\widetilde{C}) \to \widetilde{\EuScript{H}}^{[n]}$.

The map $b$ is given by the composite
\begin{align*}
   H^*(\hilb^n(C)) \stackrel{\nu^*}{\to}  H^* (\widetilde{\EuScript{H}}^{[n]}_C) & \to H^*(Z)\\
   & \to H^{*-2}(\sym^{n-2}(\widetilde{C})) \stackrel{a_C^{-1}}{\to } \mathbb{S}(\widetilde{C},n-2)[-2].
\end{align*}
Here the second map is restriction to the exceptional divisor $Z$ in the blow-up, and the third
map integration down the fibre of the $\PS^1$-bundle $Z\to \sym^{n-2}(\widetilde{C})$.

\emph{Step 1: $b\circ a_C =0$}. It suffices to show that $b\circ a_C(c)=0$ for the generators $c\in \widetilde{H}^*(C)$
of $\mathbb{S}(C,n)$. This is trivially true when $c$ has degree $\leq 1$. As to the degree two case, $\nu^* a_C(o_C)=\nu^* o_C^{[1]} = o_{\widetilde{C}}^{[1]}$ is dual to the divisor $\widetilde{\delta}_x$, the proper transform in $\widetilde{\EuScript{H}}^{[n]}_C$ of $\delta_x\subset \sym^{n}(\widetilde{C})$ (we defined $\delta_x$ to be the locus of $n$--tuples which contain $x$ in their support). Hence $\nu^* a_C(o_C)$ evaluates trivially on the fibres of $Z\to \sym^{n-2}(\widetilde{C})$ and $b\circ a_C(c)=0$.

\emph{Step 2: $a_C$ is injective.} \emph{Irreducible case.} Identifying $H^*(\hilb^n(C))$
with $ H^*(\hilb^n_\Delta(E))$, we have $i_1^*\circ a_C  = a_{E_1} \circ j_1^*$ by Lemma \ref{coh lem}. This shows that $i_1^*\circ a_C$ is injective (here we use the irreducibility of $C$). \emph{Reducible case.} Here $i_{\mathrm{crit}}^*\circ a_C$ is injective, since it equals  $a_{\widetilde{E}_0}\circ j_{\mathrm{crit}}^*$ by Lemma \ref{coh lem}.

\emph{Step 3: $b$ is surjective.} In sheaf cohomology, we have 
\[ H^*(\hilb^n(C); \nu_*\Z) \cong H^*( \widetilde{\EuScript{H}}^{[n]}_C;\Z).\]
Indeed,  $R^q\nu_*\Z=0$ for $q>0$ since the fibres of $\nu$ are zero-dimensional, so the isomorphism is implied by the Leray spectral sequence. There is a short exact sequence of sheaves of abelian groups on $\hilb^n(C)$,
\[0\to \Z_{\hilb^n(C)} \to \nu_*\Z_{\widetilde{\EuScript{H}}^{[n]}_C}\to i_{\mathrm{crit}*} \Z_{\sym^{n-1}(\widetilde{C})}\to 0.\]
(The sheaf on the right is the extension by zero of the constant sheaf $\Z$ on the singular set.) The cohomology exact triangle is
\[\xymatrix{
 H^*(\hilb^n(C)) \ar^{\nu^*}[r]
 	& H^*( \widetilde{\EuScript{H}}^{[n]}) \ar^{(\tilde{i}_+)^*- (\tilde{i}_-)^*}[d] \\
& H^*(\sym^{n-1}(\widetilde{C});\Z) \ar^{+1}[ul]
} \]
By the formula for cohomology of blow-ups,
\[ H^i (\widetilde{\EuScript{H}}^{[n]}_C) \cong H^i(\sym^n(\widetilde{C}))\oplus H^{i-2}(\sym^{n-2}(\widetilde{C})). \]
The composite map $H^i(\hilb^n(C))\stackrel{\nu^*}{\to} H^i( \widetilde{\EuScript{H}}^{[n]}_C)
\to H^{i-2}(\sym^{n-2}(\widetilde{C}))$ (where the second map is projection on a summand) is nothing but $b$.
Our task is therefore to show that $\tilde{i}_+^*- \tilde{i}_-^*$ is zero on this summand.

Now, any class in the summand $H^{i-2}(\sym^{n-2}(\widetilde{C}))\subset  H^i( \widetilde{\EuScript{H}}^{[n]}) $
is supported in the exceptional divisor $Z$. It therefore suffices to show that $(\tilde{i}_+^*- \tilde{i}_-^* )|Z=0$.
But $Z=\mathbb{P}(N^+\oplus N^-)$, a projective bundle over $\sym^{n-2}(\widetilde{C})$, having $\im(\tilde{i}_+)\cap Z$
and $\im(\tilde{i}_-)\cap Z$ as zero- and infinity-sections.
These two sections are homotopic, hence $(\tilde{i}_+^*- \tilde{i}_-^* )|Z=0$, as required.

\emph{Step 4: exactness in the middle.} \emph{Irreducible case.} Since $\widetilde{C}$ is connected, $\tilde{i}_+$
is homotopic to $\tilde{i}_-$, hence $(\tilde{i}_+)^*- (\tilde{i}_-)^*$ is identically zero.
The exact triangle breaks into short exact sequences since $(\tilde{i}_+)^*- (\tilde{i}_-)^*=0$,
and a dimension count finishes the proof.

\emph{Reducible case.} We claim that $(\tilde{i}_+)^*- (\tilde{i}_-)^*$ is surjective.
We restrict it to the summand $H^*(\sym^n(\widetilde{C}))$ in $ H^* (\widetilde{\EuScript{H}}^{[n]}_C)$,
and so consider it as a map
\[\mathbb{S}(\widetilde{C},n) \to \mathbb{S}(\widetilde{C},n-1).\]
Now, $\mathbb{S}(\widetilde{C},n)$ is generated by monomials $ m_{i,n-i}$, where the subscripts designate the summand
$\mathbb{S}(\widetilde{C}_1,i)\otimes \mathbb{S}(\widetilde{C}_2,n-i)$. We have
\[  [\tilde{i}_+^*- \tilde{i}_-^*] m _{i,n-i} =  m _{i-1,n-i} - m_{i,n-i-1}, \]
(monomials with negative subscripts are read as zero.) Granted this formula, surjectivity is an easy exercise.

So again the exact triangle breaks into short exact sequences. Now take $x\in \ker (b)$.
Then $\nu^*(x)\in H^*(\sym^n(\widetilde{C})) = \mathbb{S}(\widetilde{C},n)$, so we may write $\nu^*(x)=a_{\widetilde{C}}(y)$. We have $\nu^*\colon \mathbb{S}(C,n)\cong \mathbb{S}(\widetilde{C},n)$, so $y=\nu^*(y')$ for a unique $y' \in \mathbb{S}(C,n)$. But then $ \nu^* x= a_{\widetilde{C}} \nu^* (y') =\nu^* a_C (y') $, and since $\nu^*$ is injective, $x = a_C(y')$.
\end{pf}
\begin{Cor}\label{hilb H2}
$H^2(\hilb^n(C)) \cong H^2(C)\oplus  \Lambda^2 H^1(C) \oplus  \Z,$
where the third summand is generated by $c_1(\mathcal{O}(Z))$.
\end{Cor}
\begin{Rk}
It would be interesting to understand (along the lines of \cite[ch. 8]{Nak}.) the operations on $\bigoplus_{n\geq 0}{H^*(\hilb^n(C))}$ induced by the correspondences
\[\{(\mathcal{I}_1,\mathcal{I}_2) : \mathcal{I}_1 \subset \mathcal{I}_2, \,\supp(\mathcal{I}_2/\mathcal{I}_1) = \{ \mathrm{node} \} \}\subset \hilb^n(C)\times \hilb^{n+i}(C).\] 
\end{Rk}

We finish by giving two deferred proofs of cohomological results concerning the relative Hilbert scheme. The first, Proposition \ref{coh corr}, gave the structure of the cohomological correspondence $\mathrm{corr}^2_n$.

\begin{pf}[Proof of Prop. \ref{coh corr}]
By Corollary \ref{hilb H2},
\[H^2(\hilb^n(E_0)) \cong H^2(E_0)\oplus  \Lambda^2 H^1(E_0) \oplus  \Z,
\]
where the third summand is generated by $\zeta$, so our task is to understand $i_{\mathrm{crit}}^*$ and $i_1^*$ on each of the three summands. By Lemma \ref{coh lem}, $i_{\mathrm{crit}}^*(c^{[1]})= (n^* j_0^*c)^{[1]}$ and $i_1^*(c^{[1]}) = i_1^* c^{[1]}$.
The class $\theta_{E_1}$ is the pullback by the Abel--Jacobi map of 
\[ c_1(\Theta_{E_1})\in H^2(\Pic^n(E_1)) ,\]
where $\Theta_{E_1}$ is the standard theta line-bundle.
Likewise, $\theta_{\widetilde{E}_0}$ is the pullback of $c_1(\Theta_{\widetilde{E}_0}) \in H^2(\Pic^{n-1}(\widetilde{E}_0);\Z)$.
The theta line-bundles extend to a line bundle $\Theta_\Delta\to \Ppic^n_{\Delta}(E)$
(the natural ample class relative to the base which is part of the construction of the compactified Picard family).
This pulls back to a class on $\hilb^n_\Delta(E)$. Hence $(\theta_{E_1}, \theta_{\widetilde{E}_0})\in \mathrm{corr}_n^2$.

We need to see that the classes we have constructed span $\mathrm{corr}^2_n$.
When $E_0$ is irreducible, we have at our disposal a set $B$ of linearly independent classes in $\mathrm{corr}^2_n$,
such that $B$ bijects with a $\Z$-basis for $\hilb^n(C)$, and such that each member of $B$ is simple
(not a multiple of another class by an integer $>1$). Hence $B$ spans $\mathrm{corr}^2_n$.

When $E_0$ is reducible, the result will follow as soon as we can show that the natural map
$H^2(\hilb^n(E_0))\to \mathrm{corr}^2_n$ has non-trivial kernel.
Consider the class $c=\zeta- o_{E_0}^{[1]}$. One can see directly that  $i_{\mathrm{crit}}^*c = 0$.
Thus $\langle c , S\rangle =0$ for any 2-cycle $S$ contained in $\sym^{n-1}(\widetilde{E}_0)$.
But such homology classes such 2-cycles span $H_2(\sym^n(E_1))$, hence $i_1^*c=0 $ too.
\end{pf}

The other deferred proof concerned K\"ahler classes on the relative Hilbert scheme.

\begin{pf}[Proof of Lemma \ref{Kahler classes}]
There are several ways to go about this. One is to use the fact that, for $n \gg 0$, $\hilb^n_\Delta(E)$ is the total space of a projectivised holomorphic vector bundle $p\colon \mathbb{P}V\to \Ppic^n_\Delta(E)$. Thus $(p^*\Theta)^{\otimes M}\otimes \mathcal{O}_{\mathbb{P}V}(N)$ is an ample line bundle over $ \mathbb{P}V$ (relative to $\Delta)$, for any $M>0$, $N>0$. It is therefore represented by a closed $(1,1)$--form $\Omega$ which is positive on each fibre $\hilb^n(E_s)$, $s\in \Delta$. Adding the pullback from $\Delta$ of $\ii R  dz\wedge d\bar{z}$, $R\gg 0$, one gets a K\"ahler form. Moreover, $[\Omega] = M\eta_{\Delta}+N\theta_{\Delta}$. By convexity of the set of K\"ahler classes, we can allow $M$ and $N$ to be non-integral. This finishes the proof when $n\gg 0$.

We now argue by descending induction. Take a K\"ahler representative $\Omega'$ for $s\eta_{\Delta}+t \theta_{\Delta}$, and a holomorphic section $\sigma\colon \Delta \to E$. Then $\Omega'$ restricts to the complex submanifold $\sigma+ \hilb^{n-1}_\Delta(E)\subset\hilb^{n-1}_\Delta(E)$ as a K\"ahler form $\Omega''$. Considered as a closed form on $\hilb^{n-1}_\Delta(E)$, $[\Omega''] = s\eta_{\Delta}+t \theta_{\Delta}$. The result follows.
\end{pf}

\end{document}